\documentclass{amsart}

\usepackage[leqno]{amsmath}
\usepackage{latexsym,amsthm,amsfonts,amssymb,xypic,colonequals}
\usepackage{color,colortbl}
\usepackage{longtable}
\usepackage[tableposition=above]{caption}
\usepackage{upgreek, hyperref}

\newcommand{\numberseries}{\bfseries}   
\newlength{\thmtopspace}                
\newlength{\thmbotspace}                
\newlength{\thmheadspace}               
\newlength{\thmindent}                  
\renewcommand{\subparagraph}{\vspace*{.25\baselineskip}\noindent}
\newtheoremstyle{bfupright head,upright body}
                {\thmtopspace}{\thmbotspace}
                {\upshape}{\thmindent}{\bfseries}{.}{\thmheadspace}
                {{\numberseries \thmnumber{#2\;}}\thmnote{#3}}

\newtheoremstyle{fixed bf head,slanted body}
                {\thmtopspace}{\thmbotspace}{\slshape}
                {\thmindent}{\bfseries}{.}{\thmheadspace}
                {{\numberseries \thmnumber{#2\;}}\thmname{#1}\thmnote{ (#3)}}

\newtheoremstyle{fixed bf head,upright body}
                {\thmtopspace}{\thmbotspace}{\upshape}
                {\thmindent}{\bfseries}{.}{\thmheadspace}
                {{\numberseries \thmnumber{#2\;}}\thmname{#1}\thmnote{ (#3)}}

\newtheoremstyle{numbered paragraph}
                {\thmtopspace}{\thmbotspace}{\upshape}
                {\thmindent}{\upshape}{}{\thmheadspace}
                {{\numberseries \thmnumber{#2.}}}
\theoremstyle{bfupright head,upright body}
\newtheorem{res}{}[section]
\newtheorem{bfhpg}[res]{}               \newtheorem*{bfhpg*}{}
\theoremstyle{fixed bf head,slanted body}
\newtheorem{thm}[res]{Theorem}          \newtheorem*{thm*}{Theorem}
\newtheorem{prp}[res]{Proposition}      \newtheorem*{prp*}{Proposition}
\newtheorem{cor}[res]{Corollary}        \newtheorem*{cor*}{Corollary}
\newtheorem{lem}[res]{Lemma}            \newtheorem*{lem*}{Lemma}
\theoremstyle{fixed bf head,upright body}
       \newtheorem*{dfn*}{Definition}
\newtheorem{con}[res]{Construction}     \newtheorem*{con*}{Construction}
\newtheorem{rmk}[res]{Remark}           \newtheorem*{rmk*}{Remark}
\newtheorem{exa}[res]{Example}          \newtheorem*{exa*}{Example}
\newtheorem{stp}[res]{Setup}            \newtheorem*{stp*}{Setup}
\theoremstyle{numbered paragraph}
\newtheorem{ipg}[res]{}
\newlength{\thmlistleft}        
\newlength{\thmlistright}       
\newlength{\thmlistpartopsep}   
\newlength{\thmlisttopsep}      
\newlength{\thmlistparsep}      
\newlength{\thmlistitemsep}     
\newcounter{eqc} 
\newenvironment{eqc}{\begin{list}{\upshape (\textit{\roman{eqc}})}%
    {\usecounter{eqc}%
      \setlength{\leftmargin}{\thmlistleft}%
      \setlength{\labelwidth}{\thmlistleft}%
      \setlength{\rightmargin}{\thmlistright}%
      \setlength{\partopsep}{\thmlistpartopsep}%
      \setlength{\topsep}{\thmlisttopsep}%
      \setlength{\parsep}{\thmlistparsep}%
      \setlength{\itemsep}{\thmlistitemsep}}}%
  {\end{list}}%
\newcommand{\eqclbl}[1]{{\upshape(\textit{#1})}}
\newcounter{prt}
\newenvironment{prt}{\begin{list}{\upshape (\alph{prt})}%
    {\usecounter{prt}%
      \setlength{\leftmargin}{\thmlistleft}%
      \setlength{\labelwidth}{\thmlistleft}%
      \setlength{\rightmargin}{\thmlistright}%
      \setlength{\partopsep}{\thmlistpartopsep}%
      \setlength{\topsep}{\thmlisttopsep}%
      \setlength{\parsep}{\thmlistparsep}%
      \setlength{\itemsep}{\thmlistitemsep}}}%
  {\end{list}}%
\newcommand{\prtlbl}[1]{{\upshape(#1)}}
\newcounter{rqm}
\newenvironment{rqm}{\begin{list}{\upshape (\arabic{rqm})}%
    {\usecounter{rqm}%
      \setlength{\leftmargin}{\thmlistleft}%
      \setlength{\labelwidth}{\thmlistleft}%
      \setlength{\rightmargin}{\thmlistright}%
      \setlength{\partopsep}{\thmlistpartopsep}%
      \setlength{\topsep}{\thmlisttopsep}%
      \setlength{\parsep}{\thmlistparsep}%
      \setlength{\itemsep}{\thmlistitemsep}}}%
  {\end{list}}%
\newenvironment{prf*}[1][Proof]{%
  \begin{proof}[\bf #1]
    \setcounter{equation}{0}
    } {\end{proof} }
\newcommand{\proofoftag}[2][:]{(#2)#1}
\newcommand{\proofofimp}[3][:]{\mbox{\eqclbl{#2}$\!\implies\!$\eqclbl{#3}#1}}
\newcommand{\pgref}[1]{\ref{#1}}
\newcommand{\stpref}[2][Setup~]{#1\ref{stp:#2}}
\newcommand{\thmref}[2][Theorem~]{#1\pgref{thm:#2}}
\newcommand{\corref}[2][Corollary~]{#1\pgref{cor:#2}}
\newcommand{\prpref}[2][Proposition~]{#1\pgref{prp:#2}}
\newcommand{\lemref}[2][Lemma~]{#1\pgref{lem:#2}}
\newcommand{\conref}[2][Construction~]{#1\pgref{con:#2}}

\newcommand{\exaref}[2][Example~]{#1\pgref{exa:#2}}
\newcommand{\rmkref}[2][Remark~]{#1\pgref{rmk:#2}}
\newcommand{\secref}[2][Section~]{#1\ref{sec:#2}}
\newcommand{\partpgref}[2]{\ref{#1}\prtlbl{#2}}
\newcommand{\partthmref}[3][Theorem~]{#1\partpgref{thm:#2}{#3}}
\newcommand{\partprpref}[3][Proposition~]{#1\partpgref{prp:#2}{#3}}
\newcommand{\partlemref}[3][Lemma~]{#1\partpgref{lem:#2}{#3}}
\newcommand{\thmcite}[2][?]{\cite[Thm.~#1]{#2}}
\newcommand{\rmkcite}[2][?]{\cite[Rmk.~#1]{#2}}
\newcommand{\corcite}[2][?]{\cite[Cor.~#1]{#2}}
\newcommand{\prpcite}[2][?]{\cite[Prop.~#1]{#2}}
\newcommand{\lemcite}[2][?]{\cite[Lem.~#1]{#2}}
\newcommand{\seccite}[2][?]{\cite[Sect.~#1]{#2}}
\newcommand{\dfncite}[2][?]{\cite[Def.~#1]{#2}}
\renewcommand{\eqref}[1]{(\pgref{eq:#1})}

\definecolor{Gray}{gray}{0.925}
\def\urltilda{\kern -.15em\lower .7ex\hbox{\~{}}\kern .04em} 
\newcommand{\set}[2][\mspace{1mu}]{\{#1 #2 #1\}}
\newcommand{\setof}[3][\mspace{1mu}]{\{#1#2 \mid #3#1\}}
\newcommand{\ZZ}{\mathbb{Z}}
\newcommand{\qtext}[1]{\quad\text{#1}\quad}
\newcommand{\qqtext}[1]{\qquad\text{#1}\qquad}
\newcommand{\qand}{\qtext{and}}
\newcommand{\qqand}{\qqtext{and}}
\newcommand{\deq}{\:=\:}
\newcommand{\dsubseteq}{\:\subseteq\:}
\newcommand{\dge}{\:\ge\:}
\newcommand{\dle}{\:\le\:}
\newcommand{\dis}{\:\is\:}
\newcommand{\gra}{\alpha}
\newcommand{\mfm}{\mathfrak{m}}
\newcommand{\mfq}{\mathfrak{q}}
\newcommand{\is}{\cong}
\newcommand{\onto}{\twoheadrightarrow}
\newcommand{\lla}{\longleftarrow}
\newcommand{\lra}{\longrightarrow}
\newcommand{\xra}[2][]{\xrightarrow[#1]{\;#2\;}}
\newcommand{\poly}[2][k]{#1[#2]}
\newcommand{\pows}[2][k]{#1[\mspace{-2.3mu}[#2]\mspace{-2.3mu}]}
\newcommand{\mapdef}[4][\rightarrow]{\nobreak{#2\colon #3 #1 #4}}
\newcommand{\dmapdef}[4][\lra]{\nobreak{#2\colon #3\:#1\:#4}}
\renewcommand{\H}[2][]{\operatorname{H}_{#1}(#2)}
\newcommand{\Shift}[2][]{\mathsf{\Sigma}^{#1}{#2}}
\newcommand{\Soc}[1]{\operatorname{Soc}{#1}}
\newcommand{\SocR}{\operatorname{Soc}R}
\newcommand{\rnk}[2][k]{\operatorname{rank}_{#1}#2}
\newcommand{\lgt}[2][R]{\operatorname{length}_{#1}#2}
\newcommand{\type}[2][R]{\operatorname{type}_{#1}#2}
\newcommand{\Hom}[3][R]{\operatorname{Hom}_{#1}(#2,#3)}
\newcommand{\Ext}[4][R]{\operatorname{Ext}_{#1}^{#2}(#3,#4)}
\newcommand{\tp}[3][R]{\nobreak{#2\otimes_{#1}#3}}
\newcommand{\Tor}[4][R]{\operatorname{Tor}^{#1}_{#2}(#3,#4)}

\newcommand{\fsp}[1][1]{f\hspace{#1pt}}
\renewcommand{\le}{\leqslant}
\renewcommand{\ge}{\geqslant}
\newcommand{\prm}{$\mspace{.5mu}'\mspace{.5mu}$}
\newcommand{\It}[2]{#1_{\langle #2 \rangle}}
\newcommand{\bidegrees}[1]{\left(\begin{smallmatrix} \scriptstyle
      #1 \end{smallmatrix}\right)}
\renewcommand{\rnk}[2][\k]{\operatorname{rank}_{#1}#2}
\newcommand{\clC}[1]{\mathbf{C}(#1)}
\newcommand{\clT}{\mathbf{T}}
\newcommand{\clB}{\mathbf{B}}
\newcommand{\clG}[1]{\mathbf{G}(#1)}
\newcommand{\clH}[1]{\mathbf{H}(#1)}
\newcommand{\sfe}{\mathsf{e}}
\renewcommand{\k}{\mathsf{k}}
\newcommand{\sfx}{\mathsf{x}}
\newcommand{\sfy}{\mathsf{y}}

\newcommand{\sfu}{\mathsf{u}}

\newcommand{\sfw}{\mathsf{w}}
\newcommand{\sff}{\mathsf{f}} 
\newcommand{\sfg}{\mathsf{g}} 
\newcommand{\sfA}{\mathsf{A}} 
\newcommand{\sfB}{\mathsf{B}}

\newcommand{\sfV}{\mathsf{V}}
\newcommand{\sfU}{\mathsf{U}}
\newcommand{\sfW}{\mathsf{W}} 
\newcommand{\sfP}{\mathsf{P}}
\newcommand{\Kzl}{\mathrm{K}}
\newcommand{\betrow}[5]{\mathtt{#1} & \mathtt{#2} & \mathtt{#3} & \mathtt{#4}& \mathtt{#5}}

\numberwithin{equation}{res}
\numberwithin{table}{res}

\begin{document}

\title[Generic local rings on a spectrum]{Generic local rings on a spectrum\\
  between Golod and Gorenstein}

\author[L.\,W. Christensen]{Lars Winther Christensen}

\address{Texas Tech University, Lubbock, TX 79409, U.S.A.}

\email{lars.w.christensen@ttu.edu}

\urladdr{http://www.math.ttu.edu/\urltilda lchriste}

\author[O. Veliche]{Oana Veliche}

\address{Northeastern University, Boston, MA~02115, U.S.A.}

\email{o.veliche@northeastern.edu}

\urladdr{https://web.northeastern.edu/oveliche}

\thanks{L.W.C.\ was partly supported by NSA grant H98230-11-0214 and
  Simons Foundation collaboration grant 428308.}

\date{2 July 2023}



\dedicatory{To Lucho Avramov on the occasion of his $75^\text{\it th}$
  birthday}

\begin{abstract}
  Artinian quotients $R$ of the local ring $Q = \pows[\k]{x,y,z}$ are
  classified by multiplicative structures on
  $\sfA = \Tor[Q]{*}{R}{\k}$; in particular, $R$ is Gorenstein if and
  only if $\sfA$ is a Poincar\'e duality algebra while $R$ is Golod if
  and only if all products in $\sfA_{\ge 1}$ are trivial.  There is
  empirical evidence that generic quotient rings with small socle
  ranks fall on a spectrum between Golod and Gorenstein in a very
  precise sense: The algebra $\sfA$ breaks up as a direct sum of a
  Poincar\'e duality algebra $\sfP$ and a graded vector space $\sfV$,
  on which $\sfP_{\ge 1}$ acts trivially. That is, $\sfA$ is a trivial
  extension, $\sfA = \sfP\ltimes \sfV$, and the extremes
  $\sfA = (\k \oplus \Shift{\k})\ltimes \sfV$ and $\sfA = \sfP$
  correspond to $R$ being Golod and Gorenstein, respectively.

  We prove that this observed behavior is, indeed, the generic
  behavior for graded quotients $R$ of socle rank $2$, and we show
  that the rank of $\sfP$ is controlled by the difference between the
  order and the degree of the socle polynomial of~$R$.
\end{abstract}

\ \vspace{-2\baselineskip}

\maketitle

\thispagestyle{empty}

\vspace{-\baselineskip}

\tableofcontents

\vspace{-1\baselineskip}


\section*{Introduction}
\label{sec:Introduction}

\noindent
A commutative noetherian local ring is the abstract form of the ring
of germs of regular functions at a point on an algebraic
variety. Accordingly, textbooks order local rings in a hierarchy based
on the character of their singularity with the nonsingular rings, also
called regular rings, being the most exclusive:
\begin{equation*}
  \text{regular} \ \subset\
  \text{hypersurface} \ \subset \ \text{complete intersection} \ \subset \
  \text{Gorenstein} \:.
\end{equation*}
Most local rings, however, fall outside this hierarchy, and the
characteristics of rings within it tell us little about local rings in
general. For example, consider quotients of a regular local ring by
high powers of its maximal ideal.  Such rings do not fall within the
geometric hierarchy, but they do belong to a recognized class, that of
Golod rings, which has minimal overlap with the hierarchy.
Thus one may wonder,
\begin{center}
  \it What does a typical local ring look like?
\end{center}

\subparagraph Any meaningful answer would presumably be a partial one,
subject to restrictions on certain ring invariants. Within such
restrictions there must, at the very least, be a systematic way to
talk about \textsl{all} rings: a classification. We proceed to
identify a viable set of restrictions that does not render the
question trivial.

A fundamental invariant of a local ring is the minimal number of
generators of its maximal ideal, it is called the \emph{embedding
  dimension.} In embedding dimension $1$, every singular local ring is
a hypersurface, in particular a complete intersection, and meets the
criterion for being Golod. In embedding dimension $2$ the two notions
separate decisively: A singular ring is either complete intersection
or Golod. In embedding dimension $3$ the range widens significantly: A
singular local ring can now be Gorenstein without being complete
intersection or it may not belong to any of the classes in the
geometric hierarchy. Crucially, though, there is a classification: It
is based on multiplicative structures in homology, and we discuss it
below.

In the artinian case, another fundamental invariant of a local ring is
the rank of its socle---the annihilator of the maximal ideal---which
is called the \emph{type}. All rings of type $1$ are Gorenstein, but
in the setting of artinian local rings of embedding dimension $3$ and
type $2$ the question above is nontrivial and a terminology is
available to express an answer.  Before we discuss our answer, we make
a further restriction to graded artinian quotients $R$ of the
trivariate power series algebra $Q$ over a field~$\k$. Our answer is
stated in the terminology of a classification of quotient rings $R$ in
terms of graded-commutative algebra structures on $\Tor[Q]{*}{R}{\k}$;
it is due to Weyman~\cite{JWm89} and to Avramov, Kustin, and
Miller~\cite{AKM-88}. We recall the relevant details of the
classification in \secref{parameters}; for now it suffices to say that
it incorporates instances of two classic results of Golod: (1)~Trivial
multiplication on $\Tor[Q]{\ge 1}{R}{\k}$ characterizes Golod rings
\cite{ESG62}; this uses that $Q$ has embedding dimension $3$.
(2)~Poincar\'e duality on $\Tor[Q]{*}{R}{\k}$ characterizes Gorenstein
rings; this follows from Avramov and Golod \cite{LLAESG71} as
$\Tor[Q]{*}{R}{\k}$ is isomorphic to the Koszul homology algebra of
$R$.

Let $R$ be a graded quotient of $Q$ of type $2$ with socle generators
in degrees $s_1 \le s_2$. One says that $R$ is \emph{compressed} if
the length of $R$ is as large as possible given $s_1$ and
$s_2$. Artinian $\k$-algebras of type $2$ can be parametrized in such
a way that there is a nonempty Zariski open subset of the parameter
space whose points correspond to compressed algebras.  In this sense,
assuming that $\k$ is large, a generic quotient ring $R$ is
compressed. This point is discussed in further detail in
\secref{generic}; till then we focus on compressed rings, as that
notion is a more operational than~``generic.''

Assume that $R$ is compressed and that its defining ideal is generated
by forms of degree $2$ and higher. As $R$ has type $2$, this ideal has
an irreducible decomposition $I_1 \cap I_2$, where $I_1$ and $I_2$
define artinian Gorenstein rings. Assume that also the rings $Q/I_1$
and $Q/I_2$ are compressed---also generic artinian Gorenstein rings
are compressed.  Our main result is \thmref{main}; via \rmkref{A} it
describes the $\k$-algebra $\Tor[Q]{*}{R}{\k}$ as a trivial extension
$\sfP\ltimes \sfV$ where $\sfP$ is a Poincar\'e duality algebra and
$\sfV \ne 0$ a graded vector space. The ring $R$ is Golod if and only
if $\sfP$ is \emph{trivial}, i.e.\ $\sfP = \k \oplus \Shift{\k}$. The
key conclusions of the main theorem depend on the parity of $s_2$ and
are summarized in \corref[Corollaries~]{mixed}--\corref[]{even}.  They
can be further condensed as: \subparagraph
\begin{center}
  \begin{minipage}{.815\textwidth}
    \it If $s_2 \ge 5$, then $\Tor[Q]{*}{R}{\k}$ is a trivial
    extension $\sfP \ltimes \sfV$ of a Poincar\'e duality $\k$-algebra
    $\sfP$ by a graded $\k$-vector space $\sfV$. Moreover,
    \begin{prt}
    \item[$-$] if $s_2$ is odd there is a number $N$, depending on
      $s_2$, such that $\sfP$ is nontrivial for $s_1 < N$ and trivial
      for~$s_1 \ge N$.
    \item[$-$] if $s_2$ is even there are numbers $N_1 < N_2$, both
      depending on $s_2$, such that $\sfP$ is nontrivial for
      $s_1 < N_1$ and trivial for $s_1 \ge N_2$.
    \end{prt}
  \end{minipage}
\end{center}

\subparagraph Our pursuit of the main theorem was spurred by data
collected with the \emph{Macaulay2} implementation \cite{M2-LWCOVl} of
the classification algorithm \cite{LWCOVl14a}. When the data first
came in we were intrigued, because Avramov \cite{LLA12} had
conjectured that this kind of nontrivial trivial extensions, i.e.\
$\Tor[Q]{*}{R}{\k} = \sfP\ltimes \sfV$ with nontrivial $\sfP$ and
$\sfV\ne 0$, would not exist at all. Our attempts to prove this
conjecture first led to the discovery of sporadic
counterexamples~\cite{LWCOVl14}, and with Weyman \cite{CVW-19} we
later developed a construction of rings with this kind of Tor-algebras
and $\sfP$ of any size, but we still thought they were rare. This
trajectory shows how our perspective on these rings changed with
accrual of experimental data to eventually make a full $180^\circ$
turn.
\begin{equation*}
  \ast \ \ast \ \ast 
\end{equation*}
A brief synopsis of the paper is in place: Let $I$ be a homogeneous
ideal in the trivariate polynomial algebra $Q$ over a field $\k$ and
assume that $I$ defines an artinian ring of type $2$. An irreducible
decomposition $I = I_1 \cap I_2$ yields two graded Gorenstein rings,
$Q/I_1$ and $Q/I_2$, and associated to this data is a Mayer--Vietoris
sequence
\begin{equation*}
  \tag{$\flat$}
  0 \lra Q/I \lra Q/I_1 \oplus Q/I_2 \lra Q/(I_1+I_2) \lra 0 \:.
\end{equation*}
The foundations of the proof of the main theorem are laid in
\secref[Sections~]{1}--\secref[]{gr} with an analysis of this
sequence. Particular attention is paid to the relations imposed by
$(\flat)$ on numerical invariants of $Q/I$, the Gorenstein rings, and
$Q/(I_1+I_2)$ as they relate to compressedness. Some of the broader
conclusions of this analysis remain valid without the assumption that
$I$ is homogeneous and, indeed,
\mbox{\secref[Sections~]{1}--\secref[]{3}} deal with general artinian
local rings. The next steps are comparative analyses of the minimal
graded free resolutions (\secref{5}) and the multiplicative structures
on the Tor algebras (\secref{parameters}) of the rings in
$(\flat)$. Central to both analyses is \conref{resR}, which harnesses
an identification, as graded $Q$-modules, of the kernel of the
homomorphism $Q/I \onto Q/I_2$ with a power of the maximal ideal of
the Gorenstein ring $Q/I_1$ and the canonical module of a quotient of
$Q$ by a power of its maximal ideal. The actual proof of the main
theorem, which takes up most of \secref{7}, also draws on various
nonhomological techniques to deal with issues not covered by the
general homological analysis in \secref{parameters}. In the final
\secref{generic} we describe the experiments that inspired this
project and discuss how the main theorem and its underpinnings explain
the generic behavior they reveal.

\section*{Acknowledgments}
\noindent
This paper benefited greatly from inspiring and clarifying
conversations we shared with Tony Iarrobino, Andy Kustin, Pedro
Marcias Marques, Liana \c{S}ega, Keller Vandebogert, and Jerzy Weyman;
we are grateful to them.

Significant progress was made in the wake of the 2020 ICERM workshop
on \emph{Free Resolutions and Representation Theory;} we thank the
organizers and participants and acknowledge the National Science
Foundation's support of this event under Grant No.~DMS-1439786.

An anonymous referee provided a detailed report that led to correction
of several typos and clarification of an argument buried deep in the
proof of the main theorem.


\section{Artinian local rings of type 2}
\label{sec:1}

\noindent
As is standard, we abbreviate the statement that a ring is local with
unique maximal ideal to, say, ``$(Q,\mfq)$ is a local ring'', and when
we need the notation $\k$ for the residue field $Q/\mfq$ we say that
$(Q,\mfq,\k)$ is local.

For $\mfq$-primary Gorenstein ideals $I_1$ and $I_2$ in a regular
local ring $(Q,\mfq)$ we identify necessary and sufficient conditions
for the ideal $I_1\cap I_2$ to define a ring of type $2$.

\begin{bfhpg}[Definitions]
  \label{defs}
  Let $(R,\mfm,\k)$ be an artinian local ring.  For every element
  $x\ne 0$ in $R$ the \emph{valuation},
  \begin{equation*}
    v_R(x) \deq \max\setof{i}{x\in \mfm^i} \:,
  \end{equation*}
  is finite.  The \emph{socle} of $R$ is the annihilator of the
  maximal ideal; it is a $\k$-vector space whose rank is called the
  \emph{type} of $R$. That~is,
  \begin{equation*}
    \SocR \deq (0:_R \mfm) \qqand 
    \type[]{R} \deq \rnk{(\SocR)}\:.
  \end{equation*}
  The ring $R$ is Gorenstein if and only if $\type[]{R}=1$ holds. The
  \emph{socle degree} of $R$ is the integer $s$ with
  $\mfm^s \ne 0 = \mfm^{s+1}$. Evidently, one has
  $\mfm^s \subseteq \SocR$; if equality holds, then $R$ is called
  \emph{level}. Following Kustin, \c{S}ega, and Vraciu
  \cite[2.3(d)]{KSV-18}, the \emph{socle polynomial} of $R$ is defined
  as
  \begin{equation*}
    \sum_{i=0}^s \rnk{\left(\frac{\mfm^i \cap \SocR }
        {\mfm^{i+1} \cap  \SocR }\right)} \chi^i\:.
  \end{equation*}
  Notice that $R$ is level if and only if its socle polynomial is a
  monomial.

  Let $R \is Q/I$ be a minimal Cohen presentation of $R$, see
  \thmcite[A.21]{bruher}, where $(Q,\mfq,\k)$ is a regular local ring
  and $I$ is an ideal of $Q$ contained in $\mfq^2$.  One has
  \begin{equation*}
    \mfq^{s+1} \dsubseteq I \qand \mfq^s \:\not\subseteq\: I\:.
  \end{equation*}
  The invariant $t$ defined by
  \begin{equation*}
    I \dsubseteq \mfq^t \qand I \:\not\subseteq\: \mfq^{t+1}
  \end{equation*}
  is, with a slight abuse of terminology, called the \emph{initial
    degree} of $I$. In case $I$ is homogeneous, it truly is the
  initial degree. Notice the inequalities
  \begin{equation}
    \label{eq:2ts}
    2 \dle t \dle s+1 \:.
  \end{equation}
\end{bfhpg}

\begin{ipg}
  \label{golod}
  Recall, for example from Avramov \seccite[5.2]{ifr}, that a local
  ring $(R,\mfm,\k)$ is called \emph{Golod} if the ranks of the
  modules in the minimal free resolution of $\k$ over $R$ attain an
  upper bound established by Serre. If $R$ is artinian, then it
  follows from an observation by L\"ofwall \thmcite[2.4]{CLf86} that
  $R$ is Golod if the inequality
  \begin{equation*}
    \left\lceil \frac{s+1}{2}\right\rceil \:<\: t
  \end{equation*}
  holds. Rossi and \c{S}ega give a different argument in the proof of
  \prpcite[6.3]{MERLMS14}.
\end{ipg}

\begin{lem}
  \label{lem:socle-valuations}
  Let $(R,\mfm,\k)$ be an artinian local ring of type 2 and socle
  degree $s$. If one has $\rnk{\mfm^s} = 1$, then there exists an
  integer $v < s$ such that every nonzero element in $\Soc{R}$ has
  valuation $v$ or $s$.
\end{lem}

\begin{prf*}
  Let $x\ne 0$ be an element in $\mfm^s$, and choose an element $y$
  such that $\set{x,y}$ is a basis for the $\k$-vector space
  $\Soc{R}$. Notice that $v = v_R(y)$ satisfies $v < s$, as every
  nonzero element of $R$ has valuation at most $s$ and
  $y\not\in \mfm^s$. Consider an element $z\in \Soc{R}$ with
  $v_R(z) < s$; it is a linear combination of $x$ and $y$, so
  $v_R(z) \ge \min\set{v,s} = v$ holds. Suppose one has $v_R(z) > v$
  and choose $\alpha,\beta\in R\setminus\mfm$ such that
  $z=\alpha x + \beta y$ holds. As $\alpha$ and $\beta$ are nonzero,
  it follows that $y = \beta^{-1}(z - \alpha x)$ has valuation at
  least $v_R(z)$; a contradiction. Thus, every nonzero element of
  $\Soc{R}$ has valuation $v$~or~$s$.
\end{prf*}

\begin{ipg}
  \label{stp0}
  Let $R$ be an artinian local ring of type 2; as in \pgref{defs}
  consider a minimal Cohen presentation $R \is Q/I$. The zero ideal of
  $R$ has two irreducible components, so in $Q$ one has
  $I = I_1\cap I_2$, where $I_1$ and $I_2$ are irreducible ideals; see
  e.g.\ Gr\"obner \cite[\S6, Satz 3]{WGr35}. As $I$ is $\mfq$-primary,
  so are $I_1$ and $I_2$. It follows that $Q/I_1$ and $Q/I_2$ are
  artinian Gorenstein rings.  This leads us to consider the following
  situation.
\end{ipg}

\begin{stp}
  \label{stp:1}
  Let $(Q,\mfq,\k)$ be a regular local ring and $I_1$ and $I_2$ be
  $\mfq$-primary ideals contained in $\mfq^2$. Set
  \begin{equation*}
    I \deq I_1 \cap I_2 \qqand I' \deq I_1 + I_2 \:;
  \end{equation*}
  these are also $\mfq$-primary ideals, and we adopt the following
  notation
  \begin{equation}
    \label{eq:notation1}
    \text{
      \begin{tabular}[c]{cccc}
        Ideal & Quotient of $Q$ & Initial degree & Socle degree \\ \hline
        $I_1$ & $R_1$ & $t_1$ & $s_1$ \\
        $I_2$ & $R_2$ & $t_2$ & $s_2$ \\
        $I $ & $R$ & $t$ & $s$ \\
        $I'$ & $R'$ & $t'$ & $s'$ \\
      \end{tabular}
    }
  \end{equation}
  Denote by $e$ the common embedding dimension of $Q$, $R_1$, $R_2$,
  $R$, and $R'$. Without loss of generality, assume that $s_1 \le s_2$
  holds and set
  \begin{equation}
    \label{eq:notation2}
    a \deq \min\setof{i \ge 0}{\mfq^iI_2 \subseteq I_1} \qand 
    b \deq \min\setof{i \ge 1}{\mfq^{i+1} \cap I_2 \subseteq I_1}\;.
  \end{equation}
\end{stp}

The numbers $a$ and $b$ introduced in \eqref{notation2} capture
crucial relations between the ideals $I_1$ and $I_2$. Once we pass to
the setting of compressed rings, $b$ merges with $s_1$, see
\thmref{R1R2loc}, and with further restriction to graded rings, $a$ is
determined by the invariants from \eqref{notation1}; see
\thmref{R1R2}. The notation $a$ remains convenient and is part of the
statement of the main result, \thmref{main}.

The first step in the analysis of \stpref{1} is to record some
elementary relations between the invariants from \eqref{notation1} and
\eqref{notation2}.

\enlargethispage*{\baselineskip}
\begin{prp}
  \label{prp:soc-val}
  Adopt the setup in \stpref[]{1}.  One has:
  \begin{equation*}
    \tag{a}
    s \deq \max\set{s_1,s_2} \deq s_2 \qqand   t \dge \max\set{t_1,t_2} \:.
  \end{equation*}
  \begin{equation*}
    \tag{b}
    s' \dle \min\set{s_1,s_2} \deq s_1\qqand t' \deq \min\set{t_1,t_2} \:.
  \end{equation*}
  In particular, if the inequalities
  $\left\lceil \frac{s_1+1}{2}\right\rceil \le t_1$ and
  $\left\lceil \frac{s_2+1}{2}\right\rceil \le t_2$ hold, then one has
  \begin{equation*}
    \tag{c}
    \left\lceil
      \frac{s+1}{2}\right\rceil \dle t \qqand \left\lceil
      \frac{s'+1}{2}\right\rceil \dle t'\,.
  \end{equation*}
  Moreover, there are inequalities
  \begin{equation*}
    \tag{d}
    1 \dle b \dle s_1 \:,
  \end{equation*}
  and if $I_2 \not\subseteq I_1$ holds, then one has
  \begin{equation*}
    \tag{e}
    t_2 \dle a + t_2 - 1 \dle b \dle s_1 \:.
  \end{equation*}
\end{prp}

\begin{prf*}
  (a): By the assumptions one has
  $\mfq^{s_2+1} \subseteq \mfq^{s_1+1} \subseteq I_1$ and
  $\mfq^{s_2+1} \subseteq I_2$, so $\mfq^{s_2+1} \subseteq I$ holds,
  while one has $\mfq^{s_2} \not\subseteq I$ as
  $\mfq^{s_2} \not\subseteq I_2$. The inequality for the initial
  degree $t$ of $I$ follows as one has
  $I \subseteq \mfq^{t_1}\cap \mfq^{t_2} = \mfq^{\max\set{t_1,t_2}}$.

  (b): The inequality for $s'$ follows as one has
  $\mfq^{s_1+1} \subseteq \mfq^{s_1+1} + \mfq^{s_2+1}\subseteq I_1 +
  I_2 = I'$. Finally, one has
  $I' \subseteq \mfq^{t_1} + \mfq^{t_2} = \mfq^{\min\set{t_1,t_2}}$
  while $I'\not\subseteq \mfq^{t_1+1}$ and
  $I'\not\subseteq \mfq^{t_2+1}$.

  (c): Under the assumptions, these inequalities are immediate from
  (a) and (b).

  (d): The definition of $s_1$ yields the nontrivial inequality.


  (e): As $I_2$ is not contained in $I_1$, one has $a \ge 1$, and that
  explains the first inequality. As $I_2 = \mfq^{t_2} \cap I_2$ is not
  contained in $I_1$ one has $t_2 \le b$. Now the second inequality
  follows, as the chains
  \begin{equation*}
    \mfq^{b-t_2+1}I_2 \dsubseteq I_2 \qqand
    \mfq^{b-t_2+1}I_2 \dsubseteq \mfq^{b-t_2+1}\mfq^{t_2}
    \deq \mfq^{b+1}
  \end{equation*}
  yield
  $\mfq^{b-t_2+1}I_2 \subseteq \mfq^{b+1} \cap I_2 \subseteq I_1$. The
  third inequality holds by part (d).
\end{prf*}

\begin{lem}
  \label{lem:s1}
  Adopt the setup in {\rm \stpref[]{1}}. The following assertions
  hold.
  \begin{prt}
  \item For every $x \in \mfq^{s}\setminus I_2$ the element $x+I$
    belongs to $\SocR$ and $v_R(x+I) = s$.
  \item For every $y \in (\mfq^b \cap I_2) \setminus I_1$ the element
    $y + I$ belongs to $\SocR$ and $v_R(y+I) = b$.
  \end{prt}
\end{lem}

\begin{prf*}
  (a) For $x \in \mfq^{s}\setminus I_2$, the element $x + I$ in $R$ is
  nonzero, so it has valuation $s$ and is evidently a socle element.
  
  (b): Let $y \in (\mfq^b \cap I_2) \setminus I_1$. The element $y+I$
  in $R$ is nonzero. By the definition of $b$ one has
  $\mfq y \subseteq I_1$, so $y+I$ is a socle element in
  $R$. Moreover, one has $v_R(y+I) \ge v_Q(y) \ge b$. To prove that
  equality holds, assume that one has $y = y' + y''$ with
  $y'\in \mfq^{b+1}$ and $y'' \in I$. It follows that $y'$ is in $I_2$
  and, therefore, by the definition of $b$ in $I_1$. Now it follows
  that $y$ is in $I_1$, a contradiction.
\end{prf*}

The blanket assumption $s_1 \le s_2$ informs the definitions in
\eqref{notation2} and is responsible for the asymmetry in the next
statement.

\begin{thm}
  \label{thm:a}
  Adopt the setup in {\rm \stpref[]{1}}, let $\mfm$ denote the maximal
  ideal $\mfq/I$ of $R$, and assume that $R_1$ and $R_2$ are
  Gorenstein.  The following conditions are equivalent.
  \begin{eqc}
  \item The ring $R$ is of type $2$.
  \item One has $I_2 \not\subseteq I_1$.
  \item One has $I_1 \not\subseteq I_2 \not\subseteq I_1$.
  \end{eqc}
  Moreover, if $R$ is of type $2$, then one has $e \ge 2$ and the next
  assertions hold.
  \begin{prt}
  \item One has $\rnk{\mfm^{s}} = 1$ if and only if $b < s$ holds, in
    which case one has
    \begin{equation*}
      v_R(z) = b \ \text{ for every } \
      z \in (\SocR) \setminus \mfm^{s}\:.
    \end{equation*}
  \item One has $\rnk{\mfm^{s}} = 2$, i.e.\ $R$ is level, if and only
    if $ b = s$ holds, in which case also $s_1 = s$ holds.
  \end{prt}
  In particular, if $R$ is of type $2$, then its socle polynomial is
  $\chi^b + \chi^{s}$ and $a \ge 1$ holds.
\end{thm}

\begin{prf*}
  Condition \eqclbl{iii} trivially implies \eqclbl{ii}.
  
  \proofofimp{i}{iii} If $\type[]{R} = 2$ holds, then $R$ is not
  Gorenstein, so $R_1 \ne R \ne R_2$ and, therefore,
  $I_1 \not\subseteq I_2 \not\subseteq I_1$ hold. In particular, one
  has $a \ge 1$.

  \proofofimp{ii}{i} As $s=s_2$ holds, see \partprpref{soc-val}{a},
  one can choose $x\in\mfq^{s} \setminus I_2$. It follows from
  \partlemref{s1}{a} that $x + I$ is a socle element of $R$ of
  valuation $s$. The assumption $I_2 \not\subseteq I_1$ implies that
  the set $(\mfq^b \cap I_2) \setminus I_1$ is not empty, so one can
  choose an element $y$ in this set. It follows from
  \partlemref{s1}{b} that $y + I$ is a socle element of $R$ of
  valuation $b$.  Now, assume towards a contradiction that the
  elements $x + I$ and $y + I$ in $\Soc{R}$ are linearly
  dependent. There exists then a unit $\gra$ in $Q$, such that
  $y - \gra x$ belongs to $I\subseteq I_2$. As $y$ is in $I_2$, this
  implies $x\in I_2$, which contradicts the choice of $x$.  Thus, the
  elements $x+I$ and $y+I$ are linearly independent, whence
  $\type[]{R} \ge 2$ holds. As the rings $R_1$ and $R_2$ are artinian
  and Gorenstein, the ideals $I_1$ and $I_2$ are irreducible.  Thus
  one has $\type[]{R} \le 2$; see \cite[\S 6, Satz 3]{WGr35}.

  If $e=1$ holds, then $Q$ is a DVR, so the ideals in $Q$ are linearly
  ordered, whence one has $I=I_1$ or $I=I_2$. Thus, the assumption
  $\type[]{R} = 2$ implies $e \ge 2$.

  (a): The first assertion is immediate; indeed it was proved above
  that the basis elements $x+I$ and $y+I$ of $\SocR$ have valuations
  $s$ and $b$. The second assertion follows from
  \lemref{socle-valuations}.

  (b): By hypothesis and \partprpref{soc-val}{a,d} one has
  $b \le s_1 \le s$.  If $R$ is level, then (a) yields $s \le b$,
  whence $s_1 = b = s$ holds. Conversely, if $b = s$ then one has
  $s_1 = s$ and $\rnk{\mfm^{s}} = 2$ as $x+I$ and $y+I$ are linearly
  independent in $\Soc{R}$.
\end{prf*}

\begin{lem}
  \label{lem:t}
  Adopt the setup in \stpref[]{1} and assume that $R_1$ and $R_2$ are
  Gorenstein. If $R$ has type $2$, then the following assertions hold.
  \begin{prt}
  \item There are inequalities
    $\left\lceil \frac{s_1+1}{2}\right\rceil \ge t_1$ and
    $\left\lceil \frac{s_2+1}{2}\right\rceil \ge t_2$.
  \item I\fsp\ $\left\lceil \frac{s_2+1}{2}\right\rceil = t_2$ holds,
    then one has $t_1 \le t_2 \le s_1 \le s_2 < 2s_1$.
  \end{prt}
\end{lem}

\begin{prf*}
  (a): If the inequality
  $\left\lceil \frac{s_1+1}{2}\right\rceil < t_1$ holds, then the
  Gorenstein ring $R_1$ is Golod by \pgref{golod} and hence a
  hypersurface; see \cite[Remark after 5.2.5]{ifr}. Since $R_1$ is
  artinian, this implies that the embedding dimension $e$ is $1$,
  which by \thmref{a} contradicts the assumption $\type[]{R} =
  2$. Thus, $\left\lceil \frac{s_1+1}{2}\right\rceil \ge t_1$ holds,
  and by symmetry so does the second inequality.

  (b): The inequality $s_1 \le s_2$ holds by assumption. In view of
  part (a), one now has
  $t_1 \le \left\lceil \frac{s_1+1}{2}\right\rceil \le \left\lceil
    \frac{s_2+1}{2}\right\rceil = t_2$. Next, if the inequality
  $t_2 \ge s_1 + 1$ holds, then one has
  $I_2 \subseteq \mfq^{t_2} \subseteq \mfq^{s_1+1} \subseteq I_1$,
  which by \thmref{a} contradicts the assumption $\type[]{R} =
  2$. Thus, $t_2 \le s_1$ holds. Finally, if the inequality
  $s_2 \ge 2s_1$ holds, then one has
  $t_2 = \left\lceil \frac{s_2+1}{2}\right\rceil \ge \left\lceil
    \frac{2s_1+1}{2}\right\rceil = s_1 +1$, and as proved right above
  this contradicts the assumption $\type[]{R} = 2$.
\end{prf*}

In the context of \thmref{a}, the next example shows that even if
$R_1$ and $R_2$ are Gorenstein of the same socle degree, the ring $R$
may not be level.

\enlargethispage*{\baselineskip}
\begin{exa}
  \label{exa:collision-1}
  Let $\k$ be a field. In the regular local ring
  $Q = \pows[\k]{x,y,z}$ with maximal ideal $\mfq = (x,y,z)$ consider
  the complete intersection ideals
  \begin{equation*}
    I_1 = (x^2, xy + z^2, y^2) \qqand I_2 = (x^2, xy+z^2, y^2+z^3)\:.
  \end{equation*}
  It is straightforward to check that one has:
  \begin{align*}
    \mfq^4 \dsubseteq I_1 &\qand (I_1 : \mfq) \deq (z^3) + I_1 \:;\\
    \mfq^4 \dsubseteq I_2 &\qand (I_2 : \mfq) \deq (z^3) + I_2
                            \deq (y^2) + I_2 \:.    
  \end{align*}
  The intersection of the two ideals is
  \begin{equation*}
    I \deq I_1 \cap I_2 \deq (x^2, xy + z^2, y^3, y^2z, yz^2)\:.
  \end{equation*}
  By \lemref{s1} the elements $z^3 + I$ and $y^2 + z^3 + I$ belong to
  the socle of $R=Q/I$. That is, one has
  \begin{equation*}
    (I : \mfq) = (y^2,z^3) + I.
  \end{equation*}
  As $I$ is homogeneous, this shows that $R$ has socle polynomial
  $\chi^2 + \chi^3$; in particular $R$ is not level.
\end{exa}

For our purposes the next result completes the analysis of \stpref{1}
for $e=2$.

\begin{prp}
  \label{prp:e2}
  Adopt the setup in \stpref[]{1} and assume that $R_1$ and $R_2$ are
  Gorenstein. If $e=2$ and $R$ is of type $2$, then $R$ and $R'$ are
  Golod.
\end{prp}

\begin{prf*}
  By \prpcite[(5.3.4)]{ifr} each of the rings $R_1$, $R_2$, $R$, and
  $R'$ is either a codimension $2$ complete intersection or Golod; in
  particular $R$ is Golod as it is of type~$2$. By \thmref{a} one has
  $I_2 \not\subseteq I_1$, so the ring
  $R' = Q/(I_1+I_2) \is R_1/((I_1+I_2)/I_1)$ is a proper quotient of
  the artinian complete intersection $R_1$ and hence not a complete
  intersection.
\end{prf*}


\section{Compressed Gorenstein local rings}
\label{sec:2}

\noindent
For a local ring $(R,\mfm,\k)$ the Hilbert function, $h_R$, is defined
by
\begin{equation*}
  h_R(i) \deq \rnk(\mfm^i/\mfm^{i+1}) \ \text{ for } \ i\ge 0 \:.
\end{equation*}
Let $e$ denote the embedding dimension of $R$. Recall that if $R$ is
regular, then its Hilbert function is given by
\begin{equation}
  \label{eq:hfQ}
  h_R(i) \deq \binom{e-1+i}{e-1} \ \text{ for } \ i \ge 0\:.
\end{equation}
One has $h_R(0) = 1$ and $ h_R(1) =e$; notice the equality
$\lgt[]{(R)} = \sum_{i=0}^\infty h_R(i)$.
If $R$ has finite length, then we refer to the finite sequence
$(h_R(0),h_R(1),\ldots)$ of nonzero values of the Hilbert function as
the \emph{h-vector} of $R$.

In the balance of this section we adopt the following setup.

\begin{stp}
  \label{stp:2}
  Let $(Q,\mfq,\k)$ be a regular local ring of embedding dimension
  $e$, and $(R,\mfm,\k)$ an artinian local ring with minimal Cohen
  presentation $R \is Q/I$; in particular, $R$ also has embedding
  dimension $e$. As in \pgref{defs} let $t$ denote the initial degree
  of $I$ and $s$ the socle degree of $R$; both invariants can be
  detected from the Hilbert function:
  \begin{equation}
    \label{eq:st-def}
    s \deq \max\setof{i}{h_R(i) \ne 0} \qqand 
    t \deq \min\setof{i}{h_R(i) \ne h_Q(i)} \:.
  \end{equation}
\end{stp}

We proceed to recall the notion of a compressed artinian Gorenstein
ring.  Compressedness of algebras was introduced by Iarrobino
\cite{AIr84}; here we refer to the more recent treatment of compressed
local Gorenstein rings in \cite{MERLMS14}.

\begin{ipg}
  \label{compressed-1}
  Assume that $R$ is Gorenstein. For every \mbox{$i\ge 0$} there is an
  inequality
  \begin{equation*}
    h_R(i) 
    \dle \min\set{h_Q(i),h_Q(s-i)} 
    \deq \min\left\{ {e-1+i\choose e-1}, {e-1+s-i\choose e-1} \right\} \:.
  \end{equation*}
  If equality holds for every $i$, then $R$ is called {\it
    compressed}, see \prpcite[4.2]{MERLMS14}.  Notice that if $R$ is
  compressed, then its $h$-vector is symmetric and unimodal, and one
  has
  \begin{equation}
    \label{eq:compressed-1}
    \begin{aligned}
      t & \deq \min\setof{i}{h_Q(i) > h_Q(s-i)} \\
      & \deq \min\left\{ i \:\middle|\: \binom{e-1+i}{e-1} >
        \binom{e-1+s-i}{e-1} \right\} \deq \left\lceil
        \frac{s+1}{2}\right\rceil.
    \end{aligned}
  \end{equation}
\end{ipg}

\begin{rmk}
  \label{rmk:lows}
  If $e=1$, then $R$ is a compressed Gorenstein ring. If $R$ is
  Gorenstein of socle degree at most $2$, then $R$ is compressed: The
  $h$-vectors are $(1,1)$ and~$(1,e,1)$.
\end{rmk}

For use in later sections, we record three technical statements about
compressed Gorenstein rings.

\begin{lem}
  \label{lem:ht}
  Assume that $R$ is compressed Gorenstein.
  \begin{prt}
  \item If $e=2$, then one has
    \begin{equation*}
      h_Q(t) - h_Q(s - t) \deq 
      \begin{cases}
        1 & \text{if $s$ is odd}\\
        2 & \text{if $s$ is even} \:.
      \end{cases}
    \end{equation*}
  
  \item If $e \ge 3$, then one has
    \begin{equation*}
      h_Q(t) - h_Q(s - t) \deq \binom{e-2+t}{e-2} +
      \begin{cases}
        0 & \text{if $s$ is odd}\\
        \binom{e-3+t}{e-2} & \text{if $s$ is even} \:.
      \end{cases}
    \end{equation*}
  \end{prt}
\end{lem}

\begin{prf*}
  From \eqref{hfQ} one gets
  \begin{equation*}
    \textstyle
    h_Q(t) - h_Q(s - t) 
    = \binom{e-1 + t}{e-1} - \binom{e-1+s-t}{e-1}\:.
  \end{equation*}
  For $s$ odd one has $s - t = t -1$ by \eqref{compressed-1}, so the
  right-hand side becomes $\binom{e-2+t}{e-2}$; for $e=1$ this equals
  $1$.  For $s$ even one has $s - t = t -2$. For $e=2$ the right-hand
  side now becomes $t+1 - (s-t +1) = 2t-s = 2$, and for $e \ge3$ one
  gets
  \begin{equation*}
    \textstyle
    \binom{e-1+t}{e-1} - \binom{e-3+t}{e-1}
    = \binom{e-1+t}{e-1} - \binom{e-2+t}{e-1} + \binom{e-2+t}{e-1} - \binom{e-3+t}{e-1}
    =  \binom{e-2+t}{e-2} + \binom{e-3+t}{e-2} \:. \qedhere
  \end{equation*}
  
\end{prf*}

\begin{prp}
  \label{prp:level}
  If $R$ is compressed Gorenstein with $s \ge 2$, then the following
  assertions hold for all integers $i$ with $2 \le i \le s$.
  \begin{prt}
  \item The ring $R/\mfm^i$ is level of socle degree $i-1$.
  \item If $e \le 2$ or $(i,s) \ne (3,3)$, then the ring $R/\mfm^i$ is
    Golod.
  \end{prt}
\end{prp}

\begin{prf*}
  If $e=1$, then $R/\mfm^i$ is trivially level, and it is Golod by
  \prpcite[5.2.5]{ifr}. Now assume that $e \ge 2$ holds.

  \proofoftag{a} Two applications of \prpcite[4.2(b)]{MERLMS14} yield
  \begin{align*}
    (0:_{R/\mfm^i} \mfm/\mfm^i) %
    &\deq (\mfm^i:_R \mfm)/\mfm^i\\
    &\deq ((0:_R\mfm^{s+1-i}):_R\mfm)/\mfm^i\\
    &\deq (0:_R\mfm^{s+2-i})/\mfm^i\\
    &\deq \mfm^{i-1}/\mfm^i\:.
  \end{align*}

  \proofoftag{b} If $e=2$ holds, then $R/\mfm^i$ is Golod by
  \thmcite[(2.2)]{LWCOVl18}. If $s\ne 3$ holds, then $R/\mfm^i$ is
  Golod by \prpcite[6.3]{MERLMS14}. Finally, without assumptions on
  $s$ the ring $R/\mfm^2$ is Golod; see \prpcite[5.2.4(1)]{ifr}.
\end{prf*}

In \prpcite[(3.3)]{LWCOVl18} we give examples of compressed Gorenstein
local rings $(R,\mfm)$ of socle degree $3$ with $R/\mfm^3$ not Golod,
cf.~\partprpref{level}{b}.

The final result of this section recognizes high powers of the maximal
ideal of a compressed Gorenstein local ring as dualizing modules of
Golod rings; through \conref{resR} it play a central role in
\secref[Sections~]{5} and \secref[]{parameters}.

\begin{prp}
  \label{prp:omega}
  If $R$ is compressed Gorenstein, then for every integer $i \ge t$
  there is an isomorphism of $Q$-modules,
  \begin{equation*}
    \mfm^i \dis \Ext[Q]{e}{Q/\mfq^{s+1-i}}{Q} \:.
  \end{equation*}
\end{prp}

\begin{prf*}
  The equality in the next display holds by
  \prpcite[4.2(b)]{MERLMS14},
  \begin{equation*}
    \mfm^i \deq (0:_R\mfm^{s+1-i}) \dis \Hom{R/\mfm^{s+1-i}}{R} \dis \Hom{Q/\mfq^{s+1-i}}{R} \:.
  \end{equation*}
  The first isomorphism is standard and the second holds as $I$ is
  contained in $\mfq^t$, and one has $s+1-i \le s+1 - t \le t$ by
  \eqref{compressed-1}. Finally, since $R$ is Gorenstein there is an
  isomorphism
  $\Hom{Q/\mfq^{s+1-i}}{R} \is \Ext[Q]{e}{Q/\mfq^{s+1-i}}{Q}$; see
  \thmcite[3.3.7]{bruher}.
\end{prf*}


\section{Compressed local rings of type $2$}
\label{sec:3}

\noindent
In \cite{KSV-18} the notion of compressedness is extended beyond
Gorenstein local rings.

\begin{ipg}
  \label{compressed-2}
  Adopt \stpref{2} and assume that $R$ has type $2$; by \pgref{stp0}
  and \thmref{a} one has $e \ge 2$. Let $b$ be as in
  \eqref{notation2}; the socle polynomial of $R$ is
  \begin{equation*}
    \chi^b + \chi^s
  \end{equation*}
  with $b=s$ if $R$ is level and $b<s$ otherwise. If the equality
  \begin{equation}
    \label{eq:compressed-20}
    \begin{aligned}
      h_R(i) &\deq \min\set{h_Q(i),h_Q(b-i) + h_Q(s-i)}\\
      &\deq \min\left\{ {e-1+i\choose e-1}, {e-1+b-i\choose e-1} +
        {e-1+s-i\choose e-1} \right\}
    \end{aligned}
  \end{equation}
  holds for every $i \ge 0$, then $R$ is called {\it compressed}; see
  \dfncite[2.5]{KSV-18}. Notice that if $R$ is compressed, then one
  has
  \begin{align*}
    t &\deq \min\setof{i}{h_Q(i) > h_Q(b-i) + h_Q(s-i)}\\
      &\deq \min\left\{ i \:\left|\:  {e-1+i\choose e-1} >
        {e-1+b-i\choose e-1} + {e-1+s-i\choose e-1} \right.\right\}\;.
  \end{align*}
  Moreover, the next inequality holds by \thmcite[4.4(c)]{KSV-18},
  \begin{equation}
    \label{eq:compressed-2}
    \left\lceil
      \frac{s+1}{2}\right\rceil \dle t \:.
  \end{equation}
  It follows that the $h$-vector of $R$ is unimodal. Indeed, the
  function $h_Q(i)$ is increasing, and for $i \ge t$ the functions
  $h_Q(b-i)$ and $h_Q(s-i)$ are decreasing.
\end{ipg}

Strict inequality in \eqref{compressed-2} implies that $R$ is Golod,
see \pgref{golod}. As Golod rings sit at one extreme of the spectrum
we are interested in, the focus of our attention is on rings with
equality $\left\lceil \frac{s+1}{2}\right\rceil = t$; they may still
be Golod.

\begin{exa}
  \label{exa:collision-2}
  The complete intersections $Q/I_1$ and $Q/I_2$ from
  \exaref{collision-1} have $h$-vectors $(1,3,3,1)$; indeed $I_1$ and
  $I_2$ are both minimally generated by $3$ elements in
  $\mfq^2$. Thus, both rings are compressed Gorenstein rings. The
  ideal $I$ has two generators in $\mfq^2\setminus \mfq^3$, and the
  socle polynomial of $Q/I$ is $\chi^2 + \chi^3$. It follows that
  $Q/I$ has $h$-vector $(1,3,4,1)$ and is a compressed artinian ring
  of type $2$.
\end{exa}

In general, compressedness of the Gorenstein rings $Q/I_1$ and $Q/I_2$
does not guarantee compressedness of $Q/I$. On the other hand, the
ring $Q/I$ may be compressed though one of the Gorenstein rings is
not.

\begin{exa}
  \label{exa:r1r2}
  Let $\k$ be a field. In the regular local ring
  $Q = \pows[\k]{x,y,z}$ with maximal ideal $\mfq = (x,y,z)$ consider
  the homogeneous complete intersection ideals
  \begin{align*}
    I_1 &\deq (x^2, y^2, z^2+xy+yz) \:, \\
    I_2 &\deq (x^2, y^2, z^2) \:, \quad\text{and} \\
    I_3 &\deq (yz, x^2+xy, y^3-z^3) \:.
  \end{align*}
  The corresponding quotient rings have $h$-vectors
  \begin{equation*}
    h_{Q/I_1} \deq (1,3,3,1) \deq h_{Q/I_2} \qand
    h_{Q/I_3} \deq (1,3,4,3,1) \:,
  \end{equation*}
  so $Q/I_1$ and $Q/I_2$ are compressed but $Q/I_3$ is not.  For the
  intersection ideals
  \begin{align*}
    I &\deq I_1 \cap I_2 \deq (x^2,y^2,xz^2-z^3,yz^2) \quad\text{and} \\
    J &\deq I_2 \cap I_3 
        \deq (yz^2,y^2z, x^2z,
        y^3-z^3, x^2y+xy^2,x^3-xy^2)
  \end{align*}
  the $h$-vectors are
  \begin{equation*}
    h_{Q/I} \deq  (1,3,4,2) \qqand h_{Q/J} \deq  (1,3,6,4,1) \:,
  \end{equation*}
  so $Q/J$ is compressed but $Q/I$ is not compressed.
\end{exa}

As we explain in \rmkref{main}, the ideal $J$ in \exaref{r1r2} can
actually not be obtained as an intersection of ideals that define
compressed Gorenstein rings. We complement this example with:

\begin{exa}
  \label{exa:r1r2a}
  In the regular ring $Q = \pows[\ZZ_2]{x,y,z}$ consider the
  homogeneous ideals
  \begin{align*}
    I_1 &\deq (xz,xy+yz,x^{3}+y^{3}+y^{2}z+z^{3}) \:, \\
    I_2 &\deq (y^{2}z,x^{2}z+z^{3},y^{3}+xz^{2},x^{3},x^{2}y^{2}) \:,
          \quad\text{and} \\
    I_3 &\deq (z^{3},y^{2}z+xz^{2}+yz^{2},x^{2}z+xyz,y^{3}+xyz,xy^{2},x^{2}y,x^{3}+yz^{2}) \:.
  \end{align*}
  Per \emph{Macaulay2} \cite{M2} these are Gorenstein ideals with
  $I_1 \cap I_2 = I_2 \cap I_3$, and this common intersection defines
  a compressed ring with $h$-vector $(1,3,6,9,4,1)$. The Gorenstein
  rings have $h$-vectors
  \begin{equation*}
    h_{Q/I_1} = (1,3,4,3,1) \:, \quad h_{Q/I_2} = (1,3,6,6,3,1)
    \:,\qand\: \quad h_{Q/I_3} = (1,3,6,3,1)
  \end{equation*}
  so $I_2$ and $I_3$ define compressed rings, but $I_1$ does not.
\end{exa}

The next theorem concludes our analysis of \stpref{1} vis-\`a-vis
compressedness. In the graded case stronger statements are available,
see \thmref{R1R2}.

\begin{thm}
  \label{thm:R1R2loc}
  Adopt the setup in {\rm \stpref[]{1}}. Assume that $R_1$ and $R_2$
  are compressed Gorenstein rings and that $R$ has type $2$. There are
  inequalities,
  \begin{equation}
    \label{eq:l1}
    \lgt[]{(R)} \dle \sum_{i=0}^{s} \min\set{h_Q(i), h_{R_1}(i) +
      h_{R_2}(i)}
  \end{equation}
  and
  \begin{equation}
    \label{eq:ts1}
    2 \dle t_1 \dle t_2 \dle s_1 \dle s \:<\: 2s_1 \:.
  \end{equation}
  If equality holds in {\rm \eqref{l1}}, then $R$ is compressed and
  the next assertions hold.
  \begin{prt}
  \item $h_R(i) = \min\set{h_Q(i), h_{R_1}(i) + h_{R_2}(i)}$ for every
    $i\ge 0$.
  \item $t =\min\setof{i}{h_Q(i) > h_{R_1}(i) + h_{R_2}(i)}$.
  \item The socle polynomial of $R$ is $\chi^{s_1} + \chi^{s}$.
  \end{prt}
\end{thm}

\begin{prf*}
  The first inequality in \eqref{ts1} holds as $I_1$ is contained in
  $\mfq^2$. The remaining inequalities follow per \eqref{compressed-1}
  from \partprpref{soc-val}{a} and \partlemref{t}{b}.

  By \thmref{a} the socle polynomial of $R$ is $\chi^b +
  \chi^{s}$. For every $i \ge 0$ one has
  \begin{equation*}
    \tag{$\ast$}
    \begin{aligned}
      \min\set{h_Q(i),h_Q(b-i) + h_Q(s-i)}
      &\dle \min\set{h_Q(i),h_Q(s_1-i) + h_Q(s_2-i)}\\
      &\dle \min\set{h_Q(i),h_{R_1}(i) + h_{R_2}(i)}\:.
    \end{aligned}
  \end{equation*}
  Indeed, the first inequality holds as $b \le s_1$ and $s=s_2$ hold
  by \partprpref{soc-val}{a,d}.  For the second inequality notice
  first that since $R_2$ is compressed one has
  $h_{R_2}(i) = h_Q(i) = \min\set{h_Q(i),h_{R_1}(i) + h_{R_2}(i)}$ for
  $i < t_2$. Further, $t_1 \le t_2$ holds by \eqref{ts1}, so for
  $i \ge t_2$ one has
  $h_{R_1}(i) + h_{R_2}(i) = h_Q(s_1-i) + h_Q(s_2-i)$ as also $R_1$ is
  compressed.  In view of $(\ast)$ one now has
  \begin{equation*}
    \tag{$\ast\ast$}
    \begin{aligned}
      \lgt[]{(R)} &\dle \sum_{i=0}^{s} \min\set{h_Q(i),h_Q(b-i) + h_Q(s-i)}\\
      &\dle \sum_{i=0}^{s} \min\set{h_Q(i),h_{R_1}(i) + h_{R_2}(i)}
    \end{aligned}
  \end{equation*}
  where the first inequality holds by \thmcite[4.4(a)]{KSV-18}.

  Further, it follows from \thmcite[4.4(b)]{KSV-18} and $(\ast\ast)$
  that $R$ is compressed if the equality
  $\lgt[]{(R)} = \sum_{i=0}^{s} \min\set{h_Q(i),h_{R_1}(i) +
    h_{R_2}(i)}$ holds, and in that case equalities hold in $(\ast)$
  for every $i \ge 0$.  This establishes (a), and (b) follows per
  \eqref{st-def}. By \thmref{a} one has $I_2 \not\subseteq I_1$, so to
  prove (c) it suffices by \partprpref{soc-val}{d} to establish the
  inequality $b \ge s_1$.  Assume first that $t > s_1$ holds. As $R$
  and $R_2$ are compressed one has
  \begin{equation*}
    \tag{$\dagger$}
    \begin{aligned}
      h_Q(s_1) &\dle h_Q(b - s_1) + h_Q(s_2 - s_1)\\
      &\:<\: h_Q(b - s_1) +  h_Q(2s_1 - s_1) \\
      &\deq h_Q(b - s_1) + h_Q(s_1)
    \end{aligned}
  \end{equation*}
  where the first inequality holds by part (b), and the strict
  inequality comes from \eqref{ts1}.  Thus $h_Q(b-s_1)$ is positive,
  which forces $b \ge s_1$.  Assume next that $t \le s_1$
  holds. Equalities hold in $(\ast)$ for every $i \ge 0$, so for
  $i=s_1$ one has $h_Q(b-s_1) = h_Q(s_1-s_1) =1$ in view of part (b),
  and that forces $b \ge s_1$.
\end{prf*}

Notice from \exaref{collision-2} that the equality in \eqref{l1} need
not hold for a compressed ring of type $2$.

\begin{prp}
  \label{prp:ts1}
  Adopt the setup in {\rm \stpref[]{1}}. Assume that $R_1$ and $R_2$
  are compressed Gorenstein and $R$ is compressed of type $2$. The
  following assertions hold.
  \begin{prt}
  \item If $e \ge 3$, then one has $t \le s_1$.
  \item If $t> s_1$, then one has $e=2$ and $s_2 = 2s_1 -1$.
  \end{prt}
\end{prp}

\begin{prf*}
  If $t > s_1$ holds then, as in the proof of \thmref{R1R2loc}, one
  has $b=s_1$, and the inequalities $(\dagger)$ in that proof yield
  $h_Q(s_1) = 1 + h_Q(s_2 - s_1)$. That is,
  \begin{equation*}
    \binom{e-1+s_1}{e-1} \deq 1 + \binom{e-1+s_2 - s_1}{e-1} \:.
  \end{equation*}
  As $s_2 - s_1 < s_1$ holds, the displayed equality can only hold if
  one has $e=2$ and $s_2 = 2s_1-1$. In particular, the inequality
  $t \le s_1$ holds for $e \ge 3$.
\end{prf*}

The next example illustrates \partprpref{ts1}{b}.

\begin{exa}
  \label{exa:ts1}
  Let $\k$ be a field. In the regular local ring $\pows[\k]{x,y}$ with
  maximal ideal $\mfq = (x,y)$ consider the homogeneous complete
  intersection ideals
  \begin{equation*}
    I_1 \deq (xy, x^2+ y^2) \qqand I_2 \deq (x^2, y^3)\:.
  \end{equation*}
  It is straightforward to check that one has
  \begin{align*}
    \mfq^3 \,\subseteq\, I_1 &\qand (I_1 : \mfq) \deq (x^2) + I_1 \:.\\
    \mfq^4 \,\subseteq\, I_2 &\qand (I_2 : \mfq) \deq (xy^2) + I_2 \:.    
  \end{align*}
  The intersection of the two ideals is
  \begin{equation*}
    I \deq I_1 \cap I_2 \deq (x^3, x^2y, y^3)\:.
  \end{equation*}
  Thus, in the notation from \stpref[]{1} one has
  $t_1 = t_2 = s_1 = 2 < 3 = t = s_2 = 2s_1-1$.
\end{exa}

\begin{rmk}
  \label{rmk:tt2}
  Adopt the setup in {\rm \stpref[]{1}} and assume that $R_2$ is
  compressed Gorenstein. If the initial degree of $I$ is larger than
  the initial degree of $I_2$, then $R$ is Golod. Indeed,
  \partprpref{soc-val}{a} and \eqref{compressed-1} yield
  $\left\lceil \frac{s+1}{2}\right\rceil = \left\lceil
    \frac{s_2+1}{2}\right\rceil = t_2 \le t$, so it follows from
  \pgref{golod} that $R$ is Golod if the strict inequality $t_2 < t$
  holds.  Even if $t=t_2$ holds, $R$ may still be Golod, see
  \partthmref{main}{g}.
\end{rmk}

\begin{lem}
  \label{lem:t2=t}
  Adopt the setup in {\rm \stpref[]{1}}. Assume that $e \ge 3$ holds,
  $R_1$ and $R_2$ are compressed Gorenstein, and $R$ is compressed of
  type $2$. The equality $\left\lceil \frac{s+1}{2}\right\rceil =t$
  holds if and only if one has
  \begin{equation*}
    \binom{e-1+s_1-t_2}{e-1} \:<\: \binom{e-2+t_2}{e-2} \,+\, 
    \begin{cases}
      0 & \text{if $s$ is odd}\\
      \binom{e-3+t_2}{e-2} & \text{if $s$ is even}\:.\\
    \end{cases}
  \end{equation*}
\end{lem}

\begin{prf*}
  As in \rmkref{tt2} one has
  $\left\lceil \frac{s+1}{2}\right\rceil = \left\lceil
    \frac{s_2+1}{2}\right\rceil = t_2 \le t$.  \partthmref{R1R2loc}{b}
  yields $t = \min\setof{i}{h_Q(i) > h_{R_1}(i) + h_{R_2}(i)}$, so
  $t \le t_2$ holds in view of \pgref{compressed-1} if and only if one
  has $h_Q(t_2) > h_Q(s_1 - t_2) + h_Q(s_2 - t_2)$, equivalently
  $h_Q(s_1 - t_2) < h_Q(t_2) - h_Q(s_2 - t_2)$.  By \partlemref{ht}{b}
  applied to $R_2$ this is the asserted inequality.
\end{prf*}

Compressed level rings of large socle degree are Golod.

\begin{thm}
  \label{thm:t2=t}
  Adopt the setup in {\rm \stpref[]{1}}. Assume that $e \ge 3$ holds,
  $R_1$ and $R_2$ are compressed Gorenstein, and $R$ is compressed of
  type $2$. If $R$ is level, then it is Golod provided that one has
  \begin{equation*}
    s \dge 2e-3 \,+\, 
    \begin{cases}
      0 & \text{if $s$ is odd}\\
      \sqrt{8(e-1)^2+ 1} & \text{if $s$ is even}\:.\\
    \end{cases}
  \end{equation*}
\end{thm}

\begin{prf*}
  As $R$ is level, $s_1=s$ holds by \partthmref{R1R2loc}{c}. First
  assume that $s$ is odd. It follows from \pgref{golod} and
  \lemref{t2=t} that $R$ is Golod if the inequality
  $\binom{e-1+s-t_2}{e-1} \ge \binom{e-2+t_2}{e-2}$ holds.  By
  \partprpref{soc-val}{a} and \eqref{compressed-1} one has
  $t_2 = \frac{s+1}{2} = \frac{s_1+1}{2}$, and substituting this
  expression for $t_2$ into the inequality it reads
  \begin{equation*}
    \binom{e-1 + \frac{s-1}{2}}{e-1} \dge \binom{e-1 + \frac{s-1}{2}}{e-2} \:.
  \end{equation*}
  Clearing common factors reduces this inequality to
  \begin{equation*}
    \frac{1}{e-1} \dge \frac{1}{1+\frac{s-1}{2}} \qtext{equivalently} s \dge 2e-3 \:.
  \end{equation*}
  Now assume that $s$ is even. In this case one has
  $t_2 = \frac{s}{2} +1$ and it follows as above that $R$ is Golod if
  one has
  \begin{equation*}
    \binom{e-2+\frac{s}{2}}{e-1} \dge \binom{e-1+\frac{s}{2}}{e-2} +
    \binom{e-2+\frac{s}{2}}{e-2} \:.
  \end{equation*}
  Clearing common factors reduces this inequality to
  \begin{equation*}
    \frac{1}{e-1} \dge \frac{e-1+\frac{s}{2}}{(\frac{s}{2}+1)\frac{s}{2}}
    + \frac{1}{\frac{s}{2}}  \qtext{equivalently}
    s^2 - 2(2e-3)s -4e(e-1) \dge 0 \:.
  \end{equation*}
  The quadratic polynomial in $s$ has one positive root:
  $2e-3 + \sqrt{8(e-1)^2+ 1}$.
\end{prf*}


\section{The sum and intersection of two graded Gorenstein ideals}
\label{sec:gr}

\noindent From here on we work in a restricted version of \stpref{1}:
We only consider homogeneous quotients of a ring of power series in
three or more variables.

\begin{stp}
  \label{stp:3}
  Let $\k$ be a field and $(Q,\mfq)$ the local $\k$-algebra of power
  series in $e \ge 3$ variables with coefficients in $\k$. Let $I_1$
  and $I_2$ be homogeneous $\mfq$-primary Gorenstein ideals contained
  in $\mfq^2$. Set $I = I_1 \cap I_2$ and $I' = I_1 + I_2$ and adopt
  the notation \eqref{notation1}, the assumption $s_1 \le s_2$, and
  the notation \eqref{notation2}.  Finally, assume that the Gorenstein
  rings $R_1$ and $R_2$ are compressed and that $R$ has type $2$; by
  \thmref{a}, \eqref{compressed-1}, and \partprpref{soc-val}{a} one
  then has
  \begin{equation}
    \label{eq:stp3}
    I_2 \not \subseteq I_1 \not \subseteq I_2 \qqand
    \left\lceil \frac{s+1}{2}\right\rceil \deq t_2 \:.
  \end{equation}
\end{stp}

Our interest is in the rings $R$ and $R'$. The assumption $e \ge 3$
has been made part of the setup as the situation is trivial in lower
embedding dimensions. Indeed, the assumption that $R$ has type $2$
rules out the possibility \mbox{$e=1$}, and if $e=2$, then $R$ and
$R'$ are Golod, see \prpref{e2}.  We start by noticing that in higher
embedding dimension, compressed Gorenstein rings are rarely complete
intersections.

\begin{prp}
  \label{prp:ci}
  Let $(Q,\mfq)$ be as in {\rm \stpref{3}} and $J \subseteq \mfq^2$ be
  a homogeneous $\mfq$-primary complete intersection ideal. If the
  ring $Q/J$ is compressed, then one has $e=3$ and
  $h_{Q/J} = (1,3,3,1)$.
\end{prp}

\begin{prf*}
  Let $\tilde{s}$ be the socle degree of $Q/J$; by
  \eqref{compressed-1} the initial degree of $J$ is
  $\tilde{t} = \left\lceil \frac{\tilde{s}+1}{2}\right\rceil$. As
  $Q/J$ is complete intersection, $J$ is minimally generated by $e$
  elements; in particular, one has
  $e \ge h_Q(\tilde{t}) - h_{Q/J}(\tilde{t}) = h_Q(\tilde{t}) -
  h_Q(\tilde{s}-\tilde{t})$; see \pgref{compressed-1}. From
  \partlemref{ht}{b} one now gets
  \begin{equation*}
    e \dge \binom{e-2+\tilde{t}\,}{e-2} \dge \binom{e}{e-2} \deq \frac{e(e-1)}{2} \;,
  \end{equation*}
  where the second inequality holds as one has $\tilde{t} \ge 2$ by
  the assumption $J \subseteq \mfq^2$.  It follows that one has
  $e \le 3$, and the opposite inequality holds by assumption. Now one
  has $e = \frac{e(e-1)}{2}$, so $\tilde{t}=2$ holds and it follows
  from \partlemref{ht}{b} that $\tilde{s}$ is odd, whence
  $\tilde{s}=3$. By \pgref{compressed-1} the $h$-vector of $Q/J$ is
  thus $(1,3,3,1)$.
\end{prf*}

\begin{rmk}
  From the Mayer--Vietoris sequence
  $0 \to R \to R_1\oplus R_2\to R' \to 0$ one gets the equalities
  \begin{equation}
    \label{eq:mv}
    h_R(i) + h_{R'}(i) \deq  h_{R_1}(i) + h_{R_2}(i) \ \text{ for all } \ i\ge 0\:.
  \end{equation}
\end{rmk}

\thmref{R1R2loc} can by way of \eqref{mv} be strengthened as follows:

\begin{thm}
  \label{thm:R1R2}
  Adopt the setup in {\rm \stpref[]{3}}. There is an inequality,
  \begin{equation}
    \label{eq:l}
    \lgt[]{(R)} \dle \sum_{i=0}^{s} \min\set{h_Q(i), h_{R_1}(i) +
      h_{R_2}(i)} \;; 
  \end{equation}
  and equality holds in if and only if $R$ is compressed.
  
  Moreover, if $R$ is compressed, then there inequalities,
  \begin{equation}
    \label{eq:ts}
    2 \dle t_1 \dle t_2 \dle t \dle s_1 \dle s \:<\: 2s_1 \:,
  \end{equation}
  and the following assertions hold.
  \begin{prt}
  \item $h_R(i) = \min\set{h_Q(i), h_{R_1}(i) + h_{R_2}(i)}$ for
    $i\ge 0$.
  \item $t =\min\setof{i}{h_Q(i) > h_{R_1}(i) + h_{R_2}(i)}$.
  \item[\textnormal{(b\prm)}] If the equality
    $h_Q(i) = h_{R_1}(i) + h_{R_2}(i)$ holds, then one has $t=i+1$.
  \item The socle polynomial of $R$ is $\chi^{s_1} + \chi^{s}$.
  \item $a = s_1 - t_2 + 1$.
  \end{prt}
\end{thm}

\begin{prf*}
  \stpref{3} is a special case of \stpref{1}, so the inequality
  \eqref{l} is a special case of \eqref{l1}. If equality holds in
  \eqref{l}, then it follows from \thmref{R1R2loc} that $R$ is
  compressed and (a), (b), and (c) hold. Further,
  \partprpref[Propositions~]{soc-val}{b} and \partprpref[]{ts1}{a}
  yield $t_2 \le t \le s_1$; together with \eqref{ts1} this
  establishes the inequalities in \eqref{ts}. To complete the argument
  it suffices to prove two claims: (1) If equality holds in \eqref{l}
  then (b\prm) holds. (2) If $R$ is compressed, then (d) and equality
  in \eqref{l}~hold.
  
  (1): Assume that $h_Q(i) = h_{R_1}(i) + h_{R_2}(i)$ holds. By (b)
  and \eqref{ts} one has $i < t \le s_1 \le s_2$ and, therefore,
  $0 < h_{R_2}(i) < h_Q(i)$.  It follows that $t_2 \le i$ holds; see
  \eqref{st-def}. As $t_1 \le t_2$ holds per \eqref{ts}, one now has
  \begin{equation*}
    h_Q(i+1) \:>\: h_Q(i) \deq h_{R_1}(i) + h_{R_2}(i) \:>\: h_{R_1}(i+1) + h_{R_2}(i+1) \:,
  \end{equation*}
  whence $i+1 = t$ holds by (b).

  (2): Now assume that $R$ is compressed. By \thmref{a} the socle
  polynomial of $R$ is $\chi^b + \chi^{s}$ and $a \ge 1$ holds.  To
  prove (d), choose a homogeneous element $g$ of $I_2$ with
  $\mfq^{a-1}g \not\subseteq I_1$ and set $d = v_Q(g)$. Choose
  $f\in\mfq^{a-1}\setminus \mfq^a$ such that $fg \not\in I_1$ and
  notice that the degree $v_Q(f)$ is exactly $a-1$.  As $\mfq(fg)$ is
  in $I_1$, the coset $fg + I_1$ is a nonzero socle element in
  $R_1$. Thus one has
  \begin{equation*}
    s_1 \deq v_{R_1}(fg + I_1) \deq v_Q(fg) \deq v_Q(f) + v_Q(g) \deq a-1+d\:.
  \end{equation*}
  As $g$ is in $I_2$ one has $d\ge t_2$, and it suffices to prove that
  equality holds. One has $a = s_1 - d + 1$ and hence
  $\mfq^{s_1-d+1}(I_2)_{t_2} \subseteq I_1$. In particular, one has
  \begin{equation*}
    (I_2)_{t_2} \:\subseteq\: (I_1:\mfq^{s_1+1-d}) \deq \mfq^{d} + I_1
  \end{equation*}
  where the equality holds by \prpcite[4.2]{MERLMS14} as $R_1$ is
  compressed. Assume towards a contradiction that $d > t_2$ holds. One
  then has $(I_2)_{t_2} \subseteq (\mfq^d + I_1)_{t_2} = (I_1)_{t_2}$
  and hence $(I_2)_{t_2} = I_{t_2}$. In particular, $t \le t_2$ holds.
  The opposite inequality $t_2 \le t$, holds by
  \partprpref{soc-val}{a}, so equality holds. As $R_2$ and $R$ are
  compressed, one now has
  \begin{equation*}
    h_Q(s-t_2) \deq h_{R_2}(t_2) \deq h_R(t_2) \deq h_Q(b-t_2) + h_Q(s-t_2)  
  \end{equation*}
  where the second equality holds as one has $(I_2)_{t_2} =
  I_{t_2}$. The displayed equalities yield $h_Q(b - t_2)=0$; that is,
  $b \le t_2-1$, which contradicts \partprpref{soc-val}{e}.

  To prove that equality holds in \eqref{l}, notice first that it
  follows from (d) and \partprpref{soc-val}{e} that $b=s_1$ holds,
  i.e.\ the socle polynomial of $R$ is $\chi^{s_1} + \chi^{s}$. Next,
  recall from \thmcite[4.4]{KSV-18} that since $R$ is compressed one
  has
  \begin{equation*}
    \lgt[]{(R)} \deq \sum_{i=0}^{s} \min\set{h_Q(i),h_Q(s_1-i) +
      h_Q(s-i)}\:.
  \end{equation*}
  It thus suffices to prove that the next equality holds for all
  $i\ge 0$,
  \begin{equation*}
    \min\set{h_Q(i), h_Q(s_1-i) + h_{Q}(s-i)} 
    \deq \min\set{h_Q(i), h_{R_1}(i) + h_{R_2}(i)} \:.
  \end{equation*}
  The inequality $t_1\le t_2$ holds by \eqref{ts}, so for $i\ge t_2$
  one has $h_{R_1}(i) = h_Q(s_1-i)$ and $h_{R_2}(i) = h_Q(s-i)$. For
  $i < t_2$ one has $h_{R_2}(i) = h_Q(i) \le h_Q(s-i)$, where the
  inequality holds as $s=s_2$ by \partprpref{soc-val}{a}. Thus, both
  minima are $h_Q(i)$.
\end{prf*}

Tracking the properties of $R$ is easier when one keeps an eye on
$R'$.

\begin{thm}
  \label{thm:R'}
  Adopt the setup in {\rm \stpref[]{3}} and assume that $R$ is
  compressed. There are inequalities
  \begin{equation}
    \label{eq:ts'}
    2 \dle t' \dle s' +1 \dle t \dle s' + 2 \:.
  \end{equation}
  Moreover, one has:
  \begin{prt}
  \item $h_{R'}(i) = \max\set{0,h_{R_1}(i) + h_{R_2}(i) - h_Q(i)}$ for
    $i\ge 0$.
  \item $s' =\max\setof{i}{h_Q(i) < h_{R_1}(i) + h_{R_2}(i)}$.
  \item If $t=t_2$ holds, then one has $s' + 1 = t$.
  \item If $ t' = s'+1$ holds, then one has $t' = t_2$.
  \end{prt}
\end{thm}

\begin{prf*}
  Part (a) follows immediately from \eqref{mv} and
  \partthmref{R1R2}{a}.

  (b): The first equality below holds by the definition of $s'$ and
  the second follows from part (a):
  \begin{equation*}
    s' \deq \max\setof{i}{h_{R'}(i)\ne 0}
    \deq \max\setof{i}{h_Q(i) < h_{R_1}(i)+h_{R_2}(i)} \:.
  \end{equation*}
  
  We can now prove the inequalities in \eqref{ts'}. The first one
  holds as $I'$ is contained in $\mfq^2$, and the second inequality
  holds by the definitions of $s'$ and $t'$; see \eqref{2ts}.  By
  \partthmref{R1R2}{b} one has
  $t = \min\setof{i}{h_Q(i) > h_{R_1}(i) + h_{R_2}(i)}$, so the third
  inequality follows from part (b). For $i \le t-1$ one has
  $h_Q(i) \le h_{R_1}(i) + h_{R_2}(i)$ by \partthmref{R1R2}{b} and per
  \partthmref[]{R1R2}{b\prm} equality can only hold for $i=t-1$. The
  inequality $t \le s' + 2$ now follows from part (b).

  (c): The inequality $s'+1 \le t$ holds by \eqref{ts'}. Assume that
  $t=t_2$ holds, as $R_2$ is compressed one then has
  $h_{R_2}(t-1) = h_{Q}(t-1)$. By \eqref{ts} one has $t\le s_1$, so
  $h_{R_1}(t-1)$ is positive, whence
  $h_{R_1}(t-1) + h_{R_2}(t-1) > h_{Q}(t-1)$ holds. From part (b) one
  now gets $s' \ge t-1$.

  (d): By \eqref{ts} one has $t_1 \le t_2$, so $t'=t_1$ holds by
  \partprpref{soc-val}{b}. Thus, if $t' = s'+1$ holds, then (b) yields
  $h_Q(t_1) \ge h_{R_1}(t_1)+h_{R_2}(t_1)$. As $h_{R_1}(t_1)$ is
  positive, this implies that $h_Q(t_1) > h_{R_2}(t_1)$ holds, so one
  has $t_1 \ge t_2$ as $R_2$ is compressed.
\end{prf*}

\begin{prp}
  \label{prp:golod}
  Adopt the setup in {\rm \stpref[]{3}} and assume that $R$ is
  compressed. The next assertions hold.
  \begin{prt}
  \item If $s_1=2$ or $s_1=3=s_2$ holds, then $R'$ is Golod of type
    $e$ and level of socle degree $1$. That is, one has
    $R' \deq Q/\mfq^2$.
  \item If $s_1=3 <s_2$ holds, then $R'$ is of type $e$ and level of
    socle degree $2$, and one has $R' \deq Q/(I_1+\mfq^3)$.

  \item If $s_1 \ge 4$ holds, then $R'$ is Golod. Moreover, if $t=t_2$
    holds, then $R'$ is of type $h_{R_1}(t-1)$ and level of socle
    degree $t-1$, and one has $R' \deq Q/(I_1+\mfq^t)$.
  \end{prt}
\end{prp}

\begin{prf*}
  \proofoftag{a} It suffices to show that $R'$ has $h$-vector
  $(1,e)$. Indeed, this implies that $R'$ is $Q/\mfq^2$, so $R'$ is
  Golod, see \prpcite[5.2.4]{ifr}, and evidently of type $e$ and level
  of socle degree~$1$.  Assume first that $s_1=2$ holds. By \eqref{ts}
  the possible values of $s_2$ are $2$ and $3$. For $s_2=2$, both
  $R_1$ and $R_2$ have $h$-vector $(1,e,1)$, so \partthmref{R'}{a}
  yields $h_{R'} = (1,e)$.  In case $s_2=3$, the $h$-vector of $R_2$
  is $(1,e,e,1)$, so \partthmref[]{R'}{a} yields $h_{R'} = (1,e)$.
  Finally, if one has $s_1 = 3 = s_2$, then $R_1$ and $R_2$ both have
  $h$-vector $(1,e,e,1)$. As $e \ge 3$, one has
  $2e \le \binom{e+1}{2} = h_Q(2)$, so \partthmref[]{R'}{a} yields
  $h_{R'} = (1,e)$.
  
  \proofoftag{b} The $h$-vector of $R_1$ is $(1,e,e,1)$. By \eqref{ts}
  the possible values of $s_2$ are $4$ and $5$, so the $h$-vector of
  $R_2$ is
  \begin{equation*}\textstyle
    \left(1,e,\binom{e+1}{2},e,1\right) \qquad\text{or}\qquad 
    \left(1,e,\binom{e+1}{2},\binom{e+1}{2},e,1\right) \:.
  \end{equation*}
  In either case it follows from \partthmref{R'}{a} that $R'$ has
  $h$-vector $(1,e,e)$; in particular, it has type $e$.  Further,
  $t_2 = 3$ holds in either case, see \eqref{compressed-1}, so the
  ideal $I' = I_1 + I_2$ is contained in $I_1 + \mfq^3$. As the
  quotients $R'$ and $Q/(I_1+\mfq^3)$ have the same $h$-vector, they
  are equal as claimed. With $\mfm_1 = \mfq/I_1$ one can rewrite the
  equality of rings as an isomorphism $R' \is R_1/\mfm_1^3$. It now
  follows from \partprpref{level}{a} that $R'$ is level of socle
  degree $2$.

  \proofoftag{c} \partprpref{soc-val}{a} yields $t \ge t_2$. First we
  assume that equality holds and argue that one has
  \begin{equation*}
    \tag{$\ast$}
    h_{R'}(i) \deq
    \begin{cases}
      h_{R_1}(i) & \text{for } i \le t-1 \\
      0 & \text{for } i\ge t \:.
    \end{cases}
  \end{equation*}
  Indeed, \partthmref{R'}{c} yields $s' = t-1$, so $h_{R'}(i) = 0$
  holds for $i \ge t$. For $i \le t-1$ the assumption $t=t_2$ implies
  the equalities $h_R(i) = h_Q(i) = h_{R_2}(i)$, so
  $h_{R'}(i) = h_{R_1}(i)$ holds by \eqref{mv}.

  The ideal $I_1 + I_2$ is contained in $I_1 + \mfq^t$, still by the
  assumption $t=t_2$, and $(\ast)$ shows that the quotients $R'$ and
  $Q/(I_1+\mfq^t)$ have the same $h$-vector, so they are equal as
  claimed. As above, one now has $R' \is R_1/\mfm_1^t$. It follows
  that $R'$ is of type $h_{R_1}(t-1)$, and by \partprpref{level}{a} it
  is level of socle degree $t-1$. Finally, as $s_1 \ge 4$ holds by
  assumption, $R'$ is Golod by \partprpref{level}{b}.

  Assuming now that $t > t_2$ holds, \eqref{ts} yields
  $t_2 \le s_1-1$.  As the $h$-vectors of $R_1$ and $R_2$ are
  symmetric and unimodal, see \pgref{compressed-1}, this explains the
  first inequality in the next display. The sharp inequality holds by
  the assumptions $e\ge 3$ and $s_1 \ge 4$, and the final equality
  holds in view of \eqref{hfQ}.
  \begin{align*}
    h_{R_1}(s_1-1) + h_{R_2}(s_1-1) 
    & \dle h_{R_1}(1) + h_{R_2}(s_1-2)\\
    &\dle e + h_{Q}(s_1-2)\\
    &\:< \textstyle \binom{e+s_1-3}{e-2} + h_{Q}(s_1-2)\\
    &\deq h_{Q}(s_1-1) \:.
  \end{align*}
  This inequality, $h_{R_1}(s_1-1) + h_{R_2}(s_1-1) < h_{Q}(s_1-1)$,
  implies that $t$ is at most $s_1 -1$; see
  \partthmref{R1R2}{b}. Combining this with inequalities from
  \eqref{ts'}, \eqref{compressed-1}, \eqref{ts}, and
  \partprpref{soc-val}{b} one gets
  \begin{equation*}
    s' + 1 \dle t \dle s_1 - 1 \dle 2t_1 -2 \deq 2t'-2 \:,
  \end{equation*}
  whence $R'$ is Golod by \pgref{golod}.
\end{prf*}

\begin{rmk}
  \label{rmk:T}
  Under the assumptions in \partprpref{golod}{b} the ring $R'$ may not
  be Golod. In the notation from the proof one has
  $R' \is R_1/\mfm_1^3$, and in \thmcite[(4.2)]{LWCOVl18} we identify
  conditions under which this quotient is not Golod.
\end{rmk}

With the next theorem, which conveniently describes $R'$ as a
truncation of the Gorenstein ring $R_1$, we transition to the setup of
\secref[Sections~]{5}--\secref[]{generic} where the focus is
exclusively on rings of embedding dimension $3$. The theorem is
important to the central \conref{resR}.

\begin{thm}
  \label{thm:golod}
  Adopt the setup in {\rm \stpref[]{3}}. Assume that $e=3$ and $R$ is
  compressed. I\fsp\ $\left\lceil \frac{s+1}{2}\right\rceil = t$
  holds, then one has $R' = Q/(I_1 + \mfq^t)$.
\end{thm}

\begin{prf*}
  By \partprpref{soc-val}{a}, the assumption on $t$, and
  \eqref{compressed-1} one has $t_2 = t$. For $s_1 \ge 4$ the
  assertion now follows from \partprpref{golod}{c}. For $s_1 = 3$ one
  has $3 \le s \le 5$, see \eqref{ts}. For $s = 3$ the $h$-vector of
  $R$ is $(1,3,6,2)$, see Table~\ref{table0}, so one has
  $\left\lceil \frac{s+1}{2}\right\rceil = 2 < 3 = t$. For $s = 4$ and
  $s = 5$ one has $R' = Q/(I_1 + \mfq^3)$ by
  \partprpref{golod}{b}. For $s_1 = 2$ one has $2 \le s \le 3$ per
  \eqref{ts}. For $s = 2$ and $s = 3$ one has, respectively,
  $h_R=(1,3,2)$ and $h_R=(1,3,4,1)$, see Table~\ref{table0}, so
  $\left\lceil \frac{s+1}{2}\right\rceil = 2 = t$ holds. As $I_1$ is
  contained in $\mfq^2$ one has $I_1+\mfq^t = \mfq^2$, so
  \partprpref{golod}{a} yields $R' = Q/(I_1 + \mfq^t)$.
\end{prf*}

\begin{ipg}
  Adopt the setup in {\rm \stpref[]{3}}; assume that $e=3$ and $R$ is
  compressed. For convenience we tabulate the $h$-vectors of
  $R$ for frequently referenced combinations of $s$ and $s_1$; the
  formula for the Hilbert function comes from \eqref{compressed-20}.
  \renewcommand{\arraystretch}{1.5}
  \begin{longtable}{c|c|c}
    $(s_1,s)$ & $h_R(i)$ & $h_R$ \\ \hline $(2,2)$ &
    $\min\left\{ {2+i \choose 2}, 2{4-i \choose 2} \right\}$ &
    $(1,3,2)$ \\
    \rowcolor{Gray} $(2,3)$ &
    $\min\left\{ {2+i \choose 2}, {4-i \choose 2} + {5-i \choose 2}
    \right\}$ &
    $(1,3,4,1)$ \\
    $(3,3)$ & $\min\left\{ {2+i \choose 2}, 2{5-i \choose 2} \right\}$
    &
    $(1,3,6,2)$ \\
    \rowcolor{Gray} $(4,4)$ &
    $\min\left\{ {2+i \choose 2}, 2{6-i \choose 2} \right\}$ &
    $(1,3,6,6,2)$ \\
    $(4,5)$ &
    $\min\left\{ {2+i \choose 2}, {6-i \choose 2} + {7-i \choose 2}
    \right\}$ &
    $(1,3,6,9,4,1)$ \\
    \rowcolor{Gray} $(6,6)$ &
    $\min\left\{ {2+i \choose 2}, 2{8-i \choose 2} \right\}$ &
    $(1,3,6,10,12,6,2)$ \\
    $(6,7)$ &
    $\min\left\{ {2+i \choose 2}, {8-i \choose 2} + {9-i \choose 2}
    \right\}$ &
    $(1,3,6,10,15,9,4,1)$ \\
    \rowcolor{Gray} $(7,9)$ &
    $\min\left\{ {2+i \choose 2}, {9-i \choose 2} + {11-i \choose 2}
    \right\}$ &
    $(1,3,6,10,15,21,13,7,3,1)$ \\
    \caption{The $h$-vectors of select compressed rings $R$.}
    \label{table0}
  \end{longtable}
  \renewcommand{\arraystretch}{1}
\end{ipg}


\section{Minimal graded free resolutions}
\label{sec:5}

\noindent
Let $\k$ be a field and $(Q,\mfq)$ the local $\k$-algebra of power
series in three variables with coefficients in $\k$.  Let
$J \subseteq \mfq^2$ be a $\mfq$-primary homogeneous ideal in $Q$ and
set $S = Q/J$. The minimal graded free resolution of $S$ over $Q$ has
the form
\begin{equation*}
  Q \lla
  \bigoplus_{j \ge 1} Q^{\beta_{1j}(S)}(-j)
  \lla \bigoplus_{j \ge 1} Q^{\beta_{2j}(S)}(-j) \lla \bigoplus_{j \ge 1}
  Q^{\beta_{3j}(S)}(-j) \lla 0 \:.
\end{equation*}
As the defining ideal $J$ is contained in $\mfq^2$ and the resolution
is minimal, the graded Betti numbers $\beta_{ij}(S)$ vanish for
$j \le i$. The Hilbert series
$H_S(\chi) = \sum_{j \ge 0} h_S(j)\,\chi^j$ is a rational function,
see \lemcite[4.1.13]{bruher},
\begin{equation}
  \label{eq:HS}
  H_S(\chi) = \frac{B_S(\chi)}{(1-\chi)^3} \:,
\end{equation}
where the polynomial $B_S(\chi) = \sum_{j \ge 0}b_S(j)\,\chi^j$ has
coefficients
\begin{equation}
  \label{eq:BS}
  b_S(0) = 1 \qand
  b_S(j) \deq - \beta_{1j}(S) + \beta_{2j}(S) - \beta_{3j}(S) \ \text{ for } j\ge 1 \:.
\end{equation}

The next result is not new---it follows from work of Boij
\prpcite[3.3]{MBj99}---but included to match \thmref{resR}.  Our proof
uses the existence of algebra structures on free resolutions over $Q$,
which we recall briefly in \pgref{ms0} and in further detail in
\pgref{ms}.

\begin{prp}
  \label{prp:resR2}
  Let $J \subseteq \mfq^2$ be a $\mfq$-primary homogeneous ideal in
  $Q$ such that the quotient $S = Q/J$ is compressed Gorenstein of
  socle degree $s$ and initial degree~$t$.

  If $s$ is odd, then the minimal graded free resolution of $S$ over
  $Q$ has the form
  \begin{equation*}
    Q \lla
    \begin{matrix}
      Q^{t+1}(-t) \\ \oplus \\ Q^\upbeta(-t-1)
    \end{matrix}
    \lla
    \begin{matrix}
      Q^\upbeta(-t-1) \\ \oplus \\ Q^{t+1}(-t-2)
    \end{matrix}
    \lla Q(-s-3) \lla 0 \:,
  \end{equation*}
  for some integer $\upbeta \ge 0$.

  If $s$ is even, then the minimal graded free resolution of $S$ over
  $Q$ has the form
  \begin{equation*}
    Q \lla Q^{2t+1}(-t)
    \lla Q^{2t+1}(-t-1) \lla  Q(-s-3) \lla 0 \:.
  \end{equation*}
\end{prp}

The odd socle degree case in \prpref[]{resR2} also follows from recent
work of Vandebogert \prpcite[3.3]{KVn-1}, who additionally shows that
the integer $\upbeta$ is at most $t$.

\begin{ipg}
  \label{ms0}
  By a result of Buchsbaum and Eisenbud \cite{DABDEs77} the minimal
  free resolution of $S$ over $Q$ has a structure of a commutative
  differential graded algebra. This structure is not unique, but the
  induced graded-commutative algebra structure on $\Tor[Q]{*}{S}{\k}$
  is unique. If $S$ is Gorenstein, then $\Tor[Q]{*}{S}{\k}$ is a
  Poincar\'e duality algebra; this is due to Avramov and Golod, see
  for example \cite[1.4.2]{LLA12} or the original \cite{LLAESG71}.
\end{ipg}

\begin{bfhpg*}[Proof of \ref{prp:resR2}]
  By the definition of $t$ one has $\beta_{1j}(S) = 0$ for $j \le t-1$
  and hence $\beta_{2j}(S)=0$ for $j \le t$ and $\beta_{3j}(S)=0$ for
  $j \le t+1$. Per \eqref{BS} one thus has
  \begin{equation}
    \tag{1}
    b_S(t) \deq -\beta_{1t}(S) \qqand b_S(t+1) \deq -\beta_{1\,t+1}(S) + \beta_{2\,t+1}(S) \:.
  \end{equation}
  The assumption $J \subseteq \mfq^2$ implies that $t$ and hence $s$
  is at least $2$. By \pgref{compressed-1} one has $h_S(i) = h_Q(i)$
  for $i \le t-1$ and, therefore,
  \begin{align*}
    (1-\chi)^3H_S(\chi) %
    \deq 1 &+
             \textstyle\sum_{i=t}^s (h_S(i) - 3h_S(i-1) + 3h_S(i-2) - h_S(i-3))\chi^i \\
           &  + (-3h_S(s) + 3h_S(s-1) - h_S(s-2))\chi^{s+1} \\
           &  + (3h_S(s) - h_S(s-1))\chi^{s+2} \\
           &  - h_S(s)\chi^{s+3} \:.
  \end{align*}
  Per \eqref{compressed-1} the equality
  $t = \left\lceil \frac{s+1}{2}\right\rceil$ holds; together with
  \eqref{HS} and the formula for $h_S(i)$ from \pgref{compressed-1} it
  yields
  \begin{equation}
    \tag{2}
    \begin{aligned}
      b_S(t) & \deq h_Q(s-t) - 3h_Q(t-1) + 3h_Q(t-2) - h_Q(t-3) \\
      & \deq \textstyle \binom{s-t+2}{2}
      - 3\binom{t+1}{2} + 3\binom{t}{2}  - \binom{t-1}{2} \\
      &\deq \textstyle \binom{s-t+2}{2} - 3t - \binom{t-1}{2} \\
      & \deq -
      \begin{cases}
        t+1 & \text{if $s$ is odd} \\
        2t+1 & \text{if $s$ is even}\:.
      \end{cases}
    \end{aligned}
  \end{equation}
  For $s \ge 3$ one has $t+1 \le s$, as the equality
  $t+1 = \left\lceil \frac{s+3}{2}\right\rceil$ holds. In view of
  \pgref{compressed-1} and \eqref{HS} one now gets
  \begin{equation}
    \tag{3}
    \begin{aligned}
      b_S(t+1) &\deq h_Q(s-t-1) - 3h_Q(s-t) + 3h_Q(t-1) - h_Q(t-2) \\
      & \deq \textstyle \binom{s-t+1}{2} - 3\binom{s-t+2}{2} +
      3\binom{t+1}{2} - \binom{t}{2} \\
      & \deq
      \begin{cases}
        0 & \text{if $s$ is odd} \\
        2t+1 & \text{if $s \ge 4$ is even} \:.
      \end{cases}
    \end{aligned}
  \end{equation}
  For $s=2$ one has $t=2$; a direct computation now yields
  \begin{equation}
    \tag{4}
    b_S(t+1) \deq 5 \deq 2t+1\quad \text{if $s=2$}\:.
  \end{equation}
      
  As $S$ is Gorenstein, $\Tor[Q]{*}{S}{\k}$ is a Poincar\'e duality
  algebra; see \pgref{ms0}. In particular, for every nonzero
  homogeneous element $\sfx \in \Tor[Q]{1}{S}{\k}$ there is a
  homogeneous element $\sfy \in \Tor[Q]{2}{S}{\k}$ with
  $\sfx\sfy \ne 0$, so one has
  \begin{equation*}
    \tag{5}
    |\sfx| \ge t\,,\quad |\sfy| \ge t + 1\,, \qand 
    s+3 \deq |\sfx\sfy| \deq |\sfx| + |\sfy| \:.
  \end{equation*}

  If $s$ is odd, then $s+3 = 2t + 2$ holds, so by $(5)$ one has
  $|\sfx| = t$ and $|\sfy| = t + 2$ or $|\sfx| = t+1 = |\sfy|$.  It
  follows from $(3)$ and $(1)$ that
  $\beta_{1\,t+1}(S) = \beta_{2\,t+1}(S)$ holds, and with the
  abbreviated notation $\upbeta$ for this number it now follows from
  $(1)$ and $(2)$ that the graded minimal free resolution of $S$ has
  the asserted format.

  If $s$ is even, then $s+3 = 2t+1$ holds, and $(5)$ yields
  $|\sfx| = t$ and $|\sfy| = t + 1$. The desired conclusion now
  follows from $(1)$--$(4)$. \qed
\end{bfhpg*}

In the balance of this section we adopt \stpref{3} with $e=3$ and
further assume that also $R$ is compressed.

\begin{lem}
  \label{lem:resR}
  I\fsp\ $\left\lceil \frac{s+1}{2}\right\rceil = t$ holds, then one
  has $a = s_1-t+1$ and there are equalities
  \begin{align*}
    \beta_{1t}(R) %
    &\deq - \binom{a+1}{2} +
      \begin{cases}
        t + 1 & \text{if $s$ is odd}\\
        2t+1 & \text{if $s$ is even}
      \end{cases}
    \\[-.7\baselineskip]
    \intertext{and}
    \beta_{1\,t+1}(R) - \beta_{2\,t+1}(R) %
    &\deq a(a+2) -
      \begin{cases}
        0 & \text{if $s$ is odd}\\
        2t+1 & \text{if $s$ is even} \:.
      \end{cases}
  \end{align*}
\end{lem}

\begin{prf*}
  By the assumption on $t$ and \partthmref{R1R2}{d} one has
  $a = s_1-t+1$.  By \eqref{compressed-20} one has $h_R(i) = h_Q(i)$
  for $i \le t-1$ and, therefore,
  \begin{align*}
    (1-\chi)^3H_R(\chi) %
    \deq 1 &\:+\:
             \sum_{i=t}^s (h_R(i) - 3h_R(i-1) + 3h_R(i-2) - h_R(i-3))\chi^i \\
           &  {} \:+\: (-3h_R(s) + 3h_R(s-1) - h_R(s-2))\chi^{s+1} \\
           &  \:+\: (3h_R(s) - h_R(s-1))\chi^{s+2} \\
           & \:-\: h_R(s)\chi^{s+3} \:.
  \end{align*}
  By \partthmref{R1R2}{c} and \eqref{compressed-20} one has
  $h_R(i) = h_{Q}(s_1 - i) + h_{Q}(s-i)$ for $i \ge t$. In view of
  \eqref{HS} this yields
  \begin{equation*}
    \tag{1}
    \begin{aligned}
      b_R(t) &\deq h_{Q}(s_1 - t) + h_{Q}(s-t) - 3h_Q(t-1) + 3h_Q(t-2) - h_Q(t-3) \\
      \quad &\deq \textstyle \binom{s_1-t+2}{2} + \binom{s-t+2}{2}
      - 3\binom{t+1}{2} + 3\binom{t}{2}  - \binom{t-1}{2} \\
      \quad &\deq \textstyle \binom{a+1}{2} + \binom{s-t+2}{2} - 3t - \binom{t-1}{2} \\
      & \deq \textstyle \binom{a+1}{2} -
      \begin{cases}
        t+1 & \text{if $s$ is odd} \\
        2t+1 & \text{if $s$ is even}\:.
      \end{cases}
    \end{aligned}
  \end{equation*}
  For $s \ge 3$ one has $t+1 \le s$, as
  $t+1 = \left\lceil \frac{s+3}{2}\right\rceil$ holds by
  assumption. In view of \eqref{compressed-20} and \eqref{HS} one now
  gets
  \begin{equation*}
    \tag{2}
    \begin{aligned}
      b_R(t+1) & \deq h_Q(s_1-t-1) + h_Q(s-t-1) - 3(h_Q(s_1-t) + h_Q(s-t)) \\
      & \hspace{3pc} {} + 3h_Q(t-1) - h_Q(t-2) \\
      & \deq \textstyle \binom{s_1-t+1}{2} + \binom{s-t+1}{2} -
      3\left(\binom{s_1-t+2}{2} + \binom{s-t+2}{2}\right)
      + 3\binom{t+1}{2}  - \binom{t}{2} \\
      & \deq \textstyle \binom{a}{2} + \binom{s-t+1}{2} -
      3\left(\binom{a+1}{2} + \binom{s-t+2}{2}\right)
      + 3\binom{t+1}{2}  - \binom{t}{2} \\
      & \deq \textstyle -a(a+2) + \binom{s-t+1}{2} - 3\binom{s-t+2}{2}
      + 3\binom{t+1}{2}  - \binom{t}{2} \\
      & \deq -a(a+2) +
      \begin{cases}
        0 & \text{if $s$ is odd} \\
        2t+1 & \text{if $s \ge 4$ is even}\:.
      \end{cases}
    \end{aligned}
  \end{equation*}
  For $s=2$ one has $s_1=2=t$ and $a=1$; a direct computation yields
  \begin{equation*}
    \tag{3}
    b_R(t+1) \deq 2 \deq -a(a+2) + 2t+1\quad \text{if $s=2$}\:.
  \end{equation*}
  As shown in the first lines of the proof of \prpref{resR2} one has
  $b_R(t) = -\beta_{1t}(R)$ and
  $b_R(t+1) = -\beta_{1\,t+1}(R) + \beta_{2\,t+1}(R)$, so $(1)$--$(3)$
  yield the asserted equalities.
\end{prf*}

We can now give a detailed description of the minimal graded free
resolution of $R$ over $Q$ in the case of interest:
$\left\lceil \frac{s+1}{2}\right\rceil = t$.

\begin{thm}
  \label{thm:resR}
  Assume that $\left\lceil \frac{s+1}{2}\right\rceil = t$ holds and
  set
  \begin{equation*}\textstyle
    \quad f_0 \deq \binom{a+1}{2} \,,\quad f_1 \deq a(a+2)\,, \qand\, f_2 \deq \binom{a+2}{2} \:.
  \end{equation*}

  If $s$ is odd, then the minimal graded free resolution of $R$ over
  $Q$ has the form
  \begin{equation*}
    Q \lla
    \begin{matrix}
      Q^{t+1-f_0}(-t) \\ \oplus \\ Q^{f_1 + \upbeta}(-t-1)
    \end{matrix}
    \lla
    \begin{matrix}
      Q^\upbeta(-t-1) \\ \oplus \\ Q^{ t+1 + f_2 }(-t-2)
    \end{matrix}
    \lla
    \begin{matrix}
      Q(-s_1-3) \\ \oplus \\ Q(-s-3)
    \end{matrix}
    \lla 0
  \end{equation*}
  for some integer $\upbeta \ge 0$.

  If $s$ is even, then the minimal graded free resolution of $R$ over
  $Q$ has the form
  \begin{equation*}
    Q \lla
    \begin{matrix}
      Q^{2t+1-f_0}(-t) \\ \oplus \\ Q^{\upbeta}(-t-1)
    \end{matrix}
    \lla
    \begin{matrix}
      Q^{ 2t+1 - f_1+ \upbeta}(-t-1) \\ \oplus \\ Q^{f_2}(-t-2)
    \end{matrix}
    \lla
    \begin{matrix}
      Q(-s_1-3) \\ \oplus \\ Q(-s-3)
    \end{matrix}
    \lla 0
  \end{equation*}
  for some integer $\upbeta \ge \max\set{0,f_1-2t-1}$.
\end{thm}

We prepare for the proof with a construction that is reused in the
next section.

\begin{con}
  \label{con:resR}
  Adopt the assumptions and notation from \thmref{resR}.  The kernel
  of the surjection $R \to R_2$ is the ideal $I_2/(I_1\cap I_2)$,
  which as a $Q$-module is isomorphic to
  $(I_1 + I_2)/I_1 = (I_1 + \mfq^t)/I_1$; see \thmref{golod}. This is
  the $t^{\mathrm{th}}$ power of the maximal ideal of the compressed
  Gorenstein ring $R_1$. One has $t \ge t_1$ by
  \partprpref{soc-val}{a}, and $a = s_1 -t + 1$ holds by
  \lemref{resR}, so as a $Q$-module, $(I_1 + \mfq^t)/I_1$ is by
  \prpref{omega} isomorphic to $D = \Ext[Q]{3}{Q/\mfq^{a}}{Q}$.

  The module $D$ is the dualizing module of the Cohen--Macaulay ring
  $Q/\mfq^{a}$. The socle degree of $Q/\mfq^a$ is $a-1$, so as a
  graded module $D$ is concentrated in degree $a+2$. As a homomorphism
  of graded $Q$-modules, the surjection $R \to R_2$ has kernel
  concentrated in degree $t$. Thus one has an exact sequence of graded
  $Q$-modules, $0 \lra D(-t-a-2) \lra R \lra R_2 \lra 0$, and the Horseshoe Lemma yields an 
  exact sequence of  graded free resolutions,
  \begin{equation}
    \label{eq:FFF}
    0 \lra F' \lra \widetilde{F} \lra F'' \lra 0 \:,
  \end{equation}
  where $F'$ and $F''$ are minimal. The resolution $F''$ of $R_2$ is
  described in \prpref{resR2}, and the resolution $F'$ of $D(-t-a-2)$
  can be described in similar detail. Indeed, one has $t+a+2 = s_1+3$;
  the resolution $F'$ of $D(-s_1-3)$ is obtained from the $Q$-dual of
  the minimal graded free resolution of $Q/\mfq^{a}$. Thus, $F'$ has
  the form
  \begin{equation}
    \label{eq:resD}
    Q^{f_0}(-t) \lla
    Q^{f_1}(-t-1)
    \lla
    Q^{f_2}(-t-2)
    \lla Q(-s_1-3) \lla 0 \:.
  \end{equation}
  Indeed, in the resolution of $Q/\mfq^{a}$ the free modules in
  degrees $0$ and $1$ have ranks $1$ and $h_Q(a) = f_2$, and the rank
  of the free module in degree $3$ is the rank of the socle of
  $Q/\mfq^{a}$, i.e.\ $h_Q(a-1) = f_0$. Finally, one has
  $f_0 +f_2 - 1 = f_1$.
\end{con}

\begin{bfhpg*}[Proof of \ref{thm:resR}]
  The minimal graded free resolution $F$ of $R$ is a direct summand of
  the complex $\widetilde{F}$ in the diagram \eqref{FFF}, where $F'$
  is described in \eqref{resD} and $F''$ in \prpref{resR2}. By
  \partthmref{R1R2}{c} the socle polynomial of $R$ is
  $\chi^{s_1} + \chi^s$, so $F_3$ is a free module of rank $2$ with
  generators in degrees $s_1+3$ and $s+3$. (Notice that this means
  that one has $F_3 = \widetilde{F}_3$.)
  
  Assume first that $s$ is odd.  It follows from \eqref{FFF} and the
  descriptions of $F'$ and $F''$ that $\widetilde{F}_1 = F_1' \oplus F_1''$ and,
  therefore, $F_1$ has generators in degrees $t$ and $t+1$ only. By
  \lemref{resR} there are $t+1-f_0$ generators in degree $t$, and with
  $\upbeta = \beta_{2\,t+1}(R)$ there are $f_1 + \upbeta$ generators
  in degree $t+1$. Similarly it follows that $F_2$ is generated in
  degrees $t+1$ and $t+2$. There are $\upbeta$
  generators in degree $t+1$; that also determines the number of
  generators in degree $t+2$, as the total rank of the free module
  $F_2$ is
  \begin{equation*}
    \rnk[R]{F_1} + \rnk[R]{F_3} - \rnk[R]{F_0} \deq
    t+1-f_0 + f_1 + \upbeta + 2 -1 \deq  t+1+f_2+\upbeta \:.
  \end{equation*}

  Assume now that $s$ is even. As in the odd case, $F_1$ has
  generators in degrees $t$ and $t+1$ only. By \lemref{resR} there are
  $2t+1-f_0$ generators in degree $t$; set
  $\upbeta = \beta_{1\,t+1}(R)$. As in the odd case, $F_2$ is
  generated in degrees $t+1$ and $t+2$ with $2t+1-f_1+\upbeta$
  generators in degree $t+1$, and that also determines that there are
  $f_2$ generators in degree $t+2$. \qed
\end{bfhpg*}


\section{Parameters of multiplication on the Tor-algebra}
\label{sec:parameters}
\noindent
Throughout this section we adopt \stpref{3}; we further assume that
$e=3$ holds and that $R$, like $R_1$ and $R_2$, is compressed. Recall
from \partthmref{R1R2}{d} the equality
\begin{equation}
  \label{eq:a}
  a = s_1 - t_2 + 1 \:.
\end{equation}
As in \thmref{resR} set
\begin{equation}
  \label{eq:f}
  \quad f_0 = \binom{a+1}{2} \,,\quad f_1 = a(a+2)\,, \qand\, f_2 = \binom{a+2}{2} \:,
\end{equation}
and recall from \eqref{resD} the relation
\begin{equation}
  \label{eq:rel}
  f_0 - f_1 + f_2 - 1 = 0 \:.
\end{equation}

\begin{ipg}
  \label{ms}
  The algebra $\sfA = \Tor[Q]{*}{R}{\k}$, see \pgref{ms0}, is
  bigraded; we refer to its homogeneous components with double
  indices: $\sfA_{i\,j} = (\sfA_i)_j = \Tor[Q]{i}{R}{\k}_j$. The
  multiplicative structure on $\sfA$ can be described in terms of
  three parameters
  \begin{equation*}
    p = \rnk{(\sfA_1\!\cdot\!\sfA_1)}\,, \quad
    q = \rnk{(\sfA_1\!\cdot\!\sfA_2)}\,, \qand
    r = \rnk{\delta}
  \end{equation*}
  where $\mapdef{\delta}{\sfA_2}{\Hom[k]{\sfA_1}{\sfA_3}}$ is defined
  by $\delta(\sfy)(\sfx) = \sfx\sfy$ for $\sfx \in \sfA_1$ and
  $\sfy \in \sfA_2$.  By \thmcite[2.1]{AKM-88} and \cite[3.4.2 and
  3.4.3]{LLA12} there exist bases
  \begin{align*}
    \sfe_1,\ldots,\sfe_m \ \text{ for $\sfA_1$}\,,\quad
    \sff_1,\ldots, \sff_{m+1} \ \text{ for $\sfA_2$}\,,\qand
    \sfg_1,\sfg_2 \ \text{ for $\sfA_3$}
  \end{align*}
  such that the multiplicative structure on $\sfA$ is one of
  following:
  \begin{equation}
    \label{eq:efg}
    \begin{alignedat}{5}
      \textbf{B}:& \ &\sfe_1\sfe_2 &= \sff_3 &
      \ \sfe_i\sff_i &= \sfg_1 \ \text{ for } \ 1\le i \le 2\\
      \clG{r}:& & & &
      \sfe_i\sff_i &= \sfg_1 \ \text{ for } \ 1\le i \le r\\
      \clH{p,q}:& & \sfe_{p+1}\sfe_i &= \sff_i \ \text{ for } \ 1\le
      i\le p\qquad & \sfe_{p+1}\sff_{p+j} &= \sfg_j \ \text{ for } \
      1\le j\le q \:.
    \end{alignedat}
  \end{equation}
  Here it is understood that all products that are not listed, and not
  determined by those listed and the rules of graded commutativity,
  are zero.

  We say that $R$ is of class $\clB$ if the multiplicative structure
  on $\sfA$ is given by $\clB$ in \eqref{efg} etc.  Notice that the
  classes $\clG{1}$ and $\clH{0,1}$ coincide as do $\clG{0}$ and
  $\clH{0,0}$. This overlap is usually avoided by only using $\clG{r}$
  with $r \ge 2$, see \cite[1.3]{LLA12}, but here it is convenient to
  refer to the first class as $\clG{1}$. For the class
  $\clG{0} = \clH{0,0}$ we only use the latter symbol; the rings of
  this class are precisely the Golod rings, and they are mostly
  referred to as such.

  The relationship between the parameters $p$, $q$, and $r$ and the
  classes is simple:
  \begin{equation}
    \label{eq:pqr}
    \begin{array}{r|ccc}
      \text{Class of $R$} & p & q & r\\
      \hline
      \clB & 1 &1 &2 \\
      \clG{r} \ [r \ge 1]& 0 &1 &r \\
      \clH{p,q} \ [q \le 2] & p & q & q\\
    \end{array}
  \end{equation}  
\end{ipg}

\begin{rmk}
  \label{rmk:classes}
  The multiplicative structures described in \pgref{ms} is the subset
  of those found in \thmcite[2.1]{AKM-88} that can be realized by
  quotient rings $Q/J$ of type $2$, where $J$ is $\mfq$-primary and
  contained in $\mfq^2$; see also \rmkref{comp}. A quotient ring $Q/J$
  of type $1$, such as $R_1$ or $R_2$, is of class $\clC{3}$ if it is
  complete intersection and otherwise of class $\clG{r}$ with $r$
  equal to the minimal number of generators of $J$. A quotient ring of
  type $3$ or higher---$R'$ is typically such a ring, see
  \prpref{golod}---may in addition to the classes in \pgref{ms} be of
  what is known as class $\clT$, and that is exactly what happens in
  \thmcite[(4.2)]{LWCOVl18}---the case discussed in \rmkref{T}.
\end{rmk}

\begin{lem}
  \label{lem:ta}
  The equality $\left\lceil \frac{s+1}{2}\right\rceil = t$ holds if
  and only if one has
  \begin{equation*}
    f_0 \deq \binom{2+s_1 - \left\lceil \frac{s+1}{2}\right\rceil}{2}  \:<\: 
    \begin{cases}
      \left\lceil \frac{s+1}{2}\right\rceil + 1 & \text{if $s$ is odd}\\
      2\left\lceil \frac{s+1}{2}\right\rceil + 1 & \text{if $s$ is even}\:.\\
    \end{cases}
  \end{equation*}
\end{lem}

\begin{prf*}
  As one has $a = s_1 - t_2 + 1$ per \eqref{a} and the equality
  $\left\lceil \frac{s+1}{2}\right\rceil = t_2$ is part of the setup,
  the assertion follows from \lemref{t2=t}.
\end{prf*}

Per \pgref{golod} the next result shows that $R$ is Golod if the
difference $s-s_1$ is not too big; this is recorded in detail in
\corref[Corollaries~]{odd} and \corref[]{even}.

\begin{prp}
  \label{prp:ta}
  One has $\left\lceil \frac{s+1}{2}\right\rceil < t$ if and only if
  the next inequality holds,
  \begin{equation*}
    s - s_1 \dle
    \begin{cases}
      \frac{s+2-\sqrt{4s+13}}{2} & \text{ if $s$ is odd} \\
      \frac{s+1-\sqrt{8s+25}}{2} & \text{ if $s$ is even} \:.
    \end{cases}
  \end{equation*}
\end{prp}

\begin{prf*}
  The inequality $\left\lceil \frac{s+1}{2}\right\rceil \le t$ holds
  by \eqref{compressed-2}. Thus it follows from \lemref{ta} that the
  sharp inequality holds if and only if one has
  \begin{equation*}
    \binom{2 + s_1- \left\lceil \frac{s+1}{2}\right\rceil}{2} \dge 
    \begin{cases}
      \left\lceil \frac{s+1}{2}\right\rceil + 1 & \text{if $s$ is odd}\\
      2\left\lceil \frac{s+1}{2}\right\rceil +1 & \text{if $s$ is even}\:.\\
    \end{cases}
  \end{equation*}
  This inequality simplifies to
  \begin{equation*}\textstyle
    \left(4+2s_1-2\left\lceil \frac{s+1}{2}\right\rceil\right)
    \left(2+2s_1-2\left\lceil \frac{s+1}{2}\right\rceil\right)
    \dge 
    \begin{cases}
      4(2 + 2\left\lceil \frac{s+1}{2}\right\rceil) & \text{if $s$ is odd}\\
      8(1 + 2\left\lceil \frac{s+1}{2}\right\rceil) & \text{if $s$ is even}\:.\\
    \end{cases}
  \end{equation*}
  If $s$ is odd, then one has
  $2\left\lceil \frac{s+1}{2}\right\rceil = s + 1$, and the inequality
  simplifies to a quadratic inequality in $(s-s_1)$:
  \begin{equation*}
    \tag{$\ast$}
    4(s-s_1)^2 - 4(s+2)(s-s_1) +s^2 - 9 \dge 0 \:.
  \end{equation*}
  If $s$ is even, then one has
  $2\left\lceil \frac{s+1}{2}\right\rceil = s + 2$, and the inequality
  simplifies to:
  \begin{equation*}
    \tag{$\ast\ast$}
    4(s-s_1)^2 - 4(s+1)(s-s_1) +s^2 - 6s -24 \dge 0 \:.
  \end{equation*}
  The asserted bounds come from the meaningful (smaller) roots of the
  quadratic polynomials corresponding to $(\ast)$ and $(\ast\ast)$.
\end{prf*}

The central results of this section are \prpref[Propositions~]{q} and
\prpref[]{r} and \corref{level}. Between them they show that $R$, with
exception for a few special cases, is Golod or of class $\clG{r}$. The
proofs rely on the next two lemmas.

\begin{lem}
  \label{lem:p}
  I\fsp\ $3 \le s_1$ holds, then one has
  $\Tor[Q]{1}{R}{\k} \cdot \Tor[Q]{1}{R}{\k} = 0$.
\end{lem}

\begin{prf*}
  By \eqref{compressed-2} there is an inequality
  $\left\lceil \frac{s+1}{2}\right\rceil \le t$. If strict inequality
  holds, then $R$ is Golod by \pgref{golod}, i.e.\ of class
  $\clH{0,0}$; in particular, one has $p=0$, cf.~\eqref{pqr}.  We now
  assume that $\left\lceil \frac{s+1}{2}\right\rceil = t$
  holds. Together with the assumption $3 \le s_1$ and \prpref{ta} this
  implies that $s$ is at least $4$ and, therefore, $3 \le t$. The
  minimal graded free resolution of $R$ over $Q$ is given by
  \thmref{resR}.  In the bigraded $\k$-algebra
  $\sfA = \Tor[Q]{*}{R}{\k}$, the internal degree of a product of
  nonzero homogeneous elements from $\sfA_1$ is at least $2t$. Since
  $\sfA_2$ is concentrated in internal degrees $t+1$ and $t+2$ and the
  inequality $2t > t+2$ holds, one has $\sfA_1 \cdot \sfA_1 =0$.
\end{prf*}

\begin{lem}
  \label{lem:q}
  If one has $s_1=4$ and $s=5$, then the internal degree of a nonzero
  element in the subspace $\Tor[Q]{1}{R}{\k} \cdot \Tor[Q]{2}{R}{\k}$
  of $\Tor[Q]{3}{R}{\k}$ is $8=s+3$.
\end{lem}

\begin{prf*}
  The $h$-vector of $R$ is $(1,3,6,9,4,1)$, see Table~\ref{table0}.
  In particular, one has
  $\left\lceil \frac{s+1}{2}\right\rceil = 3 = t$, so \eqref{a} yields
  $a=2$, and the minimal graded free resolution of $R$ over $Q$ is
  given by \thmref{resR} with $f_0 = 3$, $f_1=8$, and $f_2 = 6$; see
  \eqref{f}.  As one has $t+1-f_0=1$, the ideal $I$ has a single
  generator of degree $3$, and the remaining generators have degree
  $4$. The single generator in degree $3$ generates a subspace of rank
  $3$ in $I_4$, so since one has $h_R(4) = 4$ compared to
  $h_Q(4) = 15$, there are $(15-3)-4 = 8$ generators of degree $4$.
  In particular, one has $\upbeta = 8 - f_1 = 0$.  Thus the minimal
  graded free resolution of $R$ over $Q$ has the form
  \begin{equation*}
    Q \lla Q(-3) \oplus Q^{8}(-4) \lla Q^{10}(-5) \lla  Q(-7) \oplus Q(-8) \lla 0 \:.
  \end{equation*}
  In the bigraded $\k$-algebra $\sfA = \Tor[Q]{*}{R}{\k}$, the
  internal degree of a product of non\-zero homogeneous elements
  $\sfx \in \sfA_1$ and $\sfy \in \sfA_2$ is
  $|\sfx\sfy| = |\sfx| + |\sfy| \ge 3+5 =8$, so if $\sfx\sfy$ is
  nonzero, then equality must hold as $\sfA_3$ is concentrated in
  degrees $7$ and $8$.
\end{prf*}

\begin{prp}
  \label{prp:q}
  If the inequalities $3 \le s_1 < s$ hold, then the internal degree
  of a nonzero element in the subspace
  $\Tor[Q]{1}{R}{\k} \cdot \Tor[Q]{2}{R}{\k}$ of $\Tor[Q]{3}{R}{\k}$
  is $s+3$; in particular, the subspace has rank at most
  $1$. Moreover, $R$ is Golod or of class~$\clG{r}$.
\end{prp}

\begin{prf*}
  Set $\sfA = \Tor[Q]{*}{R}{\k}$. It
  $\left\lceil \frac{s+1}{2}\right\rceil < t$ holds, then $R$ is Golod
  by \pgref{golod}, so there are no nonzero products in
  $\sfA_{\ge 1}$; see \pgref{ms}.  We may thus assume that
  $\left\lceil \frac{s+1}{2}\right\rceil =t$ holds. The minimal graded
  free resolution of $R$ over $Q$ is given by \thmref{resR}; it shows
  that the bigraded $\k$ algebra $\sfA$ decomposes as follows:
  \begin{equation*}
    \sfA \deq \sfA_{0\,0} \:\oplus\: \left(\sfA_{1\,t} \oplus \sfA_{1\,t+1} \right)
    \:\oplus\: \left(\sfA_{2\,t+1} \oplus \sfA_{2\,t+2} \right)
    \:\oplus\: \left(\sfA_{3\,s_1+3} \oplus \sfA_{3\,s+3} \right)\,.
  \end{equation*}
  The internal degree of a product of nonzero homogeneous elements
  from $\sfA_1$ and $\sfA_2$ is at least $2t+1$.  If $s$ is even, then
  one has $2t+1 = s + 3 > s_1 + 3$, so a nonzero product in
  $\sfA_1 \cdot \sfA_2$ has degree $s+3$; in particular $q\le 1$
  holds.  If $s$ is odd, then one has $2t+1 = s + 2$, and the
  assumptions imply that $s$ is at least $5$. If $s \ge s_1+2$ holds,
  then one has $2t+1 > s_1+3$, so a nonzero product in
  $\sfA_1 \cdot \sfA_2$ has degree $s+3$; in particular $q\le1$ holds.
  If $s = s_1 +1$ holds, then the assumption
  $\left\lceil \frac{s+1}{2}\right\rceil =t$, via \prpref{ta}, yields
  $s=5$; now invoke \lemref{q}.

  By \lemref{p} one has $p=0$, and as argued above $q$ is at most $1$,
  so $R$ is Golod or of class $\clG{r}$; see \pgref{ms}.
\end{prf*}

\begin{rmk}
  \label{rmk:A}
  In case $R$ is of class $\clG{r}$ with $r \ge 2$, the algebra
  $\sfA = \Tor[Q]{*}{R}{\k}$ is described in \cite[1.3]{LLA12} as a
  trivial extension of a Poincar\'e duality algebra $\sfP$, of total
  rank $2(r+1)$, by a graded $\k$-vector space $\sfV$ with the almost
  trivial $\sfP$-module structure: $\sfP_{\ge 1}\sfV = 0$. For $R$ of
  class $\clH{0,1} = \clG{1}$, see \pgref{ms}, $\sfA$ is described
  similarly with
  $\sfP = \tp[\k]{(\k \ltimes \Shift[2]{\k})}{(\k\ltimes \Shift{\k})}$
  of total rank $4$, and for $R$ of class $\clH{0,0}$ with
  $\sfP = \k \ltimes \Shift{\k}$ of rank $2$.
\end{rmk}

\rmkref{classes} explains why we use $r_2$ below to denote the minimal
number of generators of the defining ideal $I_2$ of the Gorenstein
ring $R_2$.

\begin{prp}
  \label{prp:r}
  Assume that $\left\lceil \frac{s+1}{2}\right\rceil = t$ holds and
  let $m$ and $r_2$ denote the minimal number of generators of the
  ideals $I$ and $I_2$.  The next inequality holds:
  \begin{prt}
  \item $r \dge m -f_1$.
  \end{prt}
  Moreover, if one has $s_1 < s$, then the following (in)equalities
  hold:
  \begin{prt}\stepcounter{prt}
  \item $r \dle r_2 - f_0$.
  \item $r \deq m-f_1$ if $s$ is odd.
  \item $r \dle m-f_1+f_0$ if $s$ is even.
  \item $r \deq m-f_1 = m-3$ if $s$ is even and $s_1=\frac{s}{2} + 1$.
  \end{prt}
\end{prp}

\begin{prf*}
  As in \conref{resR} let $D$ be the canonical module of $Q/\mfq^a$
  and set
  \begin{equation*}
    \sfA' \deq \Tor[Q]{*}{D(-s_1-3)}{\k}\,, \quad
    \sfA \deq \Tor[Q]{*}{R}{\k}\,, \qand \sfA'' \deq \Tor[Q]{*}{R_2}{\k} \:.
  \end{equation*}
  The ranks of the $\k$-vector spaces $\sfA'_i$, $\sfA_i$, and
  $\sfA''_i$ are given by \eqref{resD}, \thmref{resR}, and
  \prpref{resR2}; in particular one has
  \begin{equation*}
    \rnk{\sfA_1}  \deq m \deq \rnk{\sfA_2 - 1} \qqand
    \rnk{\sfA''_1} \deq r_2 \deq \rnk{\sfA''_2}
    \:.
  \end{equation*}
  The exact sequence $0 \lra D(-t-a-2) \lra R \lra R_2 \lra 0$ induces
  an exact sequence of $\k$-vector spaces
  \begin{equation*}
    \tag{$\dagger$}
    \begin{gathered}
      \xymatrix@C=1.8pc@R=2.2pc{ 0 \ar[r] & \sfA'_3
        \ar[rr]_-{(1)}^-{\phi_3} && \sfA_3 \ar[rr]_-{(1)}^-{\psi_3} &&
        \sfA_3''
        \ar `r[rd]`[l] `[llllld]^-{(0)} `[dl] [dllll] \\
        & \sfA_2' \ar[rr]_-{(f_2)}^-{\phi_2} && \sfA_2
        \ar[rr]_-{(m+1-f_2)}^-{\psi_2} && \sfA''_2
        \ar `r[rd] `[l] `[llllld]^-{(r_2 - m -1 + f_2)} `[dl] [dllll] & \\
        & \sfA'_1 \ar[rr]_-{(m-r_2+f_0)}^-{\phi_1} && \sfA_1
        \ar[rr]_-{(r_2-f_0)}^-{\psi_1} && \sfA''_1
        \ar `r[rd] `[l] `[llllld]^-{(f_0)} `[dl] [dllll] & \\
        & \sfA'_0 \ar[rr]_-{(0)}^-{\phi_0} && \sfA_0
        \ar[rr]_-{(1)}^-{\psi_0} && \sfA''_0 \ar[r] & 0 }
    \end{gathered}
  \end{equation*}
  where the numbers under the arrows indicate the ranks of the maps.
  These ranks are computed by way of \eqref{rel}. The minimal free
  resolutions $F$ and $F''$ are differential graded algebras, see
  \pgref{ms0}, and the map $R \to R_2$ lifts to a morphism
  \begin{equation*}
    \tag{$\ast$}
    F \lra F''
  \end{equation*}
  of such algebras. Therefore,
  $\psi = (\mapdef{\psi_i}{\sfA_i}{\sfA_i''})_{0\le i\le 3}$ is a
  morphism of bigraded $\k$-algebras.

  Recall the map $\mapdef{\delta}{\sfA_2}{\Hom[k]{\sfA_1}{\sfA_3}}$
  from \pgref{ms} and consider the map
  \begin{align*}
    \dmapdef{\tilde\delta}{\psi_2(\sfA_2)}{\Hom[\k]{\psi_1(\sfA_1)}{\sfA_3''}}
    \\[-.7\baselineskip]
    \intertext{given by}
    \tilde\delta(\psi_2(\sfy))(\psi_1(\sfx)) \deq \psi_1(\sfx)\psi_2(\sfy) \deq
    \psi_3(\sfx\sfy) \:.
  \end{align*}
  The rank $\tilde{r}$ of $\tilde\delta$ is a lower bound for $r$, the
  rank of $\delta$: To see this, notice that if
  $\tilde\delta(\psi_2(\sfy_1)), \ldots,
  \tilde\delta(\psi_2(\sfy_{\tilde{r}}))$ are linearly independent
  maps in $\Hom[\k]{\psi_1(\sfA_1)}{\sfA_3''}$, then
  $\delta(\sfy_1), \ldots, \delta(\sfy_{\tilde{r}})$ are linearly
  independent in $\Hom[\k]{\sfA_1}{\sfA_3}$.
  
  \subparagraph\textbf{\textit{Part} (a):} Since $\sfA''$ is a
  Poincar\'e duality algebra, see \pgref{ms0}, the subspace
  $\psi_1(\sfA_1)$ determines a subspace $\sfU$ of $\sfA_2''$ of the
  same rank, such that for every $\sfx\in\psi_1(\sfA_1)$ one has
  $\sfx\sfu \ne 0$ for some $\sfu\in\sfU$. Thus, a lower bound for
  $\tilde{r}$ is the rank of the intersection
  $\psi_2(\sfA_2) \cap \sfU$, which via \eqref{rel} is at least
  \begin{equation*}
    \rnk{\psi_1} + \rnk{\psi_2} - r_2 \deq
    m+1 -f_2 - f_0 \deq m - f_1 \:.
  \end{equation*}
  Part (a) facilitates the proof an auxiliary result akin to
  \lemref{q}.

  \begin{bfhpg*}[Sublemma]
    \sl If one has $s_1=2$ and $s=3$, then the internal degree of a
    nonzero element in the subspace $\sfA_1 \cdot \sfA_2$ of $\sfA_3$
    is $6=s+3$.
  \end{bfhpg*}

  \begin{prf*}
    The $h$-vector of $R$ is $(1,3,4,1)$, see Table~\ref{table0}, so
    $\left\lceil \frac{s+1}{2}\right\rceil = 2 = t$ holds. Now
    \eqref{a} yields $a=1$, and the ranks of the nonzero components of
    the bigraded algebra $\sfA$ are given by \thmref{resR} with
    $f_0 = 1$, $f_1 = 3$, and $f_2 = 3$; see \eqref{f}.
    \begin{alignat*}{3}
      \rnk{\sfA_{3\, 5}}  &\deq 1 & &\qand & \rnk{\sfA_{3\,6}} & \deq 1 \\
      \rnk{\sfA_{2\, 3}}  &\deq \upbeta & &\qand & \rnk{\sfA_{2\,4}} & \deq 6   \\
      \rnk{\sfA_{1\, 2}} &\deq 2 & &\qand & \rnk{\sfA_{1\,3}} & \deq 3
      + \upbeta \:.
    \end{alignat*}
    If $\upbeta=0$ holds, then a homogeneous product in
    $\sfA_1 \cdot\sfA_2 \subseteq \sfA_3$ has internal degree at least
    $2 +4 = 6 = s+3$, so a nonzero product has degree $s+3$.  If
    $\upbeta$ is positive, then one has $r \ge 3$ as part (a) yields
    $r \ge m - f_1 = 2+\upbeta$. One now has $q \le 2 < r$, which
    means that $R$ is of class $\clG{r}$, see \eqref{pqr}; in
    particular $q=1$ holds, so all nonzero products of homogeneous
    elements from $\sfA_1$ and $\sfA_2$ have the same internal
    degree. If this degree were $5$, then the parameter $r$ would be
    limited by the rank of $\sfA_{1\,2}$, which is $2$; a
    contradiction. Thus nonzero products have degree $6=s+3$.
  \end{prf*}
  
  \noindent Assume now that $s_1 < s$ holds. For $3 \le s_1$ it
  follows from \prpref{q} that a nonzero product of homogeneous
  elements from $\sfA_1$ and $\sfA_2$ has internal degree $s+3$; in
  particular, $q$ from \pgref{ms} is at most $1$. The Sublemma above
  shows that the same is true for $s_1=2$, as the assumption $s_1<s$
  and \eqref{ts} forces $s=3$.  It follows that $\sfA$ has a
  Poincar\'e duality subalgebra $\sfP$ of total rank $2r+2$, which
  captures all nontrivial products of elements from $\sfA_1$ and
  $\sfA_2$; see \eqref{efg}. For homogeneous elements
  $\sfx \in \sfP_1$ and $\sfy \in \sfP_2$ with $\sfx\sfy \ne 0$ in
  $\sfA_3$ one has $|\sfx\sfy| = s +3$. In particular, $\sfx\sfy$ is
  not in the image of $\phi_3$ whence $\psi_3(\sfx\sfy) \ne 0$. As
  $\psi$ is a morphism of graded $\k$-algebras, it follows that the
  restriction of $\psi$ to $\sfP$ is injective.

  \subparagraph\textbf{\textit{Parts} (b) \textit{and} (d):} The ranks
  of the graded components of $\psi$ are computed in $(\dagger)$, so
  one~has
  \begin{align*}
    r &\deq \rnk{\sfP_1} \dle \rnk{\psi_1} \deq r_2-f_0 \quad\text{and}\\
    r &\deq \rnk{\sfP_2} \dle \rnk{\psi_2} \deq m+1-f_2 \:.
  \end{align*}
  In view of \eqref{rel} these are the asserted bounds on $r$.

  \subparagraph\textbf{\textit{Part} (e):} As the restriction of
  $\psi$ to $\sfP$ is injective, the rank $\tilde{r}$ of
  $\tilde\delta$ equals $r$. Set
  $\sfW = \setof{\sfw\in\sfA_2''}{\psi_1(\sfA_1)\cdot \sfw = 0}$, and
  notice that because $\sfA''$ is a Poincar\'e duality algebra, one
  has
  \begin{equation*}
    \tilde{r} \dle \rnk{\psi_2} - \rnk{(\psi_2(\sfA_2) \cap \sfW)} \:.
  \end{equation*}
  Assume that $s$ is even and $s_1 = \frac{s}{2} + 1 < s$ holds. By
  \eqref{compressed-1} one has $t_2 = s_1$, which per \eqref{a}
  implies $a=1$. Now $(\dagger)$ and \eqref{f} yield
  $\rnk{\psi_2} = m+1-f_2 = m-2$, while the lower bound on $r$ from
  part (a) is $m-3$. Thus one has
  \begin{equation*}
    m-3 \dle r \deq \tilde{r} \dle m - 2 - \rnk{(\psi_2(\sfA_2) \cap \sfW)} \:,
  \end{equation*}
  so it suffices to show that the intersection
  $\psi_2(\sfA_2) \cap \sfW$ is nonzero.

  By \thmref{resR} the ideal $I$ has $2t$ minimal generators, say,
  $x_1,\ldots, x_{2t}$ of degree $t$, and without loss of generality
  one can assume that they are the first $2t$ of the $2t+1$ minimal
  generators of $I_2$; see \prpref{resR2}. By \thmcite[2.1]{DABDEs77}
  the ideal $I_2$ can be generated by the submaximal Pfaffians of a
  $(2t+1) \times (2t+1)$ skew-symmetric matrix $T$, and $R_2$ has a
  minimal free resolution $G$ over $Q$ with $\partial_2^{G} = T$; see
  also \thmcite[3.4.1]{bruher}. A standard change of basis argument
  shows one can assume that $\partial_2^{F''}$ is given by a
  $(2t+1)\times (2t+1)$ skew-symmetric matrix. Further, with
  $\sfe_i''$ and $\sff_i''$ denoting the homology classes generated by
  the basis elements for $F_1''$ and $F_2''$, the only nonzero
  products of elements in $\sfA_1''$ and $\sfA_2''$ are
  $\sfe_i''\sff_i''$ for $1 \le i \le 2t+1$. Now the image
  $\psi_1(\sfA_1)$ is spanned by $\sfe_1'',\ldots,\sfe_{2t}''$, and
  $\sfW$ is spanned by $\sff_{2t+1}''$, so the goal is to show that
  this vector is in $\psi_2(\sfA_2)$. As $\partial_2^{F''}$ is given
  by a skew-symmetric matrix, the last---i.e.\ the
  $2t+1^\mathrm{st}$---element in the basis for $F_2''$ corresponds to
  a minimal relation between the first $2t$ generators,
  $x_1,\ldots, x_{2t}$, of $I_2$.  As the entries in the
  skew-symmetric matrix are linear, see \thmref{resR}, the same
  relation is a minimal syzygy of $x_1,\ldots, x_{2t}$ as generators
  of $I$. Thus one can choose the resolution $F$ in such a way that
  one of the basis vectors in $F_2$ corresponds to this relation. Now
  the last vector in the basis for $F_2''$ is in the image of the map
  $F_2 \to F_2''$ from $(\ast)$, whence $\sff_{2t+1}''$ is in the
  image of $\psi_2$.

  \subparagraph\textbf{\textit{Part} (c):} Assume now that $s$ is
  odd. The ranks of the individual components $\sfA'_{ij}$,
  $\sfA_{ij}$, and $\sfA''_{ij}$ are also determined by \eqref{resD},
  \thmref{resR}, and \prpref{resR2}.  In greater detail, the exact
  sequence $(\dagger)$ now has the form,
  \begin{equation*}
    \xymatrix@C=1.5pc@R=2.2pc{
      0 \ar[r] &
      \sfA'_{3\,s_1+3}
      \ar[rr]_-{(1\;0)}^-{\phi_3}
      && \sfA_{3\,s_1+3} \oplus \sfA_{3\,s+3} \ar[rr]_-{\bidegrees{0\\1}}^-{\psi_3} &&
      \sfA''_{3\,s+3} 
      \ar `r[rd]`[l] `[llllld]^-{(0)} `[dl] [dllll] \\
      & \sfA'_{2\,t+2} \ar[rr]_-{(0\;f_2)}^-{\phi_2}
      && \sfA_{2\,t+1} \oplus \sfA_{2\,t+2} \ar[rr]_-{\bidegrees{\upbeta & 0 \\ 0 & t+1}}^-{\psi_2} &&
      \sfA''_{2\,t+1} \oplus \sfA''_{2\,t+2} 
      \ar `r[rd] `[l] `[llllld]^-{\binom{\upbeta''-\upbeta}{0}} `[dl] [dllll] & \\
      & \sfA'_{1\,t+1} \ar[rr]_-{(0\;f_1-\upbeta''+\upbeta)}^-{\phi_1}
      && \sfA_{1\,t} \oplus \sfA_{1\,t+1} \ar[rr]_-{\bidegrees{t+1-f_0 & 0 \\ 0 & \smash{\upbeta''}}}^-{\psi_1} &&
      \sfA''_{1\,t} \oplus \sfA''_{1\,t+1}
      \ar `r[rd] `[l] `[llllld]^-{\binom{f_0}{0}} `[dl] [dllll] & \\
      & \sfA'_{0\,t} \ar[rr]_-{(0)}^-{\phi_0}
      && \sfA_{0\,0} \ar[rr]_-{(1)}^-{\psi_0} && \sfA''_{0\,0} \ar[r] & 0
    }
  \end{equation*}
  where the (arrays of) numbers below the arrows indicate the ranks of
  the individual components of the maps. In particular, $\upbeta''$ is
  the rank of $\sfA''_{1\,t+1}$; see \prpref{resR2}. As established
  after the Sublemma, a nonzero product $\sfx\sfy$ of homogeneous
  elements $\sfx\in\sfA_1$ and $\sfy\in\sfA_2$ has degree
  $s + 3 = 2t +2$. In the next computation, the last two equalities
  follow from \thmref{resR}:
  \begin{align*}
    r &\deq \rnk{\sfP_{1\,t}} + \rnk{\sfP_{1\,t+1}} \\
      &\deq
        \rnk{\sfP_{1\,t}} + \rnk{\sfP_{2\,t+1}} \\
      &\dle \rnk{\sfA_{1\,t}} + \rnk{\sfA_{2\,t+1}} \\
      &\deq  t+1-f_0 + \upbeta \\
      & \deq m-f_1 \:.
  \end{align*}
  The opposite inequality, $r \ge m-f_1$, holds by part (a).
\end{prf*}

For certain socle polynomials, $\chi^{s_1}+\chi^{s}$, the arithmetic
constraints imposed by $s_1$ and $s$ completely determine the class of
$R$ as well as the number of minimal generators of the defining ideal
and their~degrees.

\begin{prp}
  \label{prp:fam-odd}
  If one has $s = k(k+1) - 1$ and $s_1 = \frac{1}{2}k(k+3) - 1$ for
  some integer $k \ge 2$, then $t = \frac{1}{2}k(k+1)$ holds and $I$
  is generated by one form of degree $t$ and $k(k+2)$ forms of degree
  $t+1$. Moreover, $R$ is of class $\clG{1}$.
\end{prp}

\begin{prf*}
  With $s$ and $s_1$ as given, one has
  $a \deq \frac{1}{2}k(k+3) - \frac{1}{2}k(k + 1) = k$. In the
  notation from \thmref{resR} one has
  $f_0 = \frac{1}{2}k(k + 1) = \left\lceil \frac{s+1}{2}\right\rceil$,
  whence \lemref{ta} yields
  $\left\lceil \frac{s+1}{2}\right\rceil = t$, and we notice that
  $t = f_0$ holds. From \thmref{resR} it now follows that $I$ has
  exactly one generator of degree $t$. In turn, this implies that the
  first syzygy of $I$ has no generators of degree $t+1$, i.e.\
  $\upbeta=0$ in \thmref[]{resR}, and it follows that the number of
  generators in degree $t+1$ is $f_1 = k(k + 2)$. Finally, one has
  $4 \le s_1 < s$, so \partprpref{r}{c} yields $r = m-f_1 = 1$, whence
  $R$ belongs to class $\clG{1}$ by \prpref{q} and \eqref{pqr}.
\end{prf*}

\begin{prp}
  \label{prp:fam-even}
  If $s$ is even and one has $s = \frac{1}{2}k(k+1) - 2$ and
  $s_1 = \frac{1}{4}k(k+5) - 1$ for some integer $k \ge 4$, then
  $t=\frac{1}{4}k(k + 1)$ holds and $I$ is generated by one form of
  degree $t$ and $\frac{1}{2}k(k+3) -1$ forms of degree
  $t+1$. Moreover, $R$ is Golod.
\end{prp}

\begin{prf*}
  With $s$ and $s_1$ as given, one has
  $a \deq \frac{1}{4}k(k+5) - \frac{1}{4}k(k + 1) = k$. In the
  notation from \thmref{resR} one has
  $f_0 = \frac{1}{2}k(k + 1) = 2\left\lceil
    \frac{s+1}{2}\right\rceil$, whence \lemref{ta} yields
  $\left\lceil \frac{s+1}{2}\right\rceil = t$, and we notice that
  $2t = f_0$ holds. From \thmref{resR} it now follows that $I$ has
  exactly one generator of degree $t$. In turn, this implies that the
  first syzygy of $I$ has no generators of degree $t+1$, i.e.\
  $\upbeta - f_1 + 2t + 1=0$ in \thmref[]{resR}. Thus, the number of
  generators in degree $t+1$ is
  \begin{equation*}\textstyle
    \upbeta \deq f_1 - 2t - 1 \deq k(k +
    2) - \frac{1}{2}k(k + 1) - 1 \deq \frac{1}{2}k(k + 3) -1 \:.
  \end{equation*}
  Finally, one has $8 \le s_1 < s$, so \partprpref{r}{d} yields
  \begin{equation*}
    r \dle m-f_1 +f_0 \deq \upbeta + 1 -f_1 + f_0 \deq  -2t + f_0 \deq 0 \:.
  \end{equation*}
  By \prpref{q} and \eqref{pqr} it now follows that $R$ is Golod.
\end{prf*}

The next statement is folded in to \thmref{main} but worth recording
separately.

\begin{cor}
  \label{cor:level}
  Assume that $R$ is level. If $s \ge 8$ or $s$ is odd, then $R$ is
  Golod.
\end{cor}

\begin{prf*}
  For $s \ge 10$ and odd $s \ge 3$ it follows from \thmref{t2=t} that
  $R$ is Golod, and $s \ge 2$ holds by \eqref{ts}.  For $s=8$ apply
  \prpref{fam-even} with $k = 4$.
\end{prf*}

The final result of this section does not per se deal with the setup
in \stpref[]{3}, it only invokes the local ring $(Q,\mfq)$, but it
does come in handy in the proof of the main theorem. We have adapted
the proof and notation from Herzog, Reiner, and Welker's result
\thmcite[4]{HRW-99} on Golodness of componentwise linear ideals.

\begin{lem}
  \label{lem:S}
  Let $J \subseteq \mfq^2$ be a homogeneous $\mfq$-primary ideal in
  $Q$ and set $\sfA = \Tor[Q]{*}{Q/J}{\k}$. Let $u$ be the initial
  degree of $J$ and $\It{J}{u}$ the ideal generated by $J_u$; set
  $\sfB = \Tor[Q]{*}{Q/\It{J}{u}}{\k}$. For integers $j,k \le u$ and
  $\ell \le u+1$ there are inequalities
  \begin{align*}
    \tag{a}
    \rnk{(\sfA_{1\,j} \cdot \sfA_{1\,k})} %
    &\dle
      \rnk{(\sfB_{1\,j} \cdot \sfB_{1\,k})} \\
    \tag{b}
    \rnk{(\sfA_{1\,j} \cdot \sfA_{2\,\ell})} %
    &\dle
      \rnk{(\sfB_{1\,j} \cdot \sfB_{2\,\ell})} \\
    \tag{c}
    \rnk{(\sfA_{2\,\ell} \xra{\delta'\!} \Hom[\k]{\sfA_{1\,j}}{\sfA_{3\,j+\ell}})} %
    &\dle \rnk{(\sfB_{2\,\ell} \xra{\delta''\!} \Hom[\k]{\sfB_{1\,j}}{\sfB_{3\,j+\ell}})}
  \end{align*}
  with $\delta'$ and $\delta''$ defined as in \pgref{ms}.  In
  particular, if $Q/\It{J}{u}$ is Golod, then one has
  \begin{equation*}
    (\sfA_1)_{\le u} \cdot (\sfA_1)_{\le u} \deq 0 \deq (\sfA_1)_{\le u} \cdot (\sfA_2)_{\le u+1} \:.
  \end{equation*}
\end{lem}

\begin{prf*}
  Let $\Kzl$ be the Koszul complex on a minimal set of generators for
  $\mfq$. The surjective homomorphism $\mapdef{\pi}{Q/\It{J}{u}}{Q/J}$
  of graded rings is an isomorphism in degrees from $0$ to $u$. It
  induces a surjective morphism of differential graded algebras,
  \begin{equation*}
    \dmapdef[\lra]{\tp[Q]{\pi}{\Kzl}}{\tp[Q]{Q/\It{J}{u}}{\Kzl}}{\tp[Q]{Q/J}{\Kzl}} \:.
  \end{equation*}
  There are isomorphisms of graded-commutative $\k$-algebras
  $\sfA \is \H{\tp[Q]{Q/J}{\Kzl}}$ and
  $\sfB \is \H{\tp[Q]{Q/\It{J}{u}}{\Kzl}}$; see
  \cite[(1.2.1)]{LLA12}. Thus $\H{\tp[Q]{\pi}{\Kzl}}$ induces a
  morphism, $\mapdef{\tilde\pi}{\sfB}{\sfA}$, of graded-commutative
  algebras.  Since the differential on $\Kzl$ is linear and
  $\pi_{\le u}$ is an isomorphism, it follows that
  \begin{equation*}
    \dmapdef{\tilde\pi_{i\,j}}{\sfB_{i\,j}}{\sfA_{i\,j}}
  \end{equation*}
  is an isomorphism for all $i$ and $j \le i+u-1$.

  \proofoftag{a} For integers $j,k \le u$ and elements
  $\sfx\in\sfA_{1\,j}$ and $\sfx'\in\sfA_{1\,k}$ one has
  \begin{equation*}
    \sfx\sfx' \deq \tilde\pi_{1\,j}\tilde\pi_{1\,j}^{-1}(\sfx)\cdot
    \tilde\pi_{1\,k}\tilde\pi_{1\,k}^{-1}(\sfx')
    \deq \tilde\pi_{2\,j+k}\big(\tilde\pi_{1\,j}^{-1}(\sfx)\cdot\tilde\pi^{-1}_{1\,k}(\sfx')\big) \:.
  \end{equation*}
  Thus, if the product $\sfx\sfx'$ is non-zero, then so is
  $\tilde\pi^{-1}(\sfx)\cdot\tilde\pi^{-1}(\sfx')$ in
  $\sfB_{1\,j}\cdot \sfB_{1\,k}$. It follows that every nonzero
  element in $\sfA_{1\,j}\cdot \sfA_{1\,k}$ lifts to a nonzero element
  in $\sfB_{1\,j}\cdot \sfB_{1\,k}$ and, further, that linearly
  independent elements in $\sfA_{1\,j}\cdot \sfA_{1\,k}$ lift to
  linearly independent elements in $\sfB_{1\,j}\cdot \sfB_{1\,k}$.
  
  Parts (b) and (c) follow from parallel arguments, and the final
  assertion follows as $\sfB_1\cdot \sfB_1 = 0 = \sfB_1\cdot \sfB_2$
  holds if $\It{J}{u}$ is Golod.
\end{prf*}


\section{The main theorem}
\label{sec:7}

\noindent We have hitherto focused on a setup where we are given
compressed artinian Gorenstein rings defined by ideals $I_1$ and
$I_2$, and we have analyzed the situation where the ideal
$I_1 \cap I_2$ defines a compressed ring of type $2$. In the statement
of our main theorem below, the point of view is slightly shifted: The
focus is now on a compressed ring of type $2$ whose defining ideal is
assumed to be obtainable as an intersection of ideals that define
compressed Gorenstein ring; see also \rmkref{main}.

\begin{thm}
  \label{thm:main}
  Let $\k$ be a field, set $Q = \pows[\k]{x,y,z}$ and
  $\mfq = (x,y,z)$. Let $I \subseteq \mfq^2$ be a homogeneous ideal
  such that $R = Q/I$ is compressed artinian of type $2$. Assume that
  $I$ is the intersection of homogeneous ideals $I_1$ and $I_2$ that
  define compressed Gorenstein~rings.  Denote by $m$ and $t$ the
  minimal number of generators and the initial degree of $I$~and set
  $a = \min\setof{i \ge 0}{\mfq^iI_2 \subseteq I_1}$; let
  $\chi^{s_1} + \chi^s$ be the socle polynomial of $R$
  with~$s_1 \le s$.

  \subparagraph If $s$ is odd, then the following assertions hold.
  \begin{prt}
  \item I\fsp\ $\frac{s+1}{2} < t$, then $R$ is of class $\clH{0,0}$,
    i.e.\ Golod.
  \item I\fsp\ $\frac{s+1}{2} = t$ and $s \ge 5$, then $s_1 \ne s$ and
    $R$ is of class $\clG{r}$ with
    \begin{align*}
      r &\deq m - a(a+2)  \\
        & \dge \textstyle\frac{1}{2}(s+3 - a(a+1)) 
          \dge 1  \qtext{where} a=s_1 - \frac{s-1}{2} \:.
    \end{align*}
  \item I\fsp\ $\frac{s+1}{2} = t$ and $s = 3$, then $s_1 = 2$ and $R$
    is of one of the following classes
    \begin{flalign*}&\hspace{6.75pc}
      \begin{cases}
        \clB & \text{with } \ m = 5 \\
        \clG{3} & \text{with } \ m = 6 \:.
      \end{cases}&
    \end{flalign*}
  \end{prt}

  \noindent If $s$ is even, then the following assertions
  hold.
  \begin{prt}
    \setcounter{prt}{3}
  \item I\fsp\ $\frac{s}{2} +1 < t$, then $R$ is of class $\clH{0,0}$,
    i.e.\ Golod.
  \item I\fsp\ $\frac{s}{2}+1 = s_1 \ne s$, then $\frac{s}{2}+1 = t$
    and $R$ is of class $\clG{m-3}$.
  \item I\fsp\ $\frac{s}{2}+1 = t$ and $s_1 \ne s$, then $R$ is of
    class $\clH{0,0}$ or of class $\clG{r}$ with
    \begin{align*}
      \textstyle m - \frac{1}{2}a
      (a + 3)  
      \dge r
      & \dge m - \textstyle a(a+2)  \\
      & \dge s + 3 - \textstyle\frac{3}{2}a
        (a + \frac{5}{3})  \qtext{where} a = s_1 - \frac{s}{2} \:.
    \end{align*}

  \item I\fsp\ $\frac{s}{2} + 1 = t$ and $s_1 = s$, then $s \le 8$ and
    the following assertions hold.
    \begin{rqm}
    \item[] If $s=2$, then $R$ is of one of the following classes
      \begin{flalign*}&\hspace{6.75pc}
        \begin{cases}
          \clH{3,2} & \text{with } \ m = 4 \\
          \clB & \text{with } \ m=5 \:.
        \end{cases} &
      \end{flalign*}

    \item[] If $s=4$, then $R$ is of one of the following classes
      \begin{flalign*}&\hspace{6.75pc}
        \begin{cases}
          \clH{0,0} & \text{with } \ 5 \le m \le 8\\
          \clG{r} & \text{with } \ 6 \le m \le 7 \ \text{ and } \ r \le m-5\\
          \clH{0,2} & \text{with } \ m = 7\:.
        \end{cases} &
      \end{flalign*}
    \item[] If $s=6$, then $R$ is of one of the following classes
      \begin{flalign*}&\hspace{6.75pc}
        \begin{cases}
          \clH{0,0} & \text{with } \ 9 \le m \le 11\\
          \clG{1} & \text{with } \ m = 10 \:.
        \end{cases} &
      \end{flalign*}
    \item[] If $s=8$, then $R$ is of class $\clH{0,0}$ with $m=14$.
    \end{rqm}
  \end{prt}
\end{thm}

\begin{rmk}
  \label{rmk:main}
  \noindent For an ideal $I$ as in \thmref{main} one can always find
  homogeneous Gorenstein ideals $I_1$ and $I_2$ with
  $I_1 \cap I_2 = I$, but the assumption that they define compressed
  rings is crucial. Indeed, for the ideals $I_2$ and $I_3$ from
  \exaref{r1r2} the intersection $J = I_2 \cap I_3$ is a six-generated
  ideal of initial degree $3$, and it defines a compressed ring with
  socle polynomial $\chi^3 + \chi^4$. With the \emph{Macaulay2}
  package \cite{M2-LWCOVl} one can verify that $Q/J$ is of class
  $\clG{1}$.  Had $Q/I_3$ been compressed---or had it in any way been
  possible to obtain $J$ as the intersection of two homogeneous ideals
  that define compressed Gorenstein rings---then $Q/J$ would by
  \partthmref[]{main}{e} have been of class~$\clG{3}$. Thus, it
  follows that this ideal $J$ can not be obtained as an intersection
  of ideals that define compressed Gorenstein rings. It can happen,
  though, that an ideal $I$ as in \thmref[]{main} can be obtained as
  intersections of Gorenstein ideals,
  $I_1 \cap I_2 = I = I_2 \cap I_3$, in such a way that the rings
  $Q/I_2$ and $Q/I_3$ are compressed but $Q/I_1$ is not, see
  \exaref{r1r2a}.
\end{rmk}

\begin{bfhpg*}[Proof of \ref{thm:main}]
  Let $I_1$ and $I_2$ be homogeneous Gorenstein ideals with
  $I = I_1\cap I_2$ and such that $Q/I_1$ and $Q/I_2$ are compressed;
  this means that $I_1$ and $I_2$ fit \stpref{3} with $e=3$. Further,
  the assumption that $R$ is compressed means that the results in
  \secref{parameters} apply. They take care of most of the proof, but
  certain special cases elude them. To deal with those cases, we
  enlist a result of Bigatti, Geramita, and Migliore~\cite{BGM-94} on
  growth of Hilbert functions as well as our joint work with Weyman
  \cite{CVW-20}, which is based on linkage theory.
  
  \subparagraph\textbf{\textit{Parts} (a) \textit{and} (d):} By
  \eqref{compressed-2} one has
  $\left\lceil \frac{s+1}{2}\right\rceil \le t$, and if strict
  inequality holds, then $R$ is Golod by \pgref{golod}, which
  precisely means that $R$ is of class $\clH{0,0}$; see \pgref{ms}.

  \subparagraph\textbf{\textit{Part} (b):} Let $s$ be odd and assume
  that $\frac{s+1}{2} = t$ and $s \ge 5$ hold. By \prpref{ta} one has
  $s_1 < s$ and \eqref{ts} yields $s_1 \ge 3$, so \partprpref{r}{c},
  \thmref{resR}, and \lemref{ta} yield
  \begin{equation*}
    r \deq m-f_1 \dge t+1 -f_0 \:>\: 0 \:.
  \end{equation*}
  Now it follows from \prpref{q} and \eqref{pqr} that $R$ is of class
  $\clG{r}$. By \lemref{resR} one has $a = s_1 - \frac{s-1}{2}$, and
  direct computations based on \eqref{f} yield
  \begin{equation*}
    m - f_1 =  m-a(a+2) \qqand
    t+1 -f_0 \deq  \textstyle\frac{1}{2}(s+3 - a(a+1))\:.
  \end{equation*}

  \subparagraph\textbf{\textit{Part} (c):} Assume that $s = 3$ and
  $t= 2$ hold. By \prpref{ta} and \eqref{ts} one has $s_1=2$, so
  \lemref{resR} yields $a=1$. Per \eqref{f} one now has $f_0 = 1$ and
  $f_1=3$. It follows from \thmref{resR} that the ideal $I$ is
  minimally generated by two quadratic and $\upbeta+3$ cubic
  forms. Our first step is to prove that $\upbeta$ is at most $1$.

  The $h$-vector of $R$ is $(1,3,4,1)$, see Table~\ref{table0}. Let
  $\It{I}{2}$ be the ideal generated by the two quadratic forms. As
  \begin{equation*}
    \textstyle
    h_{Q/\It{I}{2}}(2) \deq 4 \deq \binom{3}{2} + \binom{1}{1} \qtext{one has}
    h_{Q/\It{I}{2}}(3) \dle 4^{\langle 2 \rangle} \deq \binom{4}{3} + \binom{2}{2} \deq 5
  \end{equation*}
  by Macaulay's theorem, see \thmcite[4.2.10]{bruher}. Thus
  $\upbeta + 3 \le 5-1$ and hence $\upbeta \le 1$.

  If $\upbeta=0$ holds, then $m$ is $5$ and \partprpref{r}{c} yields
  $r = 2$. If $R$ were of class $\clH{p,q}$, then one would have
  $q \le 1$ by \thmcite[1.1]{CVW-20}, which is impossible as $r=q$
  holds by \eqref{pqr}. Since the type of $R$ is $2$, it now follows
  from \thmcite[4.5(b)]{CVW-20} that $R$ is of class $\clB$.

  If $\upbeta=1$ holds, then $m$ is $6$ and \partprpref{r}{c} yields
  $r =3$; since the type of $R$ is $2$, it is of class $\clG{3}$; see
  \eqref{pqr}.

  \subparagraph\textbf{\textit{Part} (e):} Let $s$ be even and assume
  that one has $\frac{s}{2} + 1 = s_1 \ne s$. \lemref[Lemmas~]{ta} and
  \lemref[]{resR} yield $\frac{s}{2} + 1 = t$, $a=1$, and
  $m \ge 2t = s+2$. By \partprpref{r}{e} one has $r = m-3$, and as
  $s \ge 4$ holds, $r$ is at least $3$, whence $R$ is of class
  $\clG{r}$, see \eqref{pqr}.
  
  \subparagraph\textbf{\textit{Part} (f):} Let $s$ be even and assume
  that $\frac{s}{2} + 1 =t$ and $s_1 < s$ hold. It follows that $s$ is
  at least $4$, so $s_1$ is at least $3$ by \eqref{ts}. From
  \prpref[Propositions~]{q} and \partprpref[]{r}{a,d} it follows that
  $R$ is Golod or of class $\clG{r}$ with
  \begin{equation*}
    m-f_1 + f_0 \dge r \dge m-f_1 \:.
  \end{equation*}
  By \lemref{resR} one has $a = s_1 - \frac{s}{2}$ and direct
  computations based on \eqref{f} yield
  \begin{align*}
    m - f_1 + f_0 %
    &\deq m - f_2 + 1 \deq m - \textstyle\frac{1}{2}
      (a+2)(a+1) + 1 \deq m - \textstyle\frac{1}{2}a(a+3)\\[-1.5ex]
    \intertext{and}\\[-6ex]
    m-f_1 & \deq m - a(a+2) \:.
  \end{align*}
  By \thmref{resR} one has $m - f_1 \ge 2t+1-f_0 -f_1$, and another
  computation yields
  \begin{align*}
    2t+1 - f_0  -f_1 &\deq s+3 - \textstyle\frac{1}{2}a(a+1) - a(a+2)
                       \deq s + 3 - \textstyle\frac{3}{2}
                       a(a+\frac{5}{3}) \:.
  \end{align*}

  \subparagraph\textbf{\textit{Part} (g):} Let $s$ be even and assume
  that $\frac{s}{2} + 1 = t$ and $s_1 = s$ hold. It follows from
  \prpref{ta} that $s$ is at most $8$, and by \lemref{resR} one has
  $a = \frac{s}{2}$.  In the balance of the proof, let $\sfA$ be the
  bigraded $\k$-algebra $\Tor[Q]{*}{R}{\k}$.

  The case $s=8$ is covered by \corref{level}. We address the
  remaining cases in descending order.

  \begin{center}
    \it Case $s=6$
  \end{center}
  One has $t=4$, $a=3$, $f_0=6$, $f_1=15$, and $f_2= 10$; see
  \eqref{a} and \eqref{f}.  Recall from \thmref{resR} that for some
  $\upbeta \ge 15 - 8 - 1 = 6$ the minimal graded free resolution of
  $R$ over $Q$ has the form
  \begin{equation*}
    \tag{$\ast$}
    Q \lla
    \begin{matrix}
      Q^{3}(-4) \\ \oplus \\ Q^{\upbeta}(-5)
    \end{matrix}
    \lla
    \begin{matrix}
      Q^{\upbeta-6}(-5) \\ \oplus \\ Q^{10}(-6)
    \end{matrix}
    \lla
    \begin{matrix}
      Q^2(-9)
    \end{matrix}
    \lla 0 \:.
  \end{equation*}
  Thus the ideal $I$ is minimally generated by three quartic forms and
  $\upbeta$ quintics. Our first step is to prove that $\upbeta$ is at
  most $8$.

  The $h$-vector of $R$ is $(1,3,6,10,12,6,2)$, see
  Table~\ref{table0}. Let $\It{I}{4}$ be the ideal generated by the
  three quartic forms. As
  \begin{equation*}
    \textstyle
    h_{Q/\It{I}{4}}(4) \deq 12 \deq \binom{5}{4} + \binom{4}{3} + \binom{3}{2} \qtext{one has}
    h_{Q/\It{I}{4}}(5) \dle 12^{\langle 4 \rangle} \deq \binom{6}{5} + \binom{5}{4}
    + \binom{4}{3} \deq 15
  \end{equation*}
  by Macaulay's theorem; thus $\upbeta \le 15-6=9$ holds. Assume
  towards a contradiction that $\upbeta$ is $9$, that is,
  $h_{Q/\It{I}{4}}(5) = 15$. This assumption implies that the Hilbert
  function of $Q/\It{I}{4}$ has maximal growth in degree $4$. It
  follows from \prpcite[2.7]{BGM-94} that the generators of
  $\It{I}{4}$ have a common cubic factor $f$, that is,
  $\It{I}{4} = \mfq f$. Since $\It{I}{4}$ is contained in $I$, the
  element $f+I$ is a socle element in $R$. As $f + I$ has degree $3$
  this contradicts the assumption that $R$ is level of socle degree
  $6$.  Thus one has $6 \le \upbeta \le 8$ and, therefore,
  $9 \le m \le 11$.

  By \lemref{p} one has $p=0$, so $R$ is of class $\clH{0,0}$,
  $\clH{0,2}$, or $\clG{r}$; see \pgref{ms}. Consider the bigraded
  $\k$-algebra $\sfA$. For homogeneous elements $\sfx \in \sfA_1$ and
  $\sfy \in \sfA_2$ with $\sfx\sfy \ne 0$ in $\sfA_3$ it follows from
  $(\ast)$ that the internal degrees are $|\sfx|=4$ and $|\sfy|=5$;
  that is, $\sfA_1\cdot\sfA_2 = \sfA_{1\,4}\cdot \sfA_{2\,5}$. Per
  $(\ast)$ one has $\rnk{\sfA_{1\,4}} = 3$ and
  $\rnk{\sfA_{2\,5}} = \upbeta -6 = m-9$.  It follows that $r$ is at
  most $m-9$, so $m$ is at least $10$ for rings of class $\clG{r}$. We
  proceed to prove that $R$ is Golod for $m=11$; this rules out the
  possibility $\clH{0,2}$, and it means that $R$ can be of class
  $\clG{r}$ only for $m=10$ and $r=1$.

  Assume that $m=11$ holds, i.e.\ $\upbeta = 8$. It suffices to show
  that the minimal free resolution of $Q/\It{I}{4}$ has the form
  $Q \lla Q^{3}(-4) \lla Q^{2}(-5) \lla 0$. It then follows that
  $Q/\It{I}{4}$ is a Golod ring, see for example \cite[5.3.4]{ifr}. As
  established above, one has
  $\sfA_1\cdot\sfA_2 = \sfA_{1\,4}\cdot \sfA_{2\,5}$, and \lemref{S}
  yields $\sfA_{1\,4}\cdot \sfA_{2\,5}=0$. Thus, $q=0$ holds and $R$
  is of class $\clH{0,0}$. To establish that $Q/\It{I}{4}$ has the
  asserted free resolution, notice first from $(\ast)$ that the three
  quartic generators of $\It{I}{4}$ have $\upbeta - 6 = 2$ linear
  syzygies. It suffices to show that they have no further syzygies,
  and since the Koszul relations are of degree $8$ this comes down to
  verifying that the minimal free resolution of $Q/\It{I}{4}$ has no
  $Q(-u)$ summand in degree $2$ for $6 \le u \le 8$. As $\upbeta = 8$
  holds, one has
  \begin{equation*}
    h_{Q/\It{I}{4}}(5)  \deq h_R(5) + \upbeta \deq 6 + 8 \deq 14
    \deq
    \textstyle
    \binom{6}{5} + \binom{5}{4} 
    + \binom{3}{3} + \binom{2}{2} + \binom{1}{1} \:,
  \end{equation*}
  so Macaulay's theorem yields
  \begin{equation*}
    \textstyle
    h_{Q/\It{I}{4}}(6) \dle \binom{7}{6} + \binom{6}{5} + 
    \binom{4}{4} + \binom{3}{3} + \binom{2}{2} \deq 16 \:.
  \end{equation*}
  A straightforward calculation based on \eqref{BS} and \eqref{HS} now
  yields
  \begin{equation*}
     \beta_{26}(Q/\It{I}{4}) - \beta_{36}(Q/\It{I}{4}) \deq b_{Q/\It{I}{4}}(6) \deq  
    h_{Q/\It{I}{4}}(6) - 16 \dle 0 \:.
  \end{equation*}
  A linear relation between the two first syzygies of the three
  quartics would also show in $(\ast)$, so one has
  $\beta_{36}(Q/\It{I}{4}) = \beta_{36}(R) = 0$, this forces
  $\beta_{26}(Q/\It{I}{4})=0$ and $h_{Q/\It{I}{4}}(6) = 16$. Further,
  this implies $\beta_{37}(Q/\It{I}{4}) = \beta_{37}(R) = 0$ as a
  quadratic relation between the first syzygies of the three quartics
  would also show in $(\ast)$. Repeating the procedure, one gets
  $h_{Q/\It{I}{4}}(7) \le 18$ and
  \begin{equation*}
     \beta_{27}(Q/\It{I}{4}) -
    \beta_{37}(Q/\It{I}{4}) \deq b_{Q/\It{I}{4}}(7) \deq h_{Q/\It{I}{4}}(7) - 18 \dle 0 \:.
  \end{equation*}
  As above this forces $\beta_{27}(Q/\It{I}{4})=0$ and
  $h_{Q/\It{I}{4}}(7) = 18$ and, in addition, one gets
  $\beta_{38}(Q/\It{I}{4}) = \beta_{38}(R) = 0$ as a cubic relation
  between the first syzygies of the three quartics would also show in
  $(\ast)$. Repeating the procedure a third time yields
  $h_{Q/\It{I}{4}}(8) \le 20$ and
  \begin{equation*}
    \beta_{28}(Q/\It{I}{4}) -
    \beta_{38}(Q/\It{I}{4}) \deq b_{Q/\It{I}{4}}(8) \deq h_{Q/\It{I}{4}}(8) - 20 \dle 0 \:,
  \end{equation*}
  which this forces $\beta_{28}(Q/\It{I}{4})=0$ as desired.

  \begin{center}
    \it Case $s=4$
  \end{center}
  One has $t=3$, $a=2$, $f_0=3$, $f_1=8$, and $f_2=6$; see \eqref{a}
  and \eqref{f}.  Recall from \thmref{resR} that for some
  $\upbeta \ge 8-6-1=1$ the minimal graded free resolution of $R$ over
  $Q$ has the form
  \begin{equation*}
    \tag{$\ast\ast$}
    Q \lla
    \begin{matrix}
      Q^{4}(-3) \\ \oplus \\ Q^{\upbeta}(-4)
    \end{matrix}
    \lla
    \begin{matrix}
      Q^{\upbeta-1}(-4) \\ \oplus \\ Q^{6}(-5)
    \end{matrix}
    \lla
    \begin{matrix}
      Q^2(-7)
    \end{matrix}
    \lla 0 \:.
  \end{equation*}
  Thus the ideal $I$ is minimally generated by four cubic forms and
  $\upbeta$ quartics. Our first step is to show that $\upbeta$ is at
  most $4$.
  
  The $h$-vector of $R$ is $(1,3,6,6,2)$, see Table~\ref{table0}. Let
  $\It{I}{3}$ denote the ideal generated by the four cubics. As
  \begin{equation*}
    \textstyle
    h_{Q/\It{I}{3}}(3) \deq 6 \deq \binom{4}{3} + \binom{2}{2} + \binom{1}{1} \qtext{one has}
    h_{Q/\It{I}{3}}(4) \dle 6^{\langle 3 \rangle} \deq \binom{5}{4} + \binom{3}{3} + \binom{2}{2} \deq 7
  \end{equation*}
  by Macaulay's theorem; thus $\upbeta \le 7-2=5$ holds. Assume
  towards a contradiction that $\upbeta$ is $5$, i.e.\
  $h_{Q/\It{I}{3}}(4) = 7$.  The assumption $h_{Q/\It{I}{3}}(4) = 7$
  implies that the Hilbert function of $Q/\It{I}{3}$ has maximal
  growth in degree $3$. It follows from \prpcite[2.7]{BGM-94} that the
  generators of $\It{I}{3}$ have a common linear factor, i.e.\ the
  four cubics have the form $l q_1,\ldots, l q_4$ for some linear form
  $l$ and quadratic forms $q_1,\ldots,q_4$. The exact sequence of
  graded $Q$-modules,
  \begin{equation*}
    0 \lra Q(-1)/(I_2:l) \lra Q/I_2 \lra Q/(I_2 + (l)) \lra 0 \:,
  \end{equation*}
  yields $h_{Q/I_2}(i) = h_{Q/(I_2:l)}(i-1) + h_{Q/(I_2 + (l))}(i)$
  for all $i$. As $I$ is contained in $I_2$, the quadratic forms
  $q_1,\ldots,q_4$ belong to $(I_2:l)$, so one has
  $h_{Q/(I_2:l)}(2) \le 2$. The ideal $(I_2:l)$ defines a Gorenstein
  ring of socle degree $s-1=3$, see \cite[Ch.~IV Thm.~35]{ZS-1} and
  \corcite[I.2.4]{pms}, so the possible $h$-vectors of this ring are
  $(1,1,1,1)$ and $(1,2,2,1)$. As one has $h_{Q/I_2} = (1,3,6,3,1)$,
  see \pgref{compressed-1}, the $h$-vector of $Q/(I_2 + (l))$ would
  have to be $(1,2,5,2)$ or $(1,2,4,1)$, but by Macaulay's theorem
  neither is a possible $h$-vector; a contradiction.  Thus one has
  $1 \le \upbeta \le 4$ and $5 \le m \le 8$.

  By \lemref{p} one has $p=0$, so $R$ is of class $\clH{0,0}$,
  $\clH{0,2}$, or $\clG{r}$, see \pgref{ms}.  Now consider the
  bigraded $\k$-algebra $\sfA$. For homogeneous elements
  $\sfx \in \sfA_1$ and $\sfy \in \sfA_2$ with $\sfx\sfy \ne 0$ in
  $\sfA_3$ it follows from $(\ast\ast)$ that the internal degrees are
  $|\sfx|=3$ and $|\sfy|=4$; that is,
  $\sfA_1\cdot\sfA_2 = \sfA_{1\,3}\cdot \sfA_{2\,4}$. Per $(\ast\ast)$
  one has $\rnk{\sfA_{1\,3}} = 4$ and
  $\rnk{\sfA_{2\,4}} = \upbeta -1 = m-5$.  It follows that $r$ is at
  most $m-5$, so $m$ is at least $6$ for rings of class $\clG{r}$ and
  at least $7$ for rings of class $\clH{0,2}$. It remains to prove
  that $R$ is Golod for $m=8$.
  
  Assume that $m=8$ holds, i.e.\ $\upbeta = 4$. It suffices to show
  that the minimal free resolution of $Q/\It{I}{3}$ has the form
  $Q \lla Q^{4}(-3) \lla Q^{3}(-4) \lla 0$. It then follows that
  $Q/\It{I}{3}$ is a Golod ring, see for example \prpcite[5.3.4]{ifr}. As
  established above, one has
  $\sfA_1\cdot\sfA_2 = \sfA_{1\,3}\cdot \sfA_{2\,4}$, and \lemref{S}
  yields $\sfA_{1\,3}\cdot \sfA_{2\,4}=0$. Thus, $q=0$ holds and $R$
  is of class $\clH{0,0}$.  To establish that $Q/\It{I}{3}$ has the
  asserted free resolution, notice first from $(\ast\ast)$ that the
  four cubic generators of $\It{I}{3}$ have $\upbeta - 1 = 3$ linear
  syzygies. It suffices to show that they have no further syzygies,
  which comes down to verifying that the minimal free resolution of
  $Q/\It{I}{3}$ has no $Q(-5)$ or $Q(-6)$ summand in degree $2$. As
  $\upbeta = 4$ holds, one has
  \begin{equation*}
    \textstyle
    h_{Q/\It{I}{3}}(4) \deq h_R(4) + \upbeta \deq 2 + 4 \deq 6 \deq
    \binom{5}{4} + \binom{3}{3}
  \end{equation*}
  and, therefore,
  $h_{Q/\It{I}{3}}(5) \le \binom{6}{5} + \binom{4}{4} = 7$ by
  Macaulay's theorem. A straightforward calculation based on
  \eqref{BS} and \eqref{HS} now yields
  \begin{equation*}
     \beta_{25}(Q/\It{I}{3}) - \beta_{35}(Q/\It{I}{3}) \deq  
    b_{Q/\It{I}{3}}(5) \deq h_{Q/\It{I}{3}}(5) - 6 \dle 1 \:.
  \end{equation*}
  A linear relation between the first syzygies of the four cubics
  would also show in $(\ast\ast)$, so one has
  $\beta_{35}(Q/\It{I}{3}) = \beta_{35}(R) = 0$, which implies that
  $\beta_{25}(Q/\It{I}{3}) \le 1$ holds. Equality would force
  $h_{Q/\It{I}{3}}(5) = 7$ and $h_{Q/\It{I}{3}}(6) \le 8$ which
  would yield
  \begin{equation*}
    \beta_{26}(Q/\It{I}{3}) - \beta_{36}(Q/\It{I}{3}) \deq
    b_{Q/\It{I}{3}}(6) \deq
    h_{Q/\It{I}{3}}(6) - 9 \,<\, 0 \:;
  \end{equation*}
  this is absurd, as one has
  $\beta_{36}(Q/\It{I}{3}) = \beta_{36}(R) =0$ since a relation in
  degree $6$ between the first syzygies of the four cubics would also
  show in $(\ast\ast)$. Thus one has $\beta_{25}(Q/\It{I}{3}) = 0$ and
  $h_{Q/\It{I}{3}}(5) = 6 = \binom{6}{5}$. Macaulay's theorem now
  yields $h_{Q/\It{I}{3}}(6) \le 6^{\langle 5 \rangle} = 7$, and as above one gets
  \begin{equation}
    \tag{$\dagger$}
    \beta_{26}(Q/\It{I}{3}) - \beta_{36}(Q/\It{I}{3}) \deq  
    b_{Q/\It{I}{3}}(6) \deq h_{Q/\It{I}{3}}(6) - 6 \dle 1 \:.
  \end{equation}
  Since $\beta_{36}(Q/\It{I}{3}) = 0$ this implies
  $\beta_{26}(Q/\It{I}{3}) \le 1$.  Assume towards a contradiction
  that equality holds. This implies $h_{Q/\It{I}{3}}(6) = 7$, i.e.\
  the Hilbert function of $Q/\It{I}{3}$ has maximal growth in degree
  $5$. It now follows from \prpcite[2.7]{BGM-94} that the generators
  of $\It{(\It{I}{3})}{5}$ have a common linear factor $l$. That is,
  one has $\mfq^2\It{I}{3} \subseteq (l)$. As $(l)$ is a prime ideal,
  one has $(l):\mfq = (l)$ and, therefore,
  $(l):\mfq^2 = ((l):\mfq):\mfq = (l)$. It follows that $\It{I}{3}$ is
  contained in $(l)$, so the four cubic generators have a common
  linear factor; as in the subcase $\upbeta = 4$ above this leads to a
  contradiction. Thus $h_{Q/\It{I}{3}}(6) \le 6$ holds, wheence
  $(\dagger)$ yields
  $\beta_{26}(Q/\It{I}{3}) - \beta_{36}(Q/\It{I}{3}) \le 0$ and as
  $\beta_{36}(Q/\It{I}{3})$ is $0$ this implies
  $\beta_{26}(Q/\It{I}{3}) = 0$ as desired.

  \begin{center}
    \it Case $s=2$
  \end{center}
  One has $t=2$, $a=1=f_0$, and $f_1=3 = f_2$; see \eqref{a}
  and~\eqref{f}.  Recall from \thmref{resR} that for some
  $\upbeta \ge 0$ the minimal graded free resolution of $R$ over $Q$
  has the form
  \begin{equation*}
    \tag{$\ast\ast\ast$}
    Q \lla
    \begin{matrix}
      Q^{4}(-2) \\ \oplus \\ Q^{\upbeta}(-3)
    \end{matrix}
    \lla
    \begin{matrix}
      Q^{ 2 + \upbeta}(-3) \\ \oplus \\ Q^{3}(-4)
    \end{matrix}
    \lla
    \begin{matrix}
      Q^2(-5)
    \end{matrix}
    \lla 0 \:.
  \end{equation*}
  Thus the ideal $I$ is minimally generated by four quadratic forms
  and $\upbeta$ cubics.  \partprpref{r}{a} yields
  $r \ge m-3 = \upbeta + 1$. Our first step is to show that
  $\upbeta \le 1$.

  The $h$-vector of $R$ is $(1,3,2)$, see Table~\ref{table0}. Let
  $\It{I}{2}$ denote the ideal generated by the four quadratics. As
  \begin{equation*}
    \textstyle
    h_{Q/\It{I}{2}}(2) \deq 2 \deq \binom{2}{2} + \binom{1}{1} \qtext{one has}
    h_{Q/\It{I}{2}}(3) \dle 2^{\langle 2 \rangle} \deq \binom{3}{3} + \binom{2}{2} \deq 2
  \end{equation*}
  by Macaulay's theorem; thus $\upbeta \le 2$ holds. Towards a
  contradiction assume $\upbeta =2$. In $(\ast\ast\ast)$ there are now
  $4$ syzygies of degree $3$. They must be linear syzygies of the four
  quadratic generators, so they appear in the minimal free resolution
  of $Q/\It{I}{2}$ which, therefore, must have length $3$. Notice also
  that $\It{I}{2}$ is not a Gorenstein ideal, since it is minimally
  generated by an even number of generators, see
  \thmcite[2.1]{DABDEs77}. Now consider the bigraded $\k$-algebra
  $\sfB = \Tor[Q]{*}{Q/\It{I}{2}}{\k}$. For homogeneous elements
  $\sfx\in\sfA_1$ and $\sfy\in\sfA_2$ with $\sfx\sfy \ne 0$ in
  $\sfA_3$ it follows from $(\ast\ast\ast)$ that their internal
  degrees are $|\sfx|=2$ and $|\sfy|=3$. Thus one has
  \begin{equation*}
    \rnk{(\sfA_{2\,3} \to \Hom[\k]{\sfA_{1\,2}}{\sfA_{3\,5}})} \deq r \dge 3 \:,
  \end{equation*}
  so by \partlemref{S}{c} the map
  $\mapdef{\delta}{\sfB_2}{\Hom[\k]{\sfB_1}{\sfB_3}}$, cf.~\pgref{ms},
  has rank at least $3$. This contradicts \thmcite[3.1]{LLA12} which
  bounds the rank of $\delta$ above by $2$.

  If $\upbeta=0$ holds, then $m$ is $4$, so $R$ is an almost complete
  intersection and hence of class $\clH{3,2}$; see \cite[3.4.2]{LLA12}
  or \thmcite[4.1]{CVW-20}.

  If $\upbeta=1$ holds, then $m$ is $5$, and \partprpref{r}{c} yields
  $r = 2$. As $R$ has type $2$ it follows from \thmcite[4.5]{CVW-20}
  that $R$ is of class $\clB$ or $\clH{p,q}$. The latter option is
  ruled out by \thmcite[1.1]{CVW-20} and \eqref{pqr}, which yield
  $q=1$ and $q=r$. \qed
\end{bfhpg*}



The statement of \thmref{main} is organized according to: first the
parity of $s$, second the relation of
$\left\lceil \frac{s+1}{2}\right\rceil$ to $t$,
cf.~\eqref{compressed-2}, and finally a comparison of $m$ to~$r$. In
\corref[Corollaries~]{mixed}--\corref[]{even} below the conclusions
from \thmref[]{main} are summarized with an emphasis on the difference
between the degrees of the socle generators of $R$.

\begin{cor}
  \label{cor:mixed}
  Let $R$ be as in \thmref{main} with $s\le 4$. One has
  $2 \le s_1 \le 4$ and the next assertions hold.
  \begin{prt}
  \item If $s_1 =2$, then $s \le 3$ and $R$ is of class
    \begin{flalign*}&\hspace{9.85pc}
      \begin{cases}
        \clH{3,2} & \text{ with } s = 2 \\
        \clB & \text{ with } 2 \le s \le 3 \\
        \clG{3} & \text{ with } s = 3
      \end{cases}&
    \end{flalign*}
  \item If $s_1 =3$, then $R$ is of class
    \begin{flalign*}&\hspace{9.85pc}
      \begin{cases}
        \clH{0,0} & \text{ with } s = 3 \\
        \clG{r} & \text{ with } s = 4
      \end{cases}&
    \end{flalign*}
  \item If $s_1 =4$, then $s=4$ and $R$ is of class $\clH{0,0}$,
    $\clG{r}$, or $\clH{0,2}$.
  \end{prt}
\end{cor}

\begin{prf*}
  Under the the assumptions in \thmref{main}, \thmref{R1R2} applies;
  in particular inequalities $2 \le s_1 \le 4$ hold by \eqref{ts}.

  \proofoftag{a} If $s_1=2$, then \eqref{ts} yields $s \le 3$. If
  $s=2$, then $h_R=(1,3,2)$ holds, and for $s=3$ one has
  $h_R=(1,3,4,1)$, see Table~\ref{table0}; in either case
  $t=2=\left\lceil \frac{s+1}{2}\right\rceil$ holds, so the assertions
  follow from \partthmref[]{main}{c,g}.

  \proofoftag{b} If $s_1=3$, then \eqref{ts} yields $3 \le s \le
  4$. If $s=3$, then $h_R=(1,3,6,2)$ holds, see Table~\ref{table0}, so
  $t = 3 > \left\lceil \frac{s+1}{2}\right\rceil$ holds; for $s=4$ one
  gets $h_R=(1,3,6,4,1)$, so
  $t=3=\left\lceil \frac{s+1}{2}\right\rceil$ holds. The assertions
  now follow from \partthmref[]{main}{a,e}.

  \proofoftag{c} If $s_1=4$, then \eqref{ts} yields $s = 4$, and the
  assertion summarizes the case $s=4$ in \partthmref[]{main}{g}.
\end{prf*}

\rmkref{A} translates \corref[Corollaries~]{odd} and \corref[]{even}
into the summary given in the introduction.

\begin{cor}
  \label{cor:odd}
  Let $R$ be as in \thmref{main} with $s$ odd and $s \ge 5$. Set
  \begin{equation*}
    N(s) \deq \frac{s-2+\sqrt{4s+13}}{2} \:.
  \end{equation*}
  There are inequalities
  \begin{equation*}
    \textstyle\frac{s+1}{2} \dle s_1 \dle s
    \qqand \textstyle\frac{s+1}{2} \: < \: N(s) \: < \: s \:,
  \end{equation*}
  and the next assertions hold.
  \begin{prt}
  \item If one has $s_1 < N(s)$, then $R$ is of class $\clG{r}$.
  \item If one has $N(s) \le s_1$, then $R$ is of class $\clH{0,0}$,
    i.e.\ Golod.
  \end{prt}
\end{cor}

\begin{prf*}
  The first set of inequalities comes from \eqref{ts}; as $s$ is at
  least $5$, the second set follows immediately from the definition of
  $N(s)$. It follows from \prpref{ta} that the inequality
  $\frac{s+1}{2} < t$ holds if and only if one has $N(s) \le s_1$, and
  the two assertions now follow immediately from
  \partthmref{main}{a,b}.
\end{prf*}

The next corollary describes the situation for even socle degree
$s \ge 6$. In the odd case, \corref{odd}, a single bound on $s_1$
determines whether $R$ is Golod or of class $\clG{r}$. In
\corref{even} it takes two bounds, $N_1(s) < N_2(s)$, to definitively
separate the classes Golod and $\clG{r}$; in the intermediate interval
both possibilities occur, see Table~\ref{table1} but also
\pgref{ss:3}(d).

\begin{cor}
  \label{cor:even}
  Let $R$ be as in \thmref{main} with $s$ even and $s \ge 6$.  Set
  \begin{equation*}
    N_1(s) \deq \frac{3s-5+\sqrt{24s + 97}}{6} \qqand
    N_2(s) \deq \frac{s-1+\sqrt{8s+25}}{2} \:.
  \end{equation*}
  There are inequalities
  \begin{equation*}
    \textstyle\frac{s}{2} + 1 \dle s_1 \dle s
    \qqand N_1(s) \:<\: N_2(s)
    \:,
  \end{equation*}
  and the next assertions hold.
  \begin{prt}
  \item If one has $s_1 < N_1(s)$, then $R$ is of class $\clG{r}$.
  \item If one has $N_2(s) \le s_1$, then $R$ is of class $\clH{0,0}$,
    i.e.\ Golod.
  \item If one has $N_1(s) \le s_1 < N_2(s)$, then $R$ is of class
    $\clH{0,0}$ or $\clG{r}$.
  \end{prt}
\end{cor}

\begin{prf*}
  The first set of inequalities comes from \eqref{ts}, and it is
  immediate from the definitions that $N_1(s) < N_2(s)$ holds. It
  follows from \prpref{ta} that the inequality $\frac{s}{2}+1 < t$
  holds if and only if one has $N_2(s) \le s_1$, so part (b) holds by
  \partthmref{main}{d}. Further, $\frac{s}{2}+1 = t$ holds for
  $s_1 < N_2(s)$, in particular for $s_1 < N_1(s)$, and in that case
  \partthmref{main}{f} yields $a = s_1 - \frac{s}{2}$ and
  \begin{equation*}
    r \dge s + 3 - \textstyle\frac{3}{2}a(a + \frac{5}{3}) 
    \deq \textstyle -\frac{3}{2}a^2
    -\frac{5}{2}a + s + 3 \:.
  \end{equation*}
  This quadratic expression is positive for
  $a < \frac{-5+\sqrt{24s + 97}}{6}$, i.e.\ for $s_1 < N_1(s)$. This
  proves (a), and (c) follows from \partthmref{main}{f,g}.
\end{prf*}

\begin{rmk}
  \label{rmk:comp}
  Let $(Q,\mfq) $ be as in \thmref{main} and $J \subseteq \mfq^2$ be a
  homogeneous $\mfq$-primary ideal such that $S = Q/J$ has type
  $2$. With $\tilde{m}$ denoting the minimal number of generators of
  $J$, it is known---see for example \cite[Thm.~1.1 and
  Sect.~4]{CVW-20}---that $S$ is of one of the following classes
  \begin{equation*}
    \begin{cases}
      \clH{3,2} & \text{with } \ \tilde{m} = 4 \\
      \clH{0,0} & \text{with } \ \tilde{m} \ge 5 \\
      \clB & \text{with } \ \tilde{m} \ge 5 \ \text{and $\tilde{m}$ odd}\\
      \clG{r} & \text{with } \ \tilde{m} \ge 6 \ \text{and } r \le \tilde{m}-3 \\
      \clH{1,2} & \text{with } \ \tilde{m} \ge 6  \ \text{and $\tilde{m}$ even} \\
      \clH{0,2} & \text{with } \ \tilde{m} \ \ge 7 \;.
    \end{cases}
  \end{equation*}
  Empirical evidence suggests that essentially all of these classes
  materialize, cf.~\cite[Conj.~7.4]{CVW-20}. In comparison, the extra
  assumptions imposed in \thmref{main} restrict rings $R$ as in
  \thmref[]{main} with $m \ge 8$ to the classes $\clH{0,0}$ and
  $\clG{r}$.
\end{rmk}


We close this section with a result that provides further evidence for
\cite[Conjecture 7.4(a)]{CVW-20} and suggests that the answer to
\cite[Question 9.7]{KVn-1} is negative.

\begin{prp}
  \label{prp:conj}
  Let $R$ be as in \thmref{main}. If $R$ is of class $\clG{r}$, then
  $r \le m -3$.
\end{prp}

\begin{prf*}
  As $R$ is of class $\clG{r}$, it follows from \eqref{compressed-2}
  and \pgref{golod} that $\left\lceil \frac{s+1}{2}\right\rceil = t$
  holds. If $s$ is odd, then it follows from
  \corref[Corollaries~]{mixed} and \corref[]{odd} that $s_1 < s$
  holds. Thus one has $r = m -f_1$ by \partprpref{r}{c}, and
  $f_1 \ge 3$ holds by \eqref{f} as $a$ is positive by
  \partprpref{soc-val}{e}. Now assume that $s$ is even. If $s_1 = s$
  holds, then \partthmref{main}{g} yields $s \le 6$ and shows that
  $r \le m-5$ holds. Assuming now that $s_1 < s$ holds,
  \partprpref{r}{d} yields $r \le m -f_1 + f_0$. For $a \ge 2$ one has
  $f_1-f_0 \ge 3$, and $a=1$ implies $s_1 = \frac{s}{2}+1$, so $r=m-3$
  holds by \partthmref[]{main}{e}.
\end{prf*}

\section{Generic behavior} 
\label{sec:generic}

\noindent
In this final section we first elaborate on the remarks made in the
introduction about generic algebras being compressed. The statement of
\thmref{main} was informed by experiments, and we share some of the
collected data in Table~\ref{table1}. This data suggests a number of
questions; we address a few of them.  Finally, we discuss in which
sense \thmref{main} explains the class of a randomly chosen graded
artinian type $2$ quotient of the trivariate power series algebra over
a field.

\begin{bfhpg}[Compressedness of generic artinian algebras of type 2]
  \label{ss:1}
  Let $\k$ be a field and $e \ge 2$ an integer. Fr\"oberg and Laksov
  \seccite[7]{RFrDLk84} prove that given a polynomial $P(\chi)$ that
  satisfies certain numerical conditions involving $e$---in
  \rmkcite[4.2]{KSV-18} they are referred to as ``legal socle
  polynomials''---there is a non-empty Zariski open set in affine
  space $\k^d$, where $d$ depends on $e$ and the coefficients of
  $P(\chi)$, whose points are in one-to-one correspondence with
  homogeneous ideals in $\poly[\k]{x_1,\ldots,x_e}$, equivalently in
  $\pows[\k]{x_1,\ldots,x_e}$, that define compressed artinian
  $\k$-algebras with socle polynomial $P(\chi)$.  For $e \ge 3$ and
  numbers $s_1 \le s < 2s_1$ the polynomial
  $P(\chi) = \chi^{s_1} + \chi^s$ satisfies the numerical conditions
  in \prpcite[5]{RFrDLk84}: In the notation of \cite{RFrDLk84} one has
  \begin{equation*}
    \textstyle r_{s_1} \deq \binom{s_1+e-2}{e-2} - 1 - \binom{s-s_1 + e-1}{e-1}
    \dge \binom{s_1+e-1}{e-1} - 1 - \binom{s_1 + e-2}{e-1} \deq \binom{s_1+e-2}{e-2} - 1 \dge 0 \;,
  \end{equation*}
  so the discriminating value, which in \prpcite[5]{RFrDLk84} is
  called $b$, is at most $s_1$. A closer look at the proof of
  \prpcite[16]{RFrDLk84} reveals that this $b$ is actually $t$ from
  \eqref{notation1}, and the verification above amounts to the
  inequality $t \le s_1$ from \partprpref{ts1}{a}.
\end{bfhpg}

\begin{bfhpg}[Random artinian algebras of type 2]
  \label{ss:2}
  In the sense discussed above in \ref{ss:1}, a generic graded
  artinian type $2$ quotient of the power series algebra
  $Q=\pows[\k]{x,y,z}$ is compressed. It is even simpler to see that a
  generic, in the same sense, artinian Gorenstein quotient of $Q$ is
  compressed; see also Boij and Laksov \thmcite[3.4]{MBjDLk94} and
  \thmcite[4.1(d)]{KSV-18}.  Thus, with $\mfq = (x,y,z)$, if one
  generates random $\mfq$-primary Gorenstein ideals
  $I_1 \subseteq \mfq^2$ and $I_2 \subseteq \mfq^2$, then one expects
  the rings $Q/I_1$, $Q/I_2$, and $Q/(I_1\cap I_2)$ to be
  compressed. In particular, the $h$-vector of $Q/(I_1\cap I_2)$
  should be determined by the socle degrees of $Q/I_1$ and $Q/I_2$,
  and since the ideals are chosen randomly, one expects $I_1\cap I_2$
  to be minimally generated by the least possible number of elements,
  given the $h$-vector. \thmref{main} should, therefore, determine the
  class of $Q/(I_1\cap I_2)$ based on the socle degrees of $Q/I_1$ and
  $Q/I_2$, and indeed it does; we explain how in \ref{ss:3}. In fact,
  the statement of \thmref{main} was informed by the outcomes of such
  experiments conducted with \emph{Macaulay2} \cite{M2}:

  {\centering \small
    \begin{longtable}[t]{cccc|c|c|c}
      \caption{\small Let $\k$ be a field and $\mfq$ the maximal ideal
        of the local ring $Q=\pows[\k]{x,y,z}$. For fixed integers
        $2 \le s_1 \le s < 2s_1$, cf.~\eqref{ts}, and various choices
        of $\k$ we generated random $\mfq$-primary Gorenstein ideals
        $I_1 \subseteq \mfq^2$ and $I_2 \subseteq \mfq^2$ with
        quotients $Q/I_1$ and $Q/I_2$ of socle degrees $s_1$ and
        $s$. Using \cite{LWCOVl14a} we classified the rings
        $Q/(I_1\cap I_2)$ and recorded the generic, i.e.\ prevalent
        class. If $Q/(I_1\cap I_2)$ was not of the generic class, we
        still recorded it if the rings $Q/I_1$, $Q/I_2$, and
        $Q/(I_1\cap I_2)$ were all compressed, cf.~\thmref{main}. Here
        we reproduce the results for $s \le 10$.}\\
      \label{table1}
      $s_1$ & $s$ & $h$-vector & $t$ & Generic class& $m$ & Other compressed classes\\
      \hline 2 & 2 & $(1,3,2)$ & 2 & $\clH{3,2}$ & 4 & %
      \textsl{Not possible, see} \partthmref[]{main}{g} \\
      & & & & & 5 & $\clB$ \\
      \hline 2 & 3 & $(1,3,4,1)$ & 2 & $\clB$ & 5 & %
      \textsl{Not possible, see} \partthmref[]{main}{c} \\
      & & & & & 6 & $\clG{3}$ \\
      \rowcolor{Gray} 3 & 3 & $(1,3,6,2)$ & 3 & $\clH{0,0}$ & 8 & %
      \textsl{Not possible, see} \partthmref[]{main}{a} \\
      \hline 3 & 4 & $(1,3,6,4,1)$ & 3 & $\clG{3}$ & 6 & %
      \textsl{Not possible, see} \partthmref[]{main}{e} \\
      & & & & & 7 & $\clG{4}$ \\
      \rowcolor{Gray} 4 & 4 & $(1,3,6,6,2)$ & 3 & $\clH{0,0}$ & 5 & %
      \textsl{Not possible, see} \partthmref[]{main}{g} \\
      \rowcolor{Gray}
      & & & & & 6 & $\clH{0,0},\:\clG{1}$ \\
      \rowcolor{Gray}
      & & & & & 7 & $\clG{1},\:\clG{2},\:\clH{0,2}$ \\
      \rowcolor{Gray}
      & & & & & 8 & $\clH{0,0}$ \\
      \hline 3 & 5 & $(1,3,6,7,3,1)$ & 3 & $\clG{3}$ & 6 & %
      \textsl{Not possible, see} \partpgref{ss:3}{c} \\
      & & & & & 7 & $\clG{4}$ \\
      & & & & & 8 & $\clG{5}$ \\
      \rowcolor{Gray} 4 & 5 & $(1,3,6,9,4,1)$ & 3 & $\clG{1}$ & 9 & %
      \textsl{Not possible, see} \prpref[]{fam-odd} \\
      5 & 5 & $(1,3,6,10,6,2)$ & 4 & $\clH{0,0}$ & 9 & %
      \textsl{Not possible, see} \partthmref[]{main}{a} \\
      & & & & & 10 & $\clH{0,0}$ \\
      \hline 4 & 6 & $(1,3,6,10,7,3,1)$ & 4 & $\clG{5}$ & 8 & %
      \textsl{Not possible, see} \partthmref[]{main}{e} \\
      & & & & & 9 & $\clG{6}$ \\
      \rowcolor{Gray}
      5 & 6  & $(1,3,6,10,9,4,1)$  & 4 & $\clG{1}$ & 6 & $\clH{0,0}$ \\
      \rowcolor{Gray}
      & & & & & 7 & $\clG{1},\:\clG{2}$ \\
      \rowcolor{Gray}
      & & & & & 8 & $\clH{0,0},\:\clG{2},\:\clG{3}$ \\
      \rowcolor{Gray}
      & & & & & 9 & $\clG{1},\:\clG{3},\:\clG{4}$ \\
      \rowcolor{Gray}
      & & & & & 10 & $\clG{2}$ \\
      6 & 6 & $(1,3,6,10,12,6,2)$ & 4 & $\clH{0,0}$ & 9 & %
      \textsl{Not possible, see} \partthmref[]{main}{g} \\
      & & & & & 10 & $\clH{0,0},\:\clG{1}$ \\
      & & & & & 11 & $\clH{0,0}$ \\
      \hline 4 & 7 & $(1,3,6,10,11,6,3,1)$ & 4 & $\clG{4}$ & 7 & %
      \textsl{Not possible, see} \partpgref{ss:3}{c} \\
      & & & & & 8 & $\clG{5}$ \\
      & & & & & 9 & $\clG{6}$ \\
      & & & & & 10 & $\clG{7}$ \\
      \rowcolor{Gray} 5 & 7 & $(1,3,6,10,13,7,3,1)$ & 4 & $\clG{2}$ &
      10 & %
      \textsl{Not possible, see} \partpgref{ss:3}{c} \\
      \rowcolor{Gray}
      & & & & & 11 & $\clG{3}$ \\
      6 & 7 & $(1,\ldots,15,9,4,1)$ & 5 & $\clH{0,0}$ & 12 & %
      \textsl{Not possible, see} \partthmref[]{main}{a} \\
      \rowcolor{Gray} 7 & 7 & $(1,\ldots,15,12,6,2)$ & 5 & $\clH{0,0}$
      & 9 & %
      \textsl{Not possible, see} \partthmref[]{main}{a} \\
      \rowcolor{Gray}
      & & & & & 10 & $\clH{0,0}$ \\
      \rowcolor{Gray}
      & & & & & 11 & $\clH{0,0}$ \\
      \hline 5 & 8 & $(1,\ldots,15,11,6,3,1)$ & 5 & $\clG{7}$ & 10 & %
      \textsl{Not possible, see} \partthmref[]{main}{e} \\
      & & & & & 11 & $\clG{8}$ \\
      \rowcolor{Gray}
      6 & 8 & $(1,\ldots,15,13,7,3,1)$ & 5 & $\clG{3}$ & 8 & $\clH{0,0},\:\clG{2}$ \\
      \rowcolor{Gray}
      & & & & & 9 & $\clG{1},\:\clG{3},\:\clG{4}$ \\
      \rowcolor{Gray}
      & & & & & 10 & $\clG{2},\:\clG{4},\:\clG{5}$ \\
      \rowcolor{Gray}
      & & & & & 11 & $\clG{3},\:\clG{5},\:\clG{6}$ \\
      7 & 8 & $(1,\ldots,15,16,9,4,1)$  & 5 & $\clH{0,0}$ & 9 &  \\
      & & & & & 10 & $\clH{0,0},\:\clG{1}$ \\
      & & & & & 11 & $\clH{0,0},\:\clG{1},\:\clG{2}$ \\
      & & & & & 12 & $\clH{0,0},\:\clG{3}$ \\
      \rowcolor{Gray} 8 & 8 & $(1,\ldots,15,20,12,6,2)$ & 5 &
      $\clH{0,0}$ & 14 & \textsl{Not possible, see} \partthmref[]{main}{g} \\
      \hline 5 & 9 & $(1,\ldots,15,16,10,6,3,1)$ & 5 & $\clG{5}$ & 8
      & %
      \textsl{Not possible, see} \partpgref{ss:3}{c} \\
      & & & & & 9  & $\clG{6}$ \\
      & & & & & 10 & $\clG{7}$ \\
      & & & & & 11 & $\clG{8}$ \\
      \rowcolor{Gray} 6 & 9 & $(1,\ldots,15,18,11,6,3,1)$ & 5 &
      $\clG{3}$ & $11$ & %
      \textsl{Not possible, see} \partpgref{ss:3}{c} \\
      \rowcolor{Gray}
      & & & & & 12 & $\clG{4}$ \\
      7 & 9 & $(1,\ldots,21,13,7,3,1)$ & 6 & $\clH{0,0}$ & 15 & %
      \textsl{Not possible, see} \partthmref[]{main}{a} \\
      \rowcolor{Gray} 8 & 9 & $(1,\ldots,21,16,9,4,1)$ & 6 &
      $\clH{0,0}$ & 12 & %
      \textsl{Not possible, see} \partthmref[]{main}{a} \\
      \rowcolor{Gray}
      & & & & & 13 & $\clH{0,0}$ \\
      \rowcolor{Gray}
      & & & & & 14 & $\clH{0,0}$ \\
      9 & 9 & $(1,\ldots,21,20,12,6,2)$ & 6 & $\clH{0,0}$ & 8 & %
      \textsl{Not possible, see} \partthmref[]{main}{a} \\
      & & & & & 9  & $\clH{0,0}$ \\
      & & & & & 10 & $\clH{0,0}$ \\
      & & & & & 11 & $\clH{0,0}$ \\
      & & & & & 12 & $\clH{0,0}$ \\
      \hline 6 & 10 & $(1,\ldots,21,16,10,6,3,1)$ & 6 & $\clG{9}$ & 12
      & %
      \textsl{Not possible, see} \partthmref[]{main}{e} \\
      & & & & & 13 & $\clG{10}$ \\
      \rowcolor{Gray}
      7 & 10 & $(1,\ldots,21,18,11,6,3,1)$ & 6 & $\clG{5}$ & 10 & $\clG{2},\:\clG{4}$ \\
      \rowcolor{Gray}
      & & & & & 11 & $\clG{3},\:\clG{5},\:\clG{6}$ \\
      \rowcolor{Gray}
      & & & & & 12 & $\clG{4},\:\clG{6},\:\clG{7}$ \\
      \rowcolor{Gray}
      & & & & & 13 & $\clG{8}$ \\
      8 & 10 & $(1,\ldots,21,21,13,7,3,1)$ & 6 & $\clH{0,0}$ & 9 &  \\
      & & & & & 10 & $\clH{0,0},\:\clG{1}$ \\
      & & & & & 11 & $\clH{0,0},\:\clG{1},\:\clG{2}$ \\
      & & & & & 12 & $\!\clH{0,0},\:\clG{1},\:\clG{2},\:\clG{3}$ \\
      & & & & & 13 & $\clH{0,0}$ \\
      \rowcolor{Gray}
      9 & 10 & $(1,\ldots,21,25,16,9,4,1)$ & 6 & $\clH{0,0}$ & 14 &  \\
      \rowcolor{Gray}
      & & & & & 15 & $\clH{0,0},\:\clG{1}$ \\
      10 & 10 & $(1,\ldots,28,20,12,6,2)$ & 7 & $\clH{0,0}$ & 16 &
      \textsl{Not possible, see} \partthmref[]{main}{d} \\
      & & & & & 17 & $\clH{0,0}$ \\
      \hline
    \end{longtable}
  } \noindent \corref[Corollaries~]{odd} and \corref[]{even} show that
  $R$ is of class $\clH{0,0}$ or $\clG{r}$ if the socle degree, $s$,
  is at least $5$. The data shows that this is the generic behavior
  for $s\ge 4$; we elaborate on this in \pgref{ss:3}.
\end{bfhpg}

\begin{bfhpg}[Socle polynomials $\chi^{s_1} + \chi^s$ that uniquely
  determine $m$]
  \label{336779}
  For $R$ as in \thmref{main} and certain socle polynomials, the class
  of $R$ as well as $m$ is uniquely
  determined. \prpref[Propositions~]{fam-odd} and \prpref[]{fam-even}
  yield some such cases---among them those where $R$ has socle
  polynomial $\chi^4 + \chi^5$ or $2\chi^8$. The data in
  Table~\ref{table1} suggests that there are more cases of this
  behavior, and below we provide ad hoc arguments on Betti tables to
  account for those that fall within the parameters of the table. We
  use the \emph{Macaulay2} convention for compact presentation of
  Betti tables.
  \begin{prt}
  \item If $R$ is as in \thmref{main} with socle polynomial $2\chi^3$,
    then $R$ is of class $\clH{0,0}$ with $m=8$. Indeed, one has
    $h_R = (1,3,6,2)$, see Table~\ref{table0}; in particular, the
    initial degree of $I$ is $3$, so $R$ is of class $\clH{0,0}$ by
    \partthmref[]{main}{a}. A direct computation, see \eqref{HS},
    yields $B_R(\chi) = 1 -8\chi^3 + 9\chi^4-2\chi^6$. As $b_R(5)=0$
    and $R$ has type $2$, one has
    $\beta_{2\,5}(R) = 0 = \beta_{3\,5}(R)$, so the Betti table of the
    minimal free resolution of $R$ over $Q$ is%
    { \scriptsize
      \begin{equation*}
        \begin{array}{c|cccc}
          \betrow{}{0}{1}{2}{3} \\
          \hline
          \betrow{0}{1}{.}{.}{.} \\
          \betrow{1}{.}{.}{.}{.} \\        
          \betrow{2}{.}{8}{9}{.} \\
          \betrow{3}{.}{.}{.}{2\rlap{ .}} \\
        \end{array}
      \end{equation*}
    }
 
  \item If $R$ is as in \thmref{main} with socle polynomial
    $\chi^7 + \chi^9$, then $R$ is of class $\clH{0,0}$ with $m=15$.
    Indeed, one has $h_R = (1,3,6,10,15,21,13,7,3,1)$, see
    Table~\ref{table0}; in particular, the initial degree of $I$ is
    $6$, so $R$ is of class $\clH{0,0}$ by \partthmref[]{main}{a}. A
    direct computation, see \eqref{HS}, yields
    $B_R(\chi) = 1 - 15\chi^6 + 16\chi^7- \chi^{10} - \chi^{12}$. As
    $b_R(11)=0$ and $R$ has type $2$, one has $\beta_{2\,11}(R) = 0$;
    this forces $\beta_{1\,10}(R) = 0$, and continuing this standard
    analysis one sees that the Betti table of the minimal free
    resolution of $R$ over $Q$ is%
    { \scriptsize
      \begin{equation*}
        \begin{array}{c|cccc}
          \betrow{}{0}{1}{2}{3} \\
          \hline
          \betrow{0}{1}{.}{.}{.} \\
          \betrow{:}{:}{:}{:}{:} \\        
          \betrow{5}{.}{15}{16}{.} \\
          \betrow{6}{.}{.}{.}{.} \\
          \betrow{7}{.}{.}{.}{1} \\
          \betrow{8}{.}{.}{.}{.} \\
          \betrow{9}{.}{.}{.}{1\rlap{ .}} 
        \end{array}
      \end{equation*}
    }
    
  \item If $R$ is as in \thmref{main} with socle polynomial
    $\chi^6 + \chi^7$, then $R$ is of class $\clH{0,0}$ with
    $m=12$. Indeed, one has $h_R = (1,3,6,10,15,9,4,1)$, see
    Table~\ref{table0}; in particular, the initial degree of $I$ is
    $5$, so $R$ is of class $\clH{0,0}$ by \partthmref[]{main}{a}. A
    direct computation, see \eqref{HS}, yields
    $B_R(\chi) = 1 -12\chi^5 + 12\chi^6 + \chi^7 - \chi^9-
    \chi^{10}$. As above one argues that for some integer
    $\upbeta \ge 0$ the Betti table of the minimal free resolution of
    $R$ over $Q$ is%
    { \scriptsize
      \begin{equation*}
        \begin{array}{c|cccc}
          \betrow{}{0}{1}{2}{3} \\
          \hline
          \betrow{0}{1}{.}{.}{.} \\
          \betrow{:}{:}{:}{:}{:} \\
          \betrow{4}{.}{12}{12 + \upbeta}{.} \\
          \betrow{5}{.}{\upbeta}{1}{.} \\
          \betrow{6}{.}{.}{.}{1} \\
          \betrow{7}{.}{.}{.}{1}
        \end{array}
      \end{equation*}
    }%
    We use Boij--S\"oderberg theory to show that $\upbeta$ is
    zero. Assume towards a contradiction that $\upbeta$ is
    positive. Performing the first two steps of the algorithm provided
    by Eisenbud and Schreyer in \seccite[1]{DEbFOS09} one gets%
    { \scriptsize
      \begin{flalign*}
        \arraycolsep=2pt \hspace{2.1pc}\left(
          \begin{array}{c|cccc}
            \betrow{}{0}{1}{2}{3} \\
            \hline
            \betrow{0}{1}{.}{.}{.} \\
            \betrow{:}{:}{:}{:}{:} \\
            \betrow{4}{.}{12}{12 + \upbeta}{.} \\
            \betrow{5}{.}{\upbeta}{1}{.} \\
            \betrow{6}{.}{.}{.}{1} \\
            \betrow{7}{.}{.}{.}{1}
          \end{array}
        \right) - \frac{1}{5} \left(
          \begin{array}{c|cccc}
            \betrow{}{0}{1}{2}{3} \\
            \hline
            \betrow{0}{2}{.}{.}{.} \\
            \betrow{:}{:}{:}{:}{:} \\
            \betrow{4}{.}{27}{30}{.} \\
            \betrow{5}{.}{.}{.}{.} \\
            \betrow{6}{.}{.}{.}{5} \\
            \betrow{7}{.}{.}{.}{.}
          \end{array}
        \right) - \frac{11}{40} \left(
          \begin{array}{c|cccc}
            \betrow{}{0}{1}{2}{3} \\
            \hline
            \betrow{0}{2}{.}{.}{.} \\
            \betrow{:}{:}{:}{:}{:} \\
            \betrow{4}{.}{24}{35}{.} \\
            \betrow{5}{.}{.}{.}{.} \\
            \betrow{6}{.}{.}{.}{.} \\
            \betrow{7}{.}{.}{.}{3}
          \end{array}
        \right) \deq \left(
          \begin{array}{c|cccc}
            \betrow{}{0}{1}{2}{3} \\
            \hline
            \betrow{0}{\frac{1}{20}}{.}{.}{.} \\
            \betrow{:}{:}{:}{:}{:} \\                  
            \betrow{4}{.}{.}{\upbeta - \frac{7}{8}}{.} \\
            \betrow{5}{.}{\upbeta}{1}{.} \\
            \betrow{6}{.}{.}{.}{.} \\
            \betrow{7}{.}{.}{.}{\frac{7}{40}}
          \end{array}
        \right)& \:.
      \end{flalign*}
    }%
    The degree sequence of the resulting table is invalid; a
    contradiction.
  \end{prt}
\end{bfhpg}

\begin{bfhpg}[Unexplained patterns]
  \thmref{main} does not explain all the patterns one can glean from
  Table~\ref{table1}. One example is the absence of rings of class
  $\clH{0,0}$ with socle polynomial $2\chi^4$ and $m=7$. A similar
  unexplained pattern is this one: Rings of class $\clG{r}$ with socle
  polynomial $\chi^6+\chi^8$ or $\chi^7+\chi^{10}$ were observed with
  both the minimal and maximal values $r=m-8$ and $r=m-5$ determined
  by \partthmref{main}{f}, but rings with $r=m-7$ never materialized.
\end{bfhpg}

We close with an explanation of the generic behavior recorded in Table
\ref{table1}.

\begin{bfhpg}[Explanation of observed generic behavior]
  \label{ss:3}
  Let $R$ be as in \thmref{main}. Implicit in some of the arguments
  below is an assumption that the field $\k$ is infinite, or at least
  large. The generic behavior is, nevertheless, also observed when the
  coefficient field is as small as $\ZZ_2$. If $R$ is random, then $m$ is as
  small as possible given the $h$-vector of $R$, which is determined
  by the socle polynomial as $R$ is compressed.  In parts (a)--(c)
  below no further assumptions are made, so they are rigorous
  statements about the class of $R$ when $m$ is minimal.  If
  $\left\lceil \frac{s+1}{2}\right\rceil = t$ holds, then this minimal
  $m$ is determined by \thmref{resR}; for $s \le 10$ some consequences
  of (a)--(c) are recorded in the right-most column of
  Table~\ref{table1}.
  \begin{prt}
  \item If $s=2$, then $R$ is of class $\clH{3,2}$ by \eqref{ts} and
    \partthmref{main}{g}.
  
  \item If $s=3$, then \partthmref{main}{c} and \corref{mixed} show
    that $R$ is of class
    \begin{flalign*}
      & \hspace{4.5pc}
      \begin{cases}
        \clB & \text{if } \ s_1 = 2 \\
        \clH{0,0} & \text{if } \ s_1 = 3 \:.
      \end{cases} &
    \end{flalign*}
  
  \item If $s \ge 5$ is odd and one sets
    $N(s) = \frac{1}{2}(s-2+\sqrt{4s+13})$, then $R$ is of class
    \begin{flalign*}
      &\hspace{4.5pc}
      \begin{cases}
        \mathbf{G}\left(\textstyle\frac{1}{2}(s + 3 - a(a+1))\right) & \text{if } \ s_1 < N(s) \\
        \clH{0,0} & \text{if } \ N(s) \le s_1
      \end{cases}
      \qtext{with} \textstyle a = s_1 - \frac{s-1}{2} \:.&
    \end{flalign*}
    Indeed, if $N(s) \le s_1$ holds, then $R$ is of class $\clH{0,0}$,
    and if $s_1 < N(s)$ holds then $R$ is of class $\clG{r}$ by
    \corref{odd}. In the second case we compute $r$ as follows: The
    equality $\frac{s+1}{2} = t$ holds, see e.g.\ \pgref{golod}, so
    $a=s_1 - \frac{s-1}{2}$ holds by \partthmref{main}{b} and
    \thmref{resR} yields $m = t+1-f_0+f_1$. \partprpref{r}{c} and
    \eqref{f} now yield
    \begin{flalign*}
      & \hspace{3pc} \textstyle r \deq m-f_1 \deq t+1-f_0 \deq
      \frac{s+3}{2}-\frac{1}{2}a(a+1) \deq \frac{1}{2}(s+3 - a(a+1))
      \:.
    \end{flalign*}
  
  \item If $s \ge 4$ is even and one sets
    $N(s) = \frac{s}{2} -1 + \sqrt{s+4}$, then $R$ is of class
    \begin{flalign*}
      & \hspace{4.5pc}
      \begin{cases}
        \clG{s-1} & \text{if } \ \frac{s}{2} + 1 = s_1\\
        \textstyle\mathbf{G}\left(s + 3 - a(a+2)\right) & \text{if } \ \frac{s}{2} +1 < s_1 < N(s) \\
        \clH{0,0} & \text{if } \ N(s) \le s_1 \:.
      \end{cases}
      \qtext{with} \textstyle a = s_1 - \frac{s}{2} \:.&
    \end{flalign*}
    Indeed, notice first that one has $\frac{s}{2}+1 \le t \le s_1$ by
    \eqref{compressed-2} and \eqref{ts}. If equalities hold, then
    \eqref{a} yields $a=1$, and it follows from \partthmref{main}{e},
    \thmref{resR}, and \eqref{f} that $R$ is of class $\clG{r}$ with
    \begin{equation*}
      r \deq m-3 \deq 2t+1-f_0-3\deq s-1 \:.
    \end{equation*}
    If $N_2(s) \le s_1$ holds, then $R$ is of class $\clH{0,0}$ by
    \corref{even}. Assuming now that $s_1 < N_2(s)$ holds, \prpref{ta}
    yields $\frac{s}{2} + 1 = t$, so \partthmref{main}{f} gives
    $a=s_1-\frac{s}{2}$, and in view of \eqref{f} one has
    \begin{equation*}
      \tag*{$\ \qquad(\#$)}
      2t+1-f_1  \textstyle\deq s+3 - a(a+2) \:.
    \end{equation*}
    From this equality it is straightforward to verify that one has
    \begin{equation*}
      \tag*{$\ \qquad(\diamond)$}
      2t+1-f_1 \dle 0 \:\iff\: N(s) \dle s_1 \:.
    \end{equation*}
    By \corref{even} the ring $R$ is of class $\clH{0,0}$ or
    $\clG{r}$. \thmref{resR} yields
    \begin{equation*}
      \tag*{$\ \qquad(\ast)$}
      m \deq \max\set{2t+1-f_0, f_1 - f_0} \:.
    \end{equation*}
    From \partprpref{r}{d} one now gets a bound on $r$,
    \begin{equation}
      \tag*{$\ \qquad(\dagger)$}
      r \dle m-f_1 +f_0 \deq \max\set{2t + 1 - f_1, 0} \:.
    \end{equation}
    Thus, if $N(s) \le s_1$ holds, then it follows from $(\diamond)$
    and $(\dagger)$ that $R$ is of class $\clH{0,0}$. Further it is
    straightforward to verify the inequality $N(s) < N_2(s)$,
    cf.~\corref{even}.

    Finally, for $\frac{s}{2} + 1 < s_1 < N(s)$ one has
    $2t + 1 - f_1 > 0$ by $(\diamond)$ and, therefore, $\upbeta = 0$
    and $m = 2t+1-f_0$; see \thmref{resR}. Per $(\#)$ the upper bound
    on $r$ from $(\dagger)$ is $s+3 - a(a+2) = 2t + 1 - f_1 > 0$.  To
    see that, generically, this bound is achieved, and $R$ hence of
    class $\clG{s+3 - a(a+2)}$, we reason along the lines of the proof
    of \partprpref{r}{e}: Recall that, in the notation from that
    proof, the upper bound on $r$, see $(\dagger)$, is $\rnk{\psi_2}$.
    By assumption $s_1 > \frac{s}{2} + 1$ holds, whence one has
    $a \ge 2$ and, therefore, $f_0 \ge 3$ by \eqref{a} and
    \eqref{f}. By \thmref{resR} the ideal $I$ is minimally generated
    by $m \le 2t-2$ elements, say, $x_1,\ldots, x_{m}$ of degree $t$,
    and without loss of generality one can assume that they are the
    first $m$ of the $2t+1$ minimal generators $x_i$ of $I_2$; see
    \prpref{resR2}. As in the proof of \partprpref[]{r}{e} one can
    assume that $\partial_2^{F''}$ is given by a $(2t+1)\times (2t+1)$
    skew-symmetric matrix with linear entries, and that the only
    nonzero products of elements in $\sfA_1''$ and and $\sfA_2''$ are
    $\sfe_i\sff_i$ for $1 \le i \le 2t+1$.  As the entries in the
    skew-symmetric matrix are linear, every relation between the
    elements $x_1,\ldots,x_m$ is a combination of the minimal
    relations between the generators $x_i$ of $I_2$. Generically, the
    skew-symmetric matrix has nonzero entries everywhere off the
    diagonal, so the $i^\mathrm{th}$ minimal relation involves all of
    the $2t+1$ generators of $I_2$ save one, namely $x_i$. It follows
    that writing a relation that involves at most $2t-2$ of the
    generators of $I_2$ in terms of those that involve $2t$ of them
    will, generically, require all $2t+1$ of those minimal
    relations. Thus, the homomorphism $F_2 \to F_2''$ from $(\ast)$ in
    the proof of \partprpref[]{r}{e} maps the basis elements in
    internal degree $t+1$ to random $\k$-linear combinations of the
    elements of the basis for $F_2''$. While $\psi_1(\sfA_1)$ is
    spanned by $\sfe_1,\ldots,\sfe_m$, it follows that the image
    $\psi_2(\sfA_2)$ is a random subspace of $\sfA_2''$. As $\sfA''$
    is a Poincar\'e duality algebra, each element of
    $\psi_2(\sfA_2)$ has nonzero products with all of
    the basis vectors $\sfe_1,\ldots,\sfe_m$. As one has
    \begin{equation*}
      \rnk{\psi_2} \deq m -f_1 + f_0 \deq m - \textstyle
      \frac{3}{2}a(a+1) \:<\: m
    \end{equation*}
    it follows that the rank $\tilde{r}$ of the map $\tilde\delta$
    from the proof of \prpref[]{r} equals $\rnk{\psi_2}$, and
    $\tilde{r}$ is a lower bound for $r$, so one has
    $r = \rnk{\psi_2}$.
  \end{prt}
\end{bfhpg}

\def\soft#1{\leavevmode\setbox0=\hbox{h}\dimen7=\ht0\advance \dimen7
  by-1ex\relax\if t#1\relax\rlap{\raise.6\dimen7
    \hbox{\kern.3ex\char'47}}#1\relax\else\if T#1\relax
  \rlap{\raise.5\dimen7\hbox{\kern1.3ex\char'47}}#1\relax \else\if
  d#1\relax\rlap{\raise.5\dimen7\hbox{\kern.9ex
      \char'47}}#1\relax\else\if D#1\relax\rlap{\raise.5\dimen7
    \hbox{\kern1.4ex\char'47}}#1\relax\else\if l#1\relax
  \rlap{\raise.5\dimen7\hbox{\kern.4ex\char'47}}#1\relax \else\if
  L#1\relax\rlap{\raise.5\dimen7\hbox{\kern.7ex
      \char'47}}#1\relax\else\message{accent \string\soft \space #1
    not defined!}#1\relax\fi\fi\fi\fi\fi\fi}
\providecommand{\MR}[1]{\mbox{\href{http://www.ams.org/mathscinet-getitem?mr=#1}{#1}}}
\renewcommand{\MR}[1]{\mbox{\href{http://www.ams.org/mathscinet-getitem?mr=#1}{#1}}}
\providecommand{\arxiv}[2][AC]{\mbox{\href{http://arxiv.org/abs/#2}{\sf
      arXiv:#2 [math.#1]}}} \def\cprime{$'$}
\providecommand{\bysame}{\leavevmode\hbox to3em{\hrulefill}\thinspace}
\providecommand{\MR}{\relax\ifhmode\unskip\space\fi MR }
\providecommand{\MRhref}[2]{%
  \href{http://www.ams.org/mathscinet-getitem?mr=#1}{#2} }
\providecommand{\href}[2]{#2}

\end{document}